
%
\font\eightrm = cmr10 scaled 800    
\font\sixrm = cmr7 scaled 850
\font\fiverm = cmr5
\font\eightbf = cmbx10 scaled 800
\font\sixbf = cmbx7 scaled 850
\font\fivebf = cmbx5
\font\eightit = cmti10 scaled 800

\font\eighti = cmmi10 scaled 800
\font\sixi = cmmi7 scaled 850
\font\fivei = cmmi5
\font\eightsy = cmsy10 scaled 800
\font\sixsy = cmsy7 scaled 850
\font\fivesy = cmsy5
\font\eightsl = cmsl10 scaled 800
\font\eighttt = cmtt10 scaled 800

\font\fivefk = eufm5                

\font\sevenfk = eufm7

\font\tenfk = eufm10

\newfam\fkfam
   \textfont\fkfam=\tenfk \scriptfont\fkfam=\sevenfk
      \scriptscriptfont\fkfam=\fivefk

\font\fivebbm  = bbmsl10 scaled 500 
\font\sixbbm   = bbmsl10 scaled 600
\font\sevenbbm = bbmsl10 scaled 700
\font\eightbbm = bbmsl10 scaled 800
\font\tenbbm   = bbmsl10

\newfam\bbmfam
   \textfont\bbmfam=\tenbbm \scriptfont\bbmfam=\sevenbbm
      \scriptscriptfont\bbmfam=\fivebbm
\def\bm{\fam\bbmfam}
\def\eightpoint{%
   \textfont0=\eightrm \scriptfont0=\sixrm \scriptscriptfont0=\fiverm
      \def\rm{\fam0\eightrm}%
   \textfont1=\eighti  \scriptfont1=\sixi  \scriptscriptfont1=\fivei
      \def\oldstyle{\fam1\eighti}%
   \textfont2=\eightsy \scriptfont2=\sixsy \scriptscriptfont2=\fivesy
   \textfont\itfam=\eightit \def\it{\fam\itfam\eightit}%
   \textfont\slfam=\eightsl \def\sl{\fam\slfam\eightsl}%
   \textfont\ttfam=\eighttt \def\tt{\fam\ttfam\eighttt}%
   \textfont\bffam=\eightbf \scriptfont\bffam=\sixbf
      \scriptscriptfont\bffam=\fivebf \def\bf{\fam\bffam\eightbf}%
   \textfont\bbmfam=\eightbbm \scriptfont\bbmfam=\sixbbm
      \scriptscriptfont\bbmfam=\fivebbm
   \rm}
\skewchar\eighti='177\skewchar\sixi='177
\skewchar\eightsy='60\skewchar\sixsy='60
\def\small{\eightpoint\baselineskip=9.6pt}%
\parskip=1ex minus .3ex \parindent=0pt 
\hoffset=0pt \voffset=0pt 
\hsize=149mm \vsize=245mm 
\def\ifundef#1{\expandafter\ifx\csname#1\endcsname\relax}
\newif\ifpdftex \ifundef{pdfoutput}\pdftexfalse\else\pdftextrue\fi
\newdimen\paperwd \newdimen\paperht
\paperwd=210mm \paperht=297mm       
\ifpdftex
   \pdfoutput=1\pdfcompresslevel=9\pdfadjustspacing=1
   \pdfpagewidth=\the\paperwd
   \pdfpageheight=\the\paperht
\else
   
\fi
\def\today{\number\day\space\ifcase\month\or
   January\or February\or March\or April\or May\or June\or
   July\or August\or September\or October\or November\or December\fi
   \space\number\year}
\catcode`@=11
\def\usefont[#1]{\font\tmp@fnt = #1\relax\tmp@fnt}
\newif\iflinenums  \linenumstrue \message{linenumstrue}
\newcount\lineno \lineno=1
\def\par{\ifvmode\relax\else\endgraf\global\advance\lineno by \prevgraf\fi}
\def\showlineno{\iflinenums%
   \llap{{\usefont[cmti5]\the\lineno\enspace}\hglue\parindent}\fi\ignorespaces}
\everypar={\showlineno}
\everydisplay={\global\advance\lineno by -2}
\newif\ifshowtag \showtagtrue \message{showtagtrue}
\def\htag[#1]#2{%
   \ifshowtag\smash{\lower 6pt\rlap{\usefont[cmr5]label (#1)}}%
      \immediate\write99{htag l.\the\inputlineno: #1 -> #2}\fi
   \ifpdftex\pdfdest name {#1} xyz\relax#2\else
   #2\fi}
\def\href[#1]#2{\leavevmode\ifpdftex
   \pdfstartlink attr {/Border [0 0 0]} goto name {#1}\relax#2\pdfendlink
   \else#2\fi
   \ifshowtag\smash{\lower 6pt\rlap{\usefont[cmr5][#1]}}%
      \immediate\write99{href l.\the\inputlineno: #1 -> #2}\fi}
\newtoks\title  \title={}
\newtoks\author \author={(Udo Hertrich-Jeromin, \today)}
\newcount\h@noA \h@noA=0 \newcount\h@noB \h@noB=0 \newcount\h@noC \h@noC=0
\def\headno#1{\ifcase#1\relax\or\the\h@noA\or\the\h@noB\or\the\h@noC\fi}
\outer\def\h#1 #2.{%
   \ifcase#1 \relax
   \or \global\advance\h@noA by 1 \global\h@noB=0 \global\h@noC=0
      \global\equano=0
   \or \global\advance\h@noB by 1 \global\h@noC=0
   \or \global\advance\h@noC by 1
   \fi
   \ifcase#1 
      \ifpdftex\pdfinfo{/Title (#2) /Author (\the\author)}\fi
      \centerline{\usefont[cmbx10 scaled 1440]#2}
      \centerline{\usefont[cmr7]\the\author}\vskip 1ex plus .5ex
   \or
      \goodbreak\vskip 3ex\noindent
      \usefont[cmbx10 scaled 1200]\headno1.~#2
      \vglue 1ex
   \or
      \goodbreak\vskip 2ex\noindent
      \usefont[cmbx10 scaled 1000]\headno1.\headno2.~\ignorespaces#2
      \vglue 1ex
   \or
      \goodbreak\vskip 2ex\noindent
      \underbar{\bf\headno1.\headno2.\headno3.~\ignorespaces#2.}%
      \enspace\ignorespaces
   \or
      \errmessage{unknown header level}
   \or
      \errmessage{unknown header level}
   \or
      \errmessage{unknown header level}
   \or
      \errmessage{unknown header level}
   \or
      \errmessage{unknown header level}
   \or
      \errmessage{unknown header level}
   \or
      \errmessage{unknown header level}
   \or
      \goodbreak\vskip 3ex\noindent
      \usefont[cmbx10 scaled 1200]#2
      \vglue 1ex
   \else\errmessage{unknown header level}
   \fi\rm}
\let\plaineqno\eqno
\newcount\equano \equano=0
\newdimen\eqn@skip
\newbox\eqn@box
\def\eqn@no{(\headno1.\the\equano)}
\def\eqno#1{%
   \global\advance\equano by 1
   \ifinner
      \htag[eqn.#1]{\rm\eqn@no}%
   \else
      \plaineqno{\htag[eqn.#1]{\rm\eqn@no)}}%
   \fi
   \ifundef{EQN@#1}%
      \expandafter\xdef\csname EQN@#1\endcsname{\eqn@no}%
   \else
      \errmessage{duplicated equation label}%
   \fi}
\def\reqn#1{\href[eqn.#1]{{\csname EQN@#1\endcsname}}}
\def\eqn@lineno{%
   {\advance\lineno by \prevgraf\usefont[cmti5]\the\lineno\enspace}%
   \global\advance\lineno by -2}
\def\eqn@extract#1\eqno#2\eqno#3\relax{\setbox\eqn@box=\hbox{$#1$}%
   \eqn@skip=.5\displaywidth\relax\advance\eqn@skip by -.5\wd\eqn@box%
   \def\eqn@label{#2}}
\def\equadisplay#1$${%
   \eqn@extract#1\eqno\eqno\relax
   \iflinenums
      \llap{\eqn@lineno\hglue\displayindent\hglue\eqn@skip}%
   \fi
   \hbox{\unhbox\eqn@box}%
   \rlap{\hglue\eqn@skip\llap{\ifx\eqn@label\empty\else\eqno{\eqn@label}\fi}}%
   $$}
\everydisplay={\equadisplay}
\def\item@exec(#1)(#2,#3,#4)#5{{\parindent=#2\relax\par\leavevmode}%
   \def\argA{#1}\relax\def\argB{#2}\relax%
   \ifx\argA\argB\hangindent=#2\else\hangindent=#3\fi\hangafter=1%
   \llap{#5\enspace}\ignorespaces}%
\def\item@optA[#1]#2{\item@exec(#1)(#1,,){#2}}%
\def\item@optB#1{\item@optA[20pt]{#1}}%
\def\item@opt{%
   \ifx\item@char[\let\item@cmd\item@optA\relax\else%
   \let\item@cmd\item@optB\fi\relax\item@cmd}%
\def\item{\futurelet\item@char\item@opt}

\newcount\foot@no \foot@no=0
\def\footnote{\global\advance\foot@no by 1
   \edef\@sf{\spacefactor\the\spacefactor}%
   \unskip{\raise.4\baselineskip\hbox{
      \usefont[cmr7]\number\foot@no)}}\@sf
   \insert\footins\bgroup\everydisplay={}\everypar={}\let\par\endgraf
      \parindent=16pt \parskip=0pt \leftskip=0pt \rightskip=0pt
      \splittopskip=10pt plus 1pt minus 1pt \floatingpenalty=20000
      \ifundef{small}\else
         \abovedisplayskip=6pt \abovedisplayshortskip=-4pt
         \belowdisplayskip=4pt \belowdisplayshortskip=4pt
         \small\fi
      \smallskip
      \textindent{\number\foot@no)}%
      \bgroup\aftergroup\@foot\let\next}
\def\@foot{\egroup}
\skip\footins=12pt plus 2pt minus 4pt
\dimen\footins=.5\vsize
\newif\ifshowcomment \showcommentfalse
\long\def\comment#1\endcomment{%
   \ifshowcomment{\linenumsfalse\everypar={}\let\par\endgraf
   \vskip 2ex\hrule\item[0pt]{$\hookrightarrow$}#1\vskip 1ex\hrule}\fi}
\catcode`@=12
\ifundef{bm}\let\bm\bf\fi
\ifundef{fk}\fi

\def\N{{\bm N}}
\def\Z{{\bm Z}}

\def\R{{\bm R}}

\def\H{{\bm H}}
\def\HP{{\bm HP}}
\def\RP{{\bm R}\!{\sl P}}

\def\Im{\mathop{\rm Im}}
\def\End{{\sl End}}
\def\Hom{{\sl Hom}}
\def\PGl{{\bm P}{\sl Gl}}
\def\O{{\sl O}}
\def\Mob{{\sl M\!\ddot ob}}

\def\span{\mathop{\rm span}}

\def\mink#1{\left({#1\atop1}\,{-#1^2\atop-#1\hfill}\right)}
\catcode`?=11
\newif\iffile \filefalse
\newread\check?file
\def\checkfile#1{\immediate\openin\check?file=#1\relax%
   \ifeof\check?file\message{file #1 does not exist.}\filefalse%
   \else\filetrue\fi\immediate\closein\check?file}
\newcount\?llx \newcount\?lly \newcount\?urx \newcount\?ury \newcount\?rwi
\newcount\?PDFxsize \newcount\?PDFysize
\newcount\?PDFscale \newcount\?PDFres \?PDFres=288
\newwrite\?PStoPDF \immediate\openout\?PStoPDF=\jobname.ps2pdf
\def\make?PStoPDF(#1x#2+#3-#4@#5)#6{%
   \immediate\write\?PStoPDF{echo "#5 #5 scale" > #6.pstopdf.tmp}%
   \immediate\write\?PStoPDF{echo ".01 .01 scale" >> #6.pstopdf.tmp}%
   \immediate\write\?PStoPDF{echo "#3 neg #4 neg translate" >> #6.pstopdf.tmp}%
   \immediate\write\?PStoPDF{cat #6.ps >> #6.pstopdf.tmp}%
   \immediate\write\?PStoPDF{gs -sDEVICE=pdfwrite -sOutputFile=#6.pdf
      -g#1x#2 -r\the\?PDFres\space -dBATCH -dNOPAUSE #6.pstopdf.tmp -c quit}%
   \immediate\write\?PStoPDF{gs -sDEVICE=ppmraw -sOutputFile=#6.ppm
      -g#1x#2 -r\the\?PDFres\space -dBATCH -dNOPAUSE #6.pstopdf.tmp -c quit}%
   \immediate\write\?PStoPDF{ppmquant 256 #6.ppm
      | ppmtogif -transparent \#fff
      > #6.gif}%
   \immediate\write\?PStoPDF{rm -f #6.pstopdf.tmp #6.ppm #6A.ppm}}
\def\psfile(#1x#2+#3-#4@#5)#6{%
   \unskip\vbox to #2pt{\vss\hbox to #1pt{\?llx=#3\relax\?lly=#4\relax%
   \?urx=#1\multiply\?urx by 100 \divide\?urx by #5\advance\?urx by \?llx%
   \?ury=#2\multiply\?ury by 100 \divide\?ury by #5\advance\?ury by \?lly%
   \?rwi=#1\multiply\?rwi by 10%
   \ifpdftex
      \checkfile{#6.pdf}%
      \iffile\pdfximage width #1pt height #2pt {#6.pdf}%
            \pdfrefximage\pdflastximage%
      \else\?PDFscale=#5%
         \?PDFxsize=\?urx \advance\?PDFxsize by -\?llx%
         \?PDFysize=\?ury \advance\?PDFysize by -\?lly%
         \multiply\?PDFxsize by \?PDFscale \divide\?PDFxsize by 100%
         \multiply\?PDFxsize by \?PDFres \divide\?PDFxsize by 72%
         \multiply\?PDFysize by \?PDFscale \divide\?PDFysize by 100%
         \multiply\?PDFysize by \?PDFres \divide\?PDFysize by 72%
         \make?PStoPDF(%
            \the\?PDFxsize x\the\?PDFysize+\the\?llx-\the\?lly%
            @\the\?PDFscale){#6}%
         \message{...check \jobname.ps2pdf}
      \fi
   \else
      \includegraphics{#6.ps}%
   \fi\hfil}}}%
\catcode`?=12

\showcommentfalse
\linenumsfalse
\showtagfalse
\headline={{\usefont[cmu10 scaled 800]%
   discrete cmc surfaces \hfil\today~-~\the\time\hfil arXiv version }}
\def\pcq{P}
\def\define#1{{\sl #1\/}}
\long\def\proof#1\endproof{\vskip 2ex\noindent{\it Proof.}\enspace
   #1\hfill\llap{$\triangleleft$}\vskip 2ex}


\author={(Fran Burstall, Udo Hertrich-Jeromin, Wayne Rossman \&
   Susana Santos; \today)}
\title={Discrete surfaces of constant mean curvature}
\h0 Discrete surfaces of constant mean curvature.
\vskip 2em
\hbox to \hsize{\hfil\vbox{\hsize=0.8\hsize{\small{\bf Abstract.}\enspace
We propose a unified definition for discrete analogues of constant mean
curvature surfaces in spaces of constant curvature as a special case
of discrete special isothermic nets.
B\"acklund transformations and Lawson's correspondence are discussed.
It is shown that the definition generalizes previous definitions
and a construction for discrete cmc surfaces of revolution in
space forms is provided.
\par}}\hfil}
\vskip 1em
\hbox to \hsize{\hfil\vbox{\hsize=0.8\hsize{\small{\bf MSC 2000.}\enspace
{\it 53A10\/}, {\it 53C42\/}, 53A30, 52C26, 37K35, 37K25
\par}}\hfil}
\vskip 1em

\h1 Introduction.
Discrete surfaces of constant mean curvature (discrete ``cmc surfaces'')
have been studied in recent years from a variety of different points of
view.
Two essentially antithetic approaches, one from variational principles
and the other from integrable systems, lead to substantially different
definitions:
normally, for example, discrete soap films and bubbles,
i.e., ``discrete variationally cmc surfaces'',
are triangulated whereas the definition of
``integrable discrete cmc surfaces''
makes use of special coordinates and,
therefore, leads to ``discrete cmc {\it nets\/}'',
i.e., quadrilateral surfaces.
Even in cases where it is sensible to compare the two approaches,
such as that of a discrete catenoid in Euclidean space, it turns out
that different notions are obtained:
each approach leads to a different class of discrete surfaces that
can be viewed as analogues of the smooth catenoid.

The present paper is concerned with the integrable systems
approach to discrete cmc surfaces.

A key feature of this approach is its compatibility with the
transformation theory of the (smooth) surface class under
consideration:
for a given class of surface, not only is a similar transformation
theory sought for the discrete case but a discrete surface in the class
should be created by repeated (B\"acklund-Darboux-)transformations of the
smooth class; or, otherwise said, every $2$-dimensional subnet of a
multidimensional net created by repeated transformation of a discrete
surface in the class should itself be a discrete net of the class.
This is what has recently been coined ``multidimensional consistency'',
see the very clear and essential description of integrable discretization
in \href[ref.bosu06a]{[6]}.

A key idea in the definition of integrable discrete cmc surfaces
has been to consider them as special discrete isothermic nets:
that is, to discretize a (conformal) curvature line net on a smooth
cmc surface --- recall that smooth surfaces of constant mean curvature
(in any space form) are isothermic, i.e., allow a parametrization
by conformal curvature line parameters.
This has been the pioneering idea in \href[ref.bopi96]{[3]},
where the authors introduced the notion of
{\it discrete minimal surfaces in Euclidean space\/}
   alongside the notion of discrete isothermic surfaces
   and their Christoffel transformation\footnote{Note the
   parallel with Christoffel's original paper \href[ref.ch67]{[12]},
   where his transformation is introduced --- motivated by an
   observation about minimal surfaces.}.
Subsequently, the notion of
{\it discrete surfaces of constant mean curvature in Euclidean space\/}
   has been introduced alongside a notion of a Darboux transformation
   for discrete isothermic nets in \href[ref.jehopi99]{[16]},
   see also \href[ref.bopi99]{[4]},
and the notion of
{\it discrete horospherical surfaces in hyperbolic space\/}
   --- as an analogue of smooth cmc $1$ surfaces ---
   has been introduced alongside a notion of a Calapso transformation
   for discrete isothermic nets in \href[ref.je00]{[17]}.
In all three cases the constant mean curvature surfaces can be
characterized as isothermic surfaces with a special behaviour of
their transformations, as we will discuss below.

Note that for all three classes of surfaces $H^2+\kappa\geq0$, where $H$
is the mean curvature of the surface and $\kappa$ is the ambient curvature.
It is straightforward to use the Calapso transformation for discrete
isothermic nets to extend this family of definitions to discrete
analogues of any constant mean curvature surfaces with
$$
   H^2+\kappa \geq 0.
$$
Here, the key observation is that, for smooth constant mean curvature
surfaces in space forms, the Calapso transformation becomes a conformal
variant of the Lawson correspondence and that Bianchi permutability
can then be used to carry over the characterization of cmc surfaces
in Euclidean space to other space forms.
However, these ideas turn out to be useless for surfaces with
$$
   H^2+\kappa < 0
$$
as, for example, for minimal surfaces in hyperbolic space.

For discrete isothermic surfaces in Euclidean space, a construction of
a mean curvature function or rather a ``mean curvature sphere congruence''
was given in \href[ref.bopi99]{[4]} and shown to be constant for discrete
minimal or constant mean curvature surfaces\footnote{In the minimal case
the reverse is in fact also true as shown in \href[ref.bopi96]{[3]}.}.
Note that this mean curvature function is defined at the {\it vertices\/}
of a discrete isothermic net.
Very recently, new ideas from \href[ref.sch07]{[22]} have led to substantial
progress in this direction:
a new definition of discrete cmc surfaces in Euclidean space relies on the
requirement that a mean curvature function --- defined via Steiner's formula
on the {\it faces\/} of a discrete (isothermic) net --- be constant,
see also \href[ref.boetal07]{[8]}.
This definition is equivalent to the one via isothermic transformations,
see \href[ref.bo07talk]{[9]}.

\vskip 3ex\hbox to \hsize{\hfil\vbox{\hsize=300pt
\hbox to \hsize{\psfile(300x233+115-50@70){babo}}\vglue 2ex
\hbox to \hsize{\hfil{\small
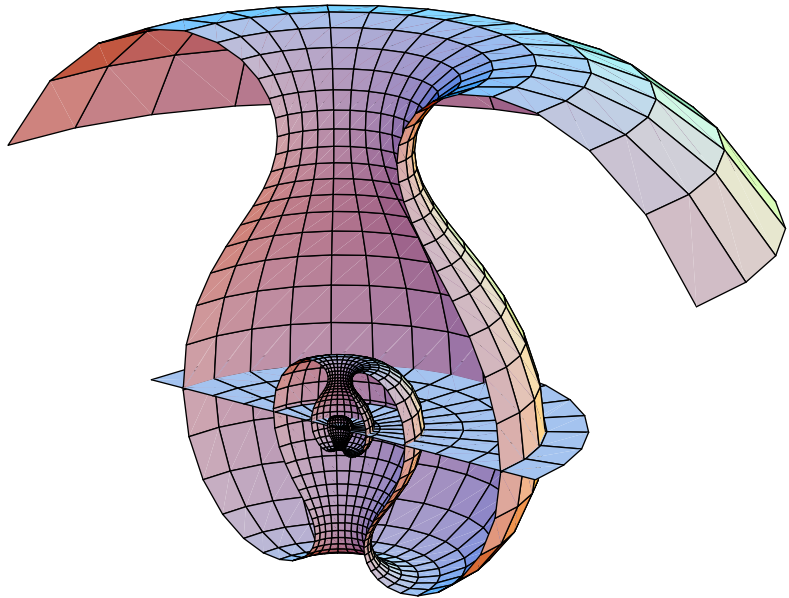
   \htag[fig.babo]{{\bf Fig.~1.}}\enspace Discrete minimal net in $H^3$
   }\hfil}}\hfil}\vskip 2ex

Our mission in the present paper will be to add another definition
of discrete cmc surfaces to the list.
However, the aim is not just to promote mathematical pluralism:
our definition provides a uniform definition of discrete cmc nets in
all space forms alike --- in particular, we also capture the previously
inaccessible case of
$$
   H^2+\kappa < 0.
$$
In fact, we define the much wider class\footnote{In contrast to the generic
terminology of ``special discrete isothermic nets'' used earlier our
``discrete special isothermic surface'' will be a technical term.}
of ``discrete special isothermic surfaces'' based on
\href[ref.buca07]{[11]} and \href[ref.sa07]{[21, Def.~2.18]},
see \href[def.special]{Def.~3.12};
these come equipped with a ``type number'' $N\in\N$
--- discrete cmc nets in space forms will be the $N=1$ case.
Hence our definition does not only provide a generalization in
allowing any ambient space form and value of the mean curvature,
but also in discussing a wider class of discrete isothermic nets
--- and we expect it to inaugurate a new direction of research
in the field.

We shall start our investigation with a short discussion of discrete
isothermic surfaces and their transformations --- not only to remind
the reader of some facts and to fix notations but also to introduce our
perspective on discrete isothermic nets via loops of flat connections,
which will be central to all that follows, see \href[thm.dic]{Lemma~2.5}.
This will set the scene for the central section of this paper:
we shall investigate the properties of polynomial loops of parallel sections,
called ``polynomial conserved quantities'', and relations to the geometry
of the underlying isothermic net.
Excluding some degenerate cases we will arrive at the notion of
``discrete special isothermic nets of type $N$'' in a natural way.
It turns out that the Darboux transformation for discrete isothermic
nets behaves nicely on these special isothermic nets, which gives
rise to a ``B\"acklund transformation'' for special isothermic nets;
in particular, we will prove a Bianchi permutability theorem that
establishes ``3D-consistency'' for special isothermic nets.
Hence our discrete special isothermic nets satisfy the two fundamental
discretization principles of the ``discrete Erlanger programme''
of \href[ref.bosu06a]{[6]}:
\item{$\bullet$} {\it Transformation group principle\/}
   --- this is built into our construction as we are working in conformal
   geometry which is the natural symmetry group for special isothermic
   (smooth) surfaces and (discrete) nets alike,
   see \href[ref.sa07]{[21, Sect.~2.2.3]};
\item{$\bullet$} {\it Consistency principle\/}
   --- which is established by our Bianchi permutability theorem for the
   B\"acklund transformation, see \href[thm.BBperm]{Thm.~4.7}.

Certain (very) special B\"acklund transforms of a special isothermic net,
its ``complementary nets'', will provide the basis for establishing the
relation of our approach with the previous approaches to discrete cmc
surfaces via their transformations as discrete isothermic surfaces
in \href[ref.bopi96]{[3]}, \href[ref.jehopi99]{[16]},
\href[ref.bopi99]{[4]} and \href[ref.je00]{[17]},
as discussed above.
Moreover, we obtain a characterization for discrete cmc surfaces
in space forms, i.e., special isothermic nets of type $1$, with
$$
   H^2+\kappa \geq 0
$$
via complementary nets --- as one may have expected; and
the lack of their existence when
$$
   H^2+\kappa < 0
$$
provides one possible explanation why the aforementioned approach
to define discrete cmc nets in space forms was doomed to failure
in this case.
In the same context we also obtain characterizations of type $2$
special isothermic nets, which discretize the classical ``special
isothermic surfaces'' of Darboux \href[ref.da99]{[13]} and Bianchi
\href[ref.bi04]{[2]}, see also \href[ref.ei23]{[14,~\S\S84--86]}.

In the final section we arrive at the main subject of this paper:
we give a definition of discrete cmc surfaces in space forms as
special isothermic nets of type $1$.
Clearly, the rich theory developed in the more general case of
special isothermic nets descends to a similarly rich theory for
discrete cmc nets in space forms --- in particular, we have the
discrete analogues of the Lawson correspondence and the B\"acklund
transformation and our discretization satisfies the two discretization
principles above for an ``integrable discretization''.
Note that we consider M\"obius geometry as the natural ambient
geometry for constant mean curvature surfaces in space forms as
these arise in ``Lawson families'' of cmc surfaces with different
ambient curvatures:
the confinement to a space form subgeometry appears as a symmetry
breaking phenomenon initiated by part of the geometric data attached
to a special isothermic net and, in particular, a discrete cmc net.

Despite the obvious merits of our definition we make a great effort
to convince the reader of its value by providing detailed analysis
of how the aforementioned previous approaches tie in with our
definition.
However, we do not conceal its problems:
we provide an example of a single discrete isothermic net which is a
cmc net in a whole family of different space forms\footnote{This net
is not spherical but it ``looks'' close to a ``wrinkled'' sphere
--- note that, in a smooth world, spheres are the only surfaces
that have constant mean curvature in different space forms.};
even though this seems to be a rather singular example we can, at
the moment, only speculate about how to obviate this anomaly.
On the positive side, we provide a method to construct discrete cmc
nets of revolution for any prescribed mean curvature $H$ and ambient
curvature $\kappa$.
In particular, we show how to explicitely construct discrete analogues of
smooth constant mean curvature surfaces that were previously unavailable:
for example, we construct (see \href[fig.babo]{Fig.~1}) the discrete
analogue of a ``hyperbolic catenoid'', that is, a minimal surface of
revolution in hyperbolic space, see \href[ref.babo93]{[1]}.

\vskip 2\parskip{\it Acknowledgements.\/}
\def\emph{\usefont[cmcsc10]}
It is our pleasure to thank our colleagues
{\emph A.~Bobenko},
{\emph T.~Hoffmann},
{\emph U.~Pinkall},
{\emph W.~Schief}
and
{\emph Y.~Suris}
for many interesting and helpful discussions about the subject.

We also gratefully acknowledge financial support for exchange visits
of the second and third authors to Japan and the UK, respectively,
from the
{\emph Daiwa Anglo-Japanese Foundation}
and the
{\emph Japanese Ministry of Education}.

The surface graphics in Figures \href[fig.babo]{1} \& \href[fig.s3torus]{2}
were produced using {\emph Mathematica}.

\h1 Discrete isothermic nets.
We consider discrete nets $f:\Z^2\supset M\to S^3$ in the (conformal)
$3$-sphere, where $$
   M = \{(m,n)\in\Z^2\,|\,m_1\leq m\leq m_2,n_1\leq n\leq n_2\}
$$
is a rectangular grid\footnote{It should be straightforward to generalize our
results to discrete nets defined on quad-graphs, making it possible to consider
discrete isothermic nets with ``umbilics'', cf.~\href[ref.ho00]{[19]}.}:

\proclaim\htag[def.din]{Def.~2.1}.
Such a net will be called a \define{discrete isothermic net} if
there is a (real) function $a$ on the edges of $M$, that is, a map
$(ij)\mapsto a_{ij}\in\R$ with $a_{ji}=a_{ij}$ for all edges $(ij)$, so
that
\item{\rm(i)}
   $a$ has equal values on opposite edges of elementary quadrilaterals
   $$
      (ijkl) = ((m,n)(m+1,n)(m+1,n+1)(m,n+1)),
   $$
   i.e., $a_{(m,n)(m+1,n)}=a_{(m,n+1),(m+1,n+1)}$ and correspondingly
   for ``vertical'' edges;
\item{\rm(ii)}
   the cross ratios\footnote{Note that the cross ratio of four points
   in $S^3$ is (up to complex conjugation) a conformal invariant.  For a
   detailed discussion see \href[ref.imdg]{[18,~Sects.~4.9,~6.5 and \S7.5.14]}.}
   $q_{ijkl}=[f_i;f_j;f_k;f_l]$ on faces factorize as
   $$
      q_{ijkl} = {a_{ij}\over a_{il}}
   $$
   into two functions of one variable.

Thus we employ the ``wide definition'' of discrete isothermic nets
\href[ref.bopi99]{[4]}, see also \href[ref.imdg]{[18,~\S5.7.2]},
which discretizes isothermic nets parametrized by curvature line coordinates
(not necessarily conformally): as all cross ratios are real, the four vertices
of any face of the net are concircular, so that a discrete isothermic net
qualifies as a discrete curvature line net or \define{discrete principal net}.
Note that the smallest domain of a discrete net where ``discrete isothermic''
imposes a condition is a $3\times3$-grid, $m_2-m_1=n_2-n_1=2$;
there the definition can be reformulated as a cross ratio $1$ condition
on four cross ratios: $$
   {q_{(m,n-1)(m+1,n-1)(m+1,n)(m,n)}\over q_{(m,n)(m+1,n)(m+1,n+1)(m,n+1)}}
   {q_{(m-1,n)(m,n)(m,n+1)(m-1,n+1)}\over q_{(m-1,n-1)(m,n-1)(m,n)(m-1,n)}}
   = 1;
$$
a cross ratio function, satisfying this condition on all $3\times3$-grids
in $M$, determines the function $a$ uniquely up to a non-zero factor.

As a mild regularity assumption, discretizing the notion of an immersed
surface parametrized by curvature lines, we will usually add the requirement
that any three of the four vertices of a face uniquely determine the circle
of the four vertices, i.e., that any three vertices are in ``general
position''.

Throughout the paper we will use the following notations:
if $g$ is a map defined on the vertices of a rectangular grid $M$,
then we let
$$
   dg_{ij} := g_j-g_i
   \quad{\rm and}\quad
   g_{ij} := {1\over2}(g_i+g_j);
\eqno differential$$
note that $(ij)\mapsto g_{ij}$ defines a function on the edges of $M$
whereas $(ij)\mapsto dg_{ij}$ defines a $1$-form,
that is, $dg_{ij}+dg_{ji}=0$.
With these notations a Leibniz rule holds:
$$
   d(g\cdot h)_{ij} = g_{ij}\cdot dh_{ij} + dg_{ij}\cdot h_{ij},
\eqno Leibniz$$
where ``$\cdot$'' denotes any product on the target space of $g$ and $h$.

\h2 The projective approach and Moutard lifts.
As we are considering nets in the {\it conformal\/} $3$-sphere it will
be helpful to consider $$
   S^3 \cong L^4/\R \subset \RP^4,
   \quad{\rm where}\quad L^4 = \{Y\in\R^{4,1}\,|\,|Y|^2=0\},
$$
as a quadric in projective $4$-space.
Recall (from \href[ref.imdg]{[18]} for example) that $2$-spheres are,
in this model, described by Minkowski $4$-spaces in $\R^{4,1}$ and
circles by Minkowski $3$-spaces or, equivalently, by their (spacelike)
orthogonal complements, and that incidence translates into a subspace
relation or as orthogonality, respectively.
For example, four points $p_n\in S^3$, $n=1,\dots,4$, generically lie
on a unique $2$-sphere $S\subset S^3$ which can be described as
$$
   S \cong \span\{P_1,\dots,P_4\} \subset \R^{4,1},
$$
where $P_n\in L^4$ with $\R P_n=p_n\in S^3\cong L^4/\R$ --- their
(complex) cross ratio\footnote{Note that $\det(\langle P_i,P_j\rangle)<0$
so that $\sqrt{\det(\langle P_i,P_j\rangle)}\in i\R$.
Using the Clifford algebra of $\R^{4,1}$, a Clifford algebra valued cross
ratio can be defined whose ``imaginary'' part encodes the $2$-sphere of
the four points \href[ref.imdg]{[18,~Sect.~6.5]}.}
is given by
$$
   q = [p_1;p_2;p_3;p_4]
     = {\langle P_1,P_2\rangle\langle P_3,P_4\rangle
      - \langle P_1,P_3\rangle\langle P_2,P_4\rangle
      + \langle P_1,P_4\rangle\langle P_2,P_3\rangle
     \pm \sqrt{\det(\langle P_i,P_j\rangle)_{i,j=1,\dots4}}
    \over2\langle P_1,P_4\rangle\langle P_2,P_3\rangle};
\eqno ProjCrossRatio$$
the cross ratio becomes real exactly when the $P_n$ become linearly
dependent, i.e., when they span a $3$-dimensional Minkowski subspace
and the four points $p_n$ are concircular.
In this case the (real) cross ratio uniquely determines the relative
position of the four points on the circle: given three of the points,
say $p_1$, $p_2$ and $p_4$, and a real cross ratio $q$ the fourth
point $p_3=\Gamma^q_{p_2,p_4}(p_1)$, where
$$
   \Gamma^q_{p,p'}(X) := X + {1\over\langle P,P'\rangle}\{
      (q-1)\,\langle X,P'\rangle P + ({1\over q}-1)\,\langle X,P\rangle P'\},
\eqno CrossRatioTrafo$$
as is easily verified from \reqn{ProjCrossRatio};
note that $\Gamma^q_{p,p'}\in \O(4,1)$ descends to a M\"obius transformation of
$S^3$ and does not depend on the choice of representatives $P,P'\in L^4$
of $p,p'\in S^3$.
Also note that
$$
   \RP^1 \cong \R\cup\{\infty\} \ni q \mapsto \Gamma^q_{p,p'}(p'') \in S^3
$$
yields a $1$-to-$1$ parametrization of the circle through three distinct
points $p,p',p''\in S^3$ in terms of the cross ratio, so that
$\Gamma^0_{p,p'}(p'')=p'$, $\Gamma^1_{p,p'}(p'')=p''$ and
$\Gamma^\infty_{p,p'}(p'')=p$.

Now suppose that $f:M\to S^3$ is an isothermic net and fix a cross ratio
factorizing function $a$.
We wish to show that there is a lift $F$ of $f$ with $$
   \langle F_i,F_j\rangle = a_{ij}
\eqno MoutardCond$$
on every edge $(ij)$ of $M$.
To this end we have to show that this scaling is compatible on any
quadrilateral: thus let $(ijkl)$ denote an elementary quadrilateral and
choose a light cone lift $F_i\in L^4$ of $f_i$; then we normalize lifts
$F_j,F_l\in L^4$ of $f_j$ and $f_l$ so that $\langle F_i,F_j\rangle=a_{ij}$
and $\langle F_i,F_l\rangle=a_{il}$.
Now we choose the lift $$
   F_k := \Gamma^{a_{ij}\over a_{il}}_{f_j,f_l}(F_i)
   = F_i + {a_{ij}-a_{il}\over\langle F_j,F_l\rangle}(F_j-F_l)
\eqno MoutardEqn$$
of $f_k$
and readily verify that $\langle F_j,F_k\rangle=a_{il}=a_{jk}$ and
$\langle F_l,F_k\rangle=a_{ij}=a_{kl}$.

Note that this lift $F$ of $f$ satisfies the discrete version
\reqn{MoutardEqn} of a Moutard equation.

Conversely, if a light cone lift $F$ of a discrete surface satisfies
a Moutard equation, $F_k-F_i\parallel F_j-F_l$ on all faces, then
it is isothermic, see \href[ref.bosu06]{[7,~Def.~9]}.
Namely, taking scalar products we learn that
$$
   \left.
   \matrix{ F_k+F_i\perp F_j-F_l \cr F_k-F_i \perp F_j+F_l \cr}
   \right\}\quad\Rightarrow\quad\left\{
   \matrix{ \langle F_j,F_k\rangle = \langle F_i,F_l\rangle \cr
            \langle F_k,F_l\rangle = \langle F_i,F_j\rangle \cr}
   \right.
$$
and hence
$$
   \langle F_i,F_k\rangle\langle F_j,F_l\rangle
   = \langle F_i,F_j-F_l\rangle\langle F_k-F_i,F_l\rangle
   = (\langle F_i,F_j\rangle-\langle F_i,F_l\rangle)^2
$$
so that \reqn{ProjCrossRatio} gives
$$
   [f_i;f_j;f_k;f_l] = {\langle F_i,F_j\rangle\over\langle F_i,F_l\rangle}.
$$

From \reqn{MoutardEqn} it is also straightforward to see that any diagonal
vertex star of a discrete isothermic net is cospherical: if $i_{(m,n)}$,
$m,n\in\{-1,0,1\}$, denote the vertices of a $3\times3$-grid then
the discrete Moutard equation \reqn{MoutardEqn} shows that the four
diagonals $F_{i_{(m,n)}}-F_{p_{(0,0)}}$, $m,n\in\{\pm1\}$, are linearly
dependent so that $$
   \dim\span\{
      F_{i_{(0,0)}},F_{i_{(1,1)}},F_{i_{(-1,1)}},F_{i_{(-1,-1)}},F_{i_{(1,-1)}}
      \} \leq 4
$$
and the five points lie on a $2$-sphere.

Assuming that the vertex star $
   \{F_{i_{(0,0)}},F_{i_{(1,0)}},F_{i_{(0,1)}},F_{i_{(-1,0)}},F_{i_{(0,-1)}}\}
$
is not cospherical the converse can also be shown, leading to two
characterizations of discrete isothermic nets,
see \href[ref.bosu06]{[7,~Sect.~3]}:

\proclaim\htag[thm.din]{Lemma 2.2}.
A discrete net $f:\Z^2\supset M\to S^3$ in the conformal $3$-sphere
is isothermic iff
\item{\rm(i)} there is a lift $F:M\to L^4$ of $f$ satisfying a discrete
   Moutard equation $F_k-F_i\parallel F_j-F_l$ on every face $(ijkl)$
   iff
\item{\rm(ii)} any diagonal vertex star is cospherical.

The sphere containing a diagonal vertex star of an isothermic net is
referred to as the \define{central sphere} of the net at the center
of the star, see \href[ref.bosu06]{[7,~Thm.~10]}.

As an example we investigate discrete surfaces of revolution:
consider the discrete net
$$
   (m,n) \mapsto f_{(m,n)} :
   = (\eta_m,\varrho_m\cos\varphi_n,\varrho_m\sin\varphi_n)
   \in \R^3 \subset \R^3\cup\{\infty\} \cong S^3
$$
where $\eta,\varrho$ and $\varphi$ are real functions of a discrete
parameter.
A straightforward cross ratio computation would reveal that
$$
   q_{(m,n)(m+1,n)(m+1,n+1)(m,n+1)}
   = -{(d\eta_{m,m+1})^2+(d\varrho_{m,m+1})^2
      \over
      4\varrho_m\varrho_{m+1}\sin^2{d\varphi_{n,n+1}\over2}}
$$
identifying the net as a discrete isothermic net.
However, we shall proceed differently to show that $f$ is an isothermic
net and to find a cross ratio factorizing function $a$ on the edges.

Consider
$$
   |(x_0,\dots,x_4)|^2 = -x_0^2 + \sum_{i=1}^4x_i^2
$$
as the quadratic form of the Minkowski scalar product of $\R^{4,1}$
and let
$$
   F^e := ({1+|f|^2\over2},f,{1-|f|^2\over2})
\eqno EuclideanLift$$
denote the \define{Euclidean lift}\footnote{The Euclidean lift into the
(flat) quadric of constant curvature (see \href[ref.imdg]{[18,~Sect.~1.4]})
$$
   {\cal Q} = \{Y\in L^4\,|\,\langle Y,Q\rangle=-1\},
   \quad{\rm where}\quad
   Q := (1,0,0,0,-1).
$$} of $f$ into the light cone $L^4\subset\R^{4,1}$.
Now observe that
$$
   F_{(m,n)} :
   = {(-1)^m\over\varrho_m}F^e_{(m,n)}
   = (-1)^m\{
     (\underbrace{\textstyle
      {1+\eta^2_m+\varrho^2_m\over2\varrho_m},
      {\eta_m\over\varrho_m},0,0,
      {1-\eta^2_m-\varrho^2_m\over2\varrho_m}
      }_{~\rlap{$=:M_m\in\R^{2,1}$}})
   + (\underbrace{
      0,0,\cos\varphi_n,\sin\varphi_n,0 
      }_{~\rlap{$=:\Phi_nC\in\R^2$}})
      \},
$$
that is, there is an orthogonal decomposition $\R^{4,1}=\R^{2,1}\oplus\R^2$
so that
$$
   F_{(m,n)} = (-1)^m(M_m+\Phi_nC) = (-1)^m\Phi_n(M_m+C),
\eqno RevolutionNet$$
where $\Phi_n$ are rotations of $\R^2$, $C\in S^1\subset\R^2$ and $M$ takes
values in the hyperbolic plane\footnote{Secretly we are using a conformal map
$\R^3\setminus\{{\rm axis}\}\to H^2\times S^1\subset\R^{2,1}\oplus\R^2=\R^{4,1}$
adapted to the rotational symmetry of the map $f$,
cf.~\href[ref.imdg]{[18,~\S1.4.16]}.}
$$
      H^2=\{Y\in\R^{2,1}\,|\,|Y|^2=-1,Y_0>0\}\subset\R^{2,1}.
$$
In particular, $M_m\perp\Phi_nC$ for all $(m,n)$.

Note that $\R^2=\span\{\Phi_nC\,|\,n\in\Z\}$ defines an elliptic sphere
pencil, hence (cf.~\href[ref.imdg]{[18,~Sect.~1.2]}) a circle
$$
   L^4\cap\R^{2,1} = L^4\cap\{\Phi_nC\,|\,n\in\Z\}^\perp,
$$
which is the axis of our discrete surface of revolution.
At the same time, it is the infinity boundary of the hyperbolic $2$-plane
$H^2$ of the meridian curve.

Clearly, $F$ satisfies the discrete Moutard equation
$$
   F_{(m+1,n+1)}-F_{(m,n)}
   = (-1)^{m+1}\{M_{m+1}+\Phi_{n+1}C+M_m+\Phi_nC\}
   = F_{(m+1,n)}-F_{(m,n+1)}
$$
and is therefore a discrete isothermic net with cross ratio
factorizing function
$$
   a_{ij} := \langle F_i,F_j\rangle
   = \cases{
      -1-\langle M_m,M_{m+1}\rangle
      = {(d\eta_{m,m+1})^2+(d\varrho_{m,m+1})^2\over2\varrho_m\varrho_{m+1}} &
      for $(ij)=((m,n)(m+1,n))$ \cr
      -1+\langle \Phi_nC,\Phi_{n+1}C\rangle
      = -2\sin^2{d\varphi_{n,n+1}\over2} &
      for $(ij)=((m,n)(m,n+1))$ \cr
      }
$$
as soon as $M_{m+1}\neq M_m$ and $\Phi_{n+1}C\neq\Phi_nC$.

\h2 Quaternions and the Calapso transformation.
The Calapso transformation, or $T$-transformation, of (discrete)
isothermic nets will be central to our investigations --- it was introduced
in \href[ref.je00]{[17]} (see also \href[ref.imdg]{[18,~\S5.7.16]})
using a quaternionic setup for M\"obius geometry.
Hence we will first briefly discuss the quaternionic approach in order to
make contact with earlier work; however, we will provide an independent
definition in the following section so that a reader unfamiliar with
previous approaches may just skip this section.

Thus, we consider $S^3\cong\Im\H\cup\{\infty\}\subset\HP^1$ and
$$
   \R^{4,1} \cong \{ X\in\End(\H^2) \,|\,
      X = \left({x\atop x_0}\,{x_\infty\atop-x}\right),
      x\in\Im\H,x_0,x_\infty\in\R\} \subset \End(\H^2)
$$
equipped with $|X|^2=-X^2=x^2+x_0x_\infty$ as the quadratic form of the
Minkowski product\footnote{This is analogous to the Vahlen matrix approach
to M\"obius geometry (see \href[ref.imdg]{[18,~Sect.~7.1]}) using the Clifford
algebra of $\R^{4,1}$.}.
In particular, we obtain an isometry
$$
   \R^3\cong\Im\H \ni x \rightarrow X=\mink x \in L^4 \subset \R^{4,1};
\eqno QEuclideanLift$$
note that, for two such ``Euclidean'' light cone lifts $X,Y\in L^4$,
$$
   -2\langle X,Y\rangle = XY+YX = -(y-x)^2 = |y-x|^2.
$$
In this setup the M\"obius group
$$
   \Mob(3)
   = \{ \left({a\atop c}\,{b\atop d}\right) \in \End(\H^2) \,|\,
   \bar ac+\bar ca=\bar bd+\bar db=0, \bar ad+\bar cb\in\R\setminus\{0\}
   \}/\R \subset \PGl(2,\H)
$$
of $S^3$ acts isometrically on $\R^{4,1}$ via
$$
   \Mob(3)\times\R^{4,1} \ni (A,X)\mapsto A\cdot X:=AXA^{-1} \in \R^{4,1}.
\eqno QuatAction$$

Now let $f:M\to S^3$ be a discrete conformal net, with cross ratio
factorizing function $a$ on the edges, and define
$$
   \tau_{ij} := {a_{ij}\over2\langle F_i,F_j\rangle}\,F_iF_j,
\eqno DefTau$$
where $F$ is any light cone lift of $f$;
assuming that $f:M\to\Im\H\subset S^3$ and using the lift \reqn{QEuclideanLift}
we find
$$
   \tau_{ij} = \left(
      {f_idf^\ast_{ij}\atop\hfill df^\ast_{ij}}
      {-f_idf^\ast_{ij}f_j\atop\hfill -df^\ast_{ij}f_j}
      \right),
      \quad{\rm where}\quad
      df^\ast_{ij} = a_{ij}(df_{ij})^{-1}
$$
is the ``derivative'' of the Christoffel transform $f^\ast$ of $f$ in $\R^3$,
see \href[ref.bopi99]{[4,~Thm.~14]} or \href[ref.imdg]{[18,~\S5.7.7]}.
Then, for $\lambda\in\R$,
$$
   (1+\lambda\tau_{ij})(1+\lambda\tau_{ji}) = 1 - \lambda a_{ij} \in \R
$$
and
$$
   (1+\lambda\tau_{ij})(1+\lambda\tau_{jk})
   = (1+\lambda\tau_{il})(1+\lambda\tau_{lk})
$$
on every elementary quadrilateral $(ijkl)$ so that
$$
   (ij) \mapsto 1+\lambda\tau_{ij}, \quad
   1+\lambda\tau_{ij}:\{j\}\times S^3\to\{i\}\times S^3,
\eqno QuatConnection$$
defines a flat $\Mob(3)$-connection\footnote{We shall make the notion
of a flat (discrete) connection precise in the following section.} on
$M\times S^3$, as long as
$$
   1-\lambda a_{ij}\neq0
   \quad\Leftrightarrow\quad
   1+\lambda\tau_{ij} \in \Mob(3)
$$
for all edges $(ij)$.
Hence there is a gauge transformation
$$
   T^\lambda:M\to\Mob(3), \quad
   T^\lambda_j = T^\lambda_i(1+\lambda\tau_{ij}),
\eqno QuadCalapso$$
which identifies the $(1+\lambda\tau)$-connection on $M\times S^3$
with the trivial connection.

The $T^\lambda$ play a key role in the transformation theory of (discrete)
isothermic nets; in particular, it turns out that every $T^\lambda f$
defines\footnote{Here $\Mob(3)$ acts on $S^3\cong\Im\H\cup\{\infty\}$
by M\"obius transformations, i.e., by fractional linear transformations.}
a discrete isothermic net:
$T^\lambda:M\to\Mob(3)$ are the \define{Calapso transformations} of $f$ and
the discrete isothermic nets $T^\lambda f$ are its \define{Calapso transforms},
see \href[ref.imdg]{[18,~\S5.7.16]}.

The connections $(ij)\mapsto1+\lambda\tau_{ij}$ lift to flat
$\O(4,1)$-connections $(ij)\mapsto\Gamma^\lambda_{ij}$ on the
(discrete) vector bundle $M\times\R^{4,1}$ via \reqn{QuatAction}
to give
$$
   X \mapsto \Gamma^\lambda_{ij}\cdot X :
   = {1\over 1-\lambda a_{ij}}(1+\lambda\tau_{ij})\,X\,(1+\lambda\tau_{ji})
   = \Gamma^{1-\lambda a_{ij}}_{f_i,f_j}(X).
\eqno QuatConnectionLift$$

\h2 The vector bundle approach.
Clearly, the flat connections $\Gamma^\lambda$ on $M\times\R^{4,1}$
in \reqn{QuatConnectionLift} can be defined without reference to
the quaternionic approach.
First we define our setup:

\proclaim\htag[def.dfc]{Def.~2.3}.
A \define{connection} on a (discrete) fibre bundle $F\to M$, where
the base $M$ is a rectangular grid as before, is a map that assigns
to each directed edge $(ij)$ in $M$ an isomorphism
$$
   \Gamma_{ij}:F_j\to F_i
   \quad\hbox{\sl so that}\quad
   \Gamma_{ij}\Gamma_{ji}=1;
$$
it will be said to be a \define{flat connection} if its holonomies
around all elementary quadilaterals $(ijkl)$ are trivial,
$$
   \Gamma_{ij}\Gamma_{jk}\Gamma_{kl}\Gamma_{li} = 1.
$$

With these notions we can now formulate the key definition:

\proclaim\htag[def.dic]{Def.~2.4}.
Let $f:M\to S^3$ be a discrete isothermic net with cross ratio factorizing
function $a$.
We say that
$$
   (\lambda,ij)\mapsto\Gamma^\lambda_{ij} :
   = \Gamma^{1-\lambda a_{ij}}_{f_i,f_j}
   \in \Hom(\{j\}\times\R^{4,1},\{i\}\times\R^{4,1})
\eqno IsothConnection$$
defines the \define{isothermic family of connections}\footnote{We have
$\Gamma^\lambda_{ij}\Gamma^\lambda_{ji}=1$ on all edges $(ij)$, so that
the $\Gamma^\lambda$ qualify as discrete (linear) connections.} of $f$,
where $\lambda\in\R$ so that $1\neq\lambda a_{ij}$ for all edges $(ij)$.

Note that the $\Gamma^\lambda$ are \define{metric} connections on
$M\times\R^{4,1}$, as all $\Gamma^\lambda_{ij}$ are isometries and
hence descend to connections on $M\times S^3$.

We already know that, if $f$ is an isothermic net, then the isothermic
family of connections \reqn{IsothConnection} is flat,
see \href[ref.imdg]{[18,~\S5.7.5]}.
Here we shall give an independent proof, not relying on the quaternionic
setup, as well as a certain converse of this fact
(cf.~\href[ref.je00]{[17,~Thm.~3.14]}):

\proclaim\htag[thm.dic]{Lemma 2.5}.
Let $f:M\to S^3$ be a regular discrete net, i.e., any three vertices of
a face are in general position, and let $a$ be a function on the edges.
Then the connection given by \reqn{IsothConnection} is flat if and
only if $f$ is isothermic with cross ratio factorizing function $a$.

\proof
First note that $\Gamma^{1-\lambda a_{ij}}_{f_i,f_j}(X)=X\bmod f_i$
when $X\perp f_i$, so that $\Gamma^{1-\lambda a_{ij}}_{f_i,f_j}$ projects
to the identity on $f_i^\perp/f_i$ and similarly for $f_j$.
Consequently,
$$
   \Gamma^{1-\lambda a_{ij}}_{f_i,f_j}\Gamma^{1-\lambda a_{jk}}_{f_j,f_k}(X)
   = X \bmod f_j \quad{\rm if}\quad X\perp f_j.
$$
The same is true for the product
$\Gamma^{1-\lambda a_{il}}_{f_i,f_l}\Gamma^{1-\lambda a_{lk}}_{f_l,f_k}$
so that, if we now assume flatness of the connection,
$$
   \Gamma^{1-\lambda a_{ij}}_{f_i,f_j}
   \Gamma^{1-\lambda a_{jk}}_{f_j,f_k}
   =
   \Gamma^{1-\lambda a_{il}}_{f_i,f_l}
   \Gamma^{1-\lambda a_{lk}}_{f_l,f_k}
   =: \Gamma^\lambda,
$$
we learn that $\Gamma^\lambda=id$ on $(f_j\oplus f_l)^\perp$
since $f_j\neq f_l$.
Moreover, $f_j$ and $f_l$ are eigendirections of $\Gamma^\lambda$ and
$$
   \Gamma^\lambda(X) = \cases{\hfill
      {1-\lambda a_{jk}\over1-\lambda a_{ij}}\,X & if $X \in f_j$, \cr\hfill
      X & if $X\perp f_j,f_l$, \cr\hfill
      {1-\lambda a_{lk}\over1-\lambda a_{il}}\,X & if $X \in f_l$. \cr
      }
$$
As $\Gamma^\lambda$ is, along with $\Gamma^{1-\lambda a_{ij}}_{f_i,f_j}$
and $\Gamma^{1-\lambda a_{jk}}_{f_j,f_k}$, an orientation preserving
orthogonal transformation we infer that, for all $\lambda$,
$$
   {1-\lambda a_{jk}\over1-\lambda a_{ij}}
   {1-\lambda a_{lk}\over1-\lambda a_{il}}
   = 1
   \quad{\rm and}\quad
   \Gamma^\lambda = \Gamma_{f_j,f_l}^{1/q(\lambda)},
   \quad q(\lambda):={1-\lambda a_{ij}\over1-\lambda a_{jk}}.
$$
Hence $a_{ij}=a_{jk}$ and $a_{il}=a_{lk}$ (in which case $\Gamma^\lambda=id$
for all $\lambda$) or $a_{ij}=a_{lk}$ and $a_{jk}=a_{il}$.

Now decompose $F_k\in f_k\setminus\{0\}$ as
$F_k=F_i+F_j+F_k^\perp\in f_i\oplus f_j\oplus(f_i\oplus f_j)^\perp$
and observe that
$$
   \Gamma^\lambda(F_k)
   = {1\over1-\lambda a_{jk}}\Gamma_{f_i,f_j}^{1-\lambda a_{ij}}(F_k)
   = {1-\lambda a_{ij}\over1-\lambda a_{jk}}\,F_i
   + {1\over1-\lambda a_{jk}}\,F_k^\perp
   + {1\over(1-\lambda a_{ij})(1-\lambda a_{jk})}\,F_j
   \longrightarrow {a_{ij}\over a_{jk}}\,F_i
$$
as $\lambda\to\infty$.
This shows that, since $f_k\neq f_i$, we cannot have $a_{ij}=a_{jk}$;
hence the function $a$ has equal values on opposite edges of an
elementary quadrilateral.
Moreover, we learn that
$$
   f_k = \Gamma_{f_j,f_l}^{q(\infty)}(f_i),
$$
showing that the four vertices $f_i$, $f_j$, $f_k$ and $f_l$ are concircular
and the edge function $a$ factorizes their cross ratio,
$$
   [f_i;f_j;f_k;f_l] = q(\infty) = {a_{ij}\over a_{jk}}.
$$
Hence $f$ is a discrete isothermic net.

Conversely, suppose that $f$ is discrete isothermic with cross ratio
factorizing function $a$; we wish to show that
$$
   \Gamma^\lambda_{ij}\Gamma^\lambda_{jk}
   = \Gamma^{1/q(\lambda)}_{f_j,f_l}
   = \Gamma^{q(\lambda)}_{f_l,f_j}
   = \Gamma^\lambda_{il}\Gamma^\lambda_{lk},
$$
where $q(\lambda)={1-\lambda a_{ij}\over1-\lambda a_{jk}}$.
As the second equation holds and the third is obtained from the first
by exchanging the roles of $j$ and $l$ it suffices to prove the first
of these equations.
Also,
$$
   \Gamma^{q(\lambda)}_{f_j,f_l}\Gamma^\lambda_{ij}\Gamma^\lambda_{jk}(X)
   = \cases{
      X & if $X\in f_j$ \cr
      X \bmod f_j & if $X\perp f_j$ \cr
      }
$$
so that flatness of the family of isothermic connections of $f$ follows as
soon as $\Gamma^{q(\lambda)}_{f_j,f_l}\Gamma^\lambda_{ij}\Gamma^\lambda_{jk}$
has another isotropic eigendirection --- we shall show that $f_k$ serves this
purpose:
consider
$$
   \RP^1 \cong \R\cup\{\infty\} \ni \lambda
   \mapsto
   \Gamma^\lambda_{ij}\Gamma^\lambda_{jk}(f_k)=\Gamma^\lambda_{ij}(f_k),
   \Gamma^{1/q(\lambda)}_{f_j,f_l}(f_k) \in S^3.
$$
Both maps parametrize the same circle in terms of a certain cross ratio,
given as a linear fractional transformation of $\lambda$;
in particular,
$$\matrix{
   \lambda = 0 \hfill&\Rightarrow&
      \Gamma^\lambda_{ij}(f_k)
      = \Gamma^1_{f_i,f_j}(f_k)
      = f_k \hfill
      = \Gamma^1_{f_j,f_l}(f_k)
      = \Gamma^{1/q(\lambda)}_{f_j,f_l}(f_k),
      \cr
   \lambda = \infty \hfill&\Rightarrow&
      \Gamma^\lambda_{ij}(f_k)
      = \Gamma^\infty_{f_i,f_j}(f_k)
      = f_i \hfill
      = \Gamma^{a_{jk}\over a_{ij}}_{f_j,f_l}(f_k)
      = \Gamma^{1/q(\lambda)}_{f_j,f_l}(f_k),
      \cr
   \lambda = {1\over a_{ij}} \hfill&\Rightarrow&
      \Gamma^\lambda_{ij}(f_k)
      = \Gamma^0_{f_i,f_j}(f_k)
      = f_j \hfill
      = \Gamma^\infty_{f_j,f_l}(f_k)
      = \Gamma^{1/q(\lambda)}_{f_j,f_l}(f_k).
      \cr
}$$
As two M\"obius transformations of a circle coincide as soon as they
coincide at three points we conclude that
$\Gamma^\lambda_{ij}(f_k) = \Gamma^{1/q(\lambda)}_{f_j,f_l}(f_k)$
for all $\lambda$.
\endproof

Note that, freeing the second part of the proof from the specific notations
of the situation, we have proved the following:

\proclaim\htag[thm.dicCor]{Lemma 2.6}.
Write $[p_1;p_2;p_3;p_4]={a\over b}$ with $a,b\in\R$ for the cross ratio
of four concircular points $p_i\in S^3$, $i=1,\dots,4$; then, for all
$\lambda\in\R$,
$$
     \Gamma^{1-a\lambda}_{p_1,p_2}\Gamma^{1-b\lambda}_{p_2,p_3}
   = \Gamma^{(1-b\lambda)/(1-a\lambda)}_{p_2,p_4}
   = \Gamma^{(1-a\lambda)/(1-b\lambda)}_{p_4,p_2}
   = \Gamma^{1-b\lambda}_{p_1,p_4}\Gamma^{1-a\lambda}_{p_4,p_3}.
$$

Thus, for a discrete isothermic net $f$, there are gauge transformations
$T^\lambda:M\to\O(4,1)$ identifying the connections $\Gamma^\lambda$
on $M\times\R^{4,1}$ with the trivial connection:

\proclaim\htag[def.calapso]{Lemma \& Def.~2.7}.
Let $f:M\to S^3$ be a discrete isothermic net with its isothermic family
of connections $\Gamma^\lambda$.
Then the gauge transformations
$$
   T^\lambda:M\to\O(4,1) \quad{\rm with}\quad
   T^\lambda_j=T^\lambda_i\Gamma^\lambda_{ij}
$$
are the \define{Calapso transformations} of $f$;
the isothermic nets $f^\lambda:=T^\lambda f$ are its
\define{Calapso transforms}.

Note that $a^\mu={a\over1-\mu a}$ is a cross ratio factorizing function
for the Calapso transform $f^\mu$ of $f$ with cross ratio factorizing
function $a$, see \href[ref.imdg]{[18,~\S5.7.16]};
the isothermic family of connections of $f^\mu$ is given by
$$
   \Gamma^{\mu,\lambda}_{ij}
   = T^\mu_i\Gamma^{\mu+\lambda}_{ij}(T^\mu_j)^{-1},
\eqno IsothConnectionT$$
which shows that the Calapso transformations of a discrete isothermic net
satisfy a $1$-parameter group property, see \href[ref.je00]{[17]} or
\href[ref.imdg]{[18,~\S5.7.30]}:
$$
   T^{\mu,\lambda}T^\mu = T^{\mu+\lambda}.
\eqno 1ParameterT$$

\h1 Polynomial conserved quantities.
The second key notion in our definition of discrete cmc nets in space
forms will be that of polynomial conserved quantities:

\proclaim\htag[def.pcq]{Def.~3.1}.
Let $f:M\to S^3$ be an isothermic net.
A \define{polynomial conserved quantity} of $f$ is a map
$$
   \R\times M \ni (\lambda,i) \mapsto \pcq_i(\lambda)
   = \sum_{k=0}^N\pcq^{(k)}_i\lambda^k \in \R^{4,1}[\lambda]
$$
so that, for every fixed $\lambda$,
$$
   T^\lambda\pcq(\lambda) \equiv const.
$$

Hence, a polynomial conserved quantity of an isothermic net can be thought
of as a polynomial family of parallel sections of the vector bundle
$M\times\R^{4,1}$ equipped with the isothermic family of connections:
$$
   T^\lambda\pcq(\lambda) \equiv const.
   \quad\Leftrightarrow\quad
   \pcq_i(\lambda) = \Gamma^\lambda_{ij}\pcq_j(\lambda)
\eqno PCQcondition$$
on all edges $(ij)$ of $M$.

\h2 Basic properties.
Clearly, as $T^\lambda$ acts linearly on $\R^{4,1}$, the polynomial
conserved quantities of a given discrete isothermic net $f$ can
be superposed:

\proclaim\htag[thm.pcqSpace]{Lemma 3.2}.
The space of polynomial conserved quantities of $f$ is a vector space.

As a consequence, we can construct new polynomial conserved quantities
from a given one by multiplying with real polynomials $p(\lambda)$,
thereby raising the degree; for example, if $\pcq$ is a polynomial
conserved quantity of $f$ then $(1+\lambda)\pcq(\lambda)$ will
be a new polynomial conserved quantity of higher degree.
Thus we will be interested in (non-vanishing) polynomial conserved
quantities of lowest possible degree.
The following lemma provides a criterion:

\proclaim\htag[thm.pcqReduction]{Lemma 3.3}.
Let $\pcq(\mu)=0$ for a polynomial conserved quantity
$\pcq(\lambda):M\to\R^{4,1}[\lambda]$ of $f$; then
$$
   \tilde\pcq(\lambda) := {1\over\lambda-\mu}\,\pcq(\lambda)
$$
is a polynomial conserved quantity of $f$ of lower degree.

Note that, if $\pcq_i(\mu)=0$ for some $i\in M$ then $\pcq(\mu)\equiv0$
on $M$ since $T^\mu\pcq(\mu)\equiv const$.

\proof
Writing $\pcq_i(\lambda)\in\R^{4,1}[\lambda]$ in terms of a basis of
$\R^{4,1}$ shows that $\mu$ is a common zero for all (real) component
polynomials, which are therefore divisible by $(\lambda-\mu)$.
Hence $\tilde\pcq_i(\lambda)$ is polynomial at any $i\in M$.

Clearly
$$
   T^\lambda\tilde\pcq(\lambda)
   = {1\over\lambda-\mu}\,T^\lambda\pcq(\lambda)
   \equiv const
$$
for any fixed $\lambda$, showing that $\tilde\pcq(\lambda)$ is a polynomial
conserved quantity of $f$.
\endproof

As a direct consequence of the previous two lemmas we learn that,
if two distinct polynomial conserved quantities $\pcq(\lambda)$ and
$\tilde\pcq(\lambda)$ of degree $N\in\N$ of an isothermic net have the
same value at some point $(\mu,i)\in\R\times M$, then there is a polynomial
conserved quantity of degree $\leq N-1$.
In particular:

\proclaim\htag[thm.pcqUniq]{Cor.~3.4}.
A non-zero polynomial conserved quantity of lowest possible degree $N$,
$$
   \pcq(\lambda) = \lambda^NZ + \dots + \lambda^0Q:M\to\R^{4,1}[\lambda],
$$
is uniquely determined by either its top or bottom coefficient $Z$ or $Q$,
respectively.

Since the $T^\lambda$ are orthogonal transformations, there is another
obvious property of a polynomial conserved quantity which will become
important later:

\proclaim\htag[thm.pcqNorm]{Lemma 3.5}.
If $\pcq(\lambda)$ is a polynomial conserved quantity of $f$,
then $|\pcq(\lambda)|^2$ depends only on $\lambda$;
in particular, $|Z|^2$ and $|Q|^2$ are constants.
If $\pcq(\lambda)=\lambda Z+Q$ is a linear conserved quantity,
then also $\langle Z,Q\rangle\equiv const$.

Finally note that the equation \reqn{PCQcondition} for a polynomial
conserved quantity depends crucially on the choice of a cross ratio
factorizing function $a$: however, if $\tilde a:=\alpha a$ is a
new cross ratio factorizing function, then
$$
   \tilde\Gamma^\lambda = \Gamma^{\alpha\lambda}
$$
by \reqn{IsothConnection};
hence $\tilde\pcq(\lambda)=\pcq(\alpha\lambda)$ is a new polynomial
conserved quantity satisfying \reqn{PCQcondition} with the new
isothermic family of connections.
Consequently:

\proclaim\htag[thm.pcqScale]{Lemma 3.6}.
If $\pcq(\lambda)$ is a polynomial conserved quantity of $f$ with respect
to $a$ as a cross ratio factorizing function then
$$
   \tilde\pcq(\lambda) := \pcq(\alpha\lambda)
$$
is a polynomial conserved quantity of $f$ with respect to $\tilde a:=\alpha a$
as a new cross ratio factorizing function.

\h2 Geometric properties.
We now turn to a more detailed analysis of \reqn{PCQcondition} and its
geometric consequences:
using \reqn{CrossRatioTrafo} the fact that a polynomial conserved quantity
$\pcq(\lambda)$ is a family of parallel sections of the isothermic family
of connections \reqn{IsothConnection} of $f$ reads
$$\matrix{
   \pcq_i(\lambda)
   &=& \Gamma^\lambda_{ij}\pcq_j(\lambda)
   &=& \pcq_j(\lambda) + {\lambda a_{ij}\over\langle F_i,F_j\rangle} \{
      {1\over1-\lambda a_{ij}}\langle\pcq_j(\lambda),F_i\rangle F_j
      - \langle\pcq_j(\lambda),F_j\rangle F_i \} \cr
}$$
on any edge $(ij)$ of $M$;
exchanging the roles of the endpoints $i$ and $j$ of the edge we obtain
a similar equation which, as $f_i\neq f_j$, yields two equations
$$\matrix{
   d\pcq_{ij}(\lambda)
   &=&\hfill {\lambda a_{ij}\over\langle F_i,F_j\rangle} \{
      \langle\pcq_j(\lambda),F_j\rangle F_i
    - \langle\pcq_i(\lambda),F_i\rangle F_j \} \cr
   &=& {\lambda a_{ij}\over1-\lambda a_{ij}}{1\over\langle F_i,F_j\rangle} \{
      \langle\pcq_i(\lambda),F_j\rangle F_i
    - \langle\pcq_j(\lambda),F_i\rangle F_j \} \cr
}$$
for some light cone lift $F$ of $f$.
Note that the second equality follows from the first by taking scalar
products with $F_i$ and $F_j$, respectively.
Hence, we also obtain the converse:

\proclaim\htag[thm.PCQcondition]{Lemma 3.7}.
$\pcq(\lambda)$ is a polynomial conserved quantity of $f$ if and only if,
for all edges $(ij)$ in $M$,
$$
   d\pcq_{ij}(\lambda)
   = {\lambda a_{ij}\over\langle F_i,F_j\rangle} \{
      \langle\pcq_j(\lambda),F_j\rangle F_i
    - \langle\pcq_i(\lambda),F_i\rangle F_j \}.
\eqno dPCQ$$

Note that, in case $F$ is a Moutard lift of $f$ satisfying \reqn{MoutardCond},
then \reqn{dPCQ} simplifies to
$$
   d\pcq_{ij}(\lambda) = \lambda\,\{p_j(\lambda)F_i-p_i(\lambda)F_j\},
   \quad{\rm where}\quad
   p(\lambda) := \langle\pcq(\lambda),F\rangle.
$$
The integrability $d^2\pcq(\lambda)=0$ of this equation then yields
$$
     (p_k(\lambda)-p_i(\lambda))(F_j-F_l)
   = (p_j(\lambda)-p_l(\lambda))(F_k-F_i),
$$
and hence $p(\lambda)$ satisfies the very same Moutard equation
\reqn{MoutardEqn} as $F$ does.

Now the key observation from \reqn{dPCQ} is that this equates a polynomial
of degree $N$ and a polynomial of degree $N+1$ with vanishing constant
coefficient.
Hence, looking at the degree $0$ and degree $N$ and $N+1$ terms
we obtain the following two corollaries:

\proclaim\htag[thm.ConstantQ]{Cor.~3.8}.
If $\pcq(\lambda)=\lambda^NZ+\dots+Q$ is a polynomial conserved quantity
of $f$, then $Q\equiv const$.

Thus a polynomial conserved quantity naturally provides an ambient quadric
${\cal Q}$ of constant curvature $\kappa=-|Q|^2$ for the isothermic net,
see \href[ref.imdg]{[18,~Sect.~1.4]}:
$$
   {\cal Q} = \{Y\in L^4\,|\,\langle Y,Q\rangle=-1\}.
\eqno SpaceForm$$

\proclaim\htag[thm.EnvelopingZ]{Cor.~3.9}.
If $\pcq(\lambda)=\lambda^NZ+\lambda^{N-1}Y+\dots+Q$ is a polynomial conserved
quantity of $f$ then:
\item{\rm(i)} $Z_i\perp f_i$ at all points $i$ in $M$;
\item{\rm(ii)} $
     Z_i+a_{ij}{\langle Y_j,F_j\rangle\over\langle F_i,F_j\rangle}\,F_i
   = Z_j+a_{ij}{\langle Y_i,F_i\rangle\over\langle F_i,F_j\rangle}\,F_j
   $ for all edges $(ij)$ in $M$;
\item{\rm(iii)} $|Z|^2\geq0$ and $|Z|^2=0$ if and only if $Z\in f$.

Here, (iii) follows from (i) since $f^\perp$ carries a positive semi-definite
metric with only $f$ as a null direction.
Also note that, since $|Z|^2$ is constant, either $Z\parallel F$ at all points
or, without loss of generality, $|Z|^2\equiv1$ as the space of polynomial
conserved quantities is linear.

We will be mostly interested in the latter case.
To interpret this situation geometrically first note that, if $|Z|^2\equiv1$,
then $i\mapsto Z_i$ defines a discrete sphere congruence so that every sphere
$Z_i$ contains the point $f_i$, by (i); moreover, the sphere
$$
   S_{ij} :
   = Z_i+a_{ij}{\langle Y_j,F_j\rangle\over\langle F_i,F_j\rangle}\,F_i
   = Z_j+a_{ij}{\langle Y_i,F_i\rangle\over\langle F_i,F_j\rangle}\,F_j
\eqno CurvatureSphere$$
belongs to both contact elements\footnote{Recall that, in contrast to M\"obius
geometry, Lie geometry considers oriented spheres; thus $\pm Z+f$ will define
two contact elements which differ by the orientation of their spheres.}
defined by $Z$ at the endpoints of an edge so that
$$
   M \ni i \mapsto Z_i+f_i := \{Z_i+\alpha F_i\,|\,\alpha\in\R\}
$$
defines a discrete principal net in Lie geometry,
see \href[ref.bosu06a]{[6,~Sect.~4.1]},
with curvature spheres $S_{ij}$.
On the other hand, the existence of the curvature spheres $S_{ij}$,
which touch both spheres $Z_i$ and $Z_j$, can be interpreted as a
discrete version of the enveloping condition for a sphere congruence
defined at the vertices of a discrete net.
Thinking of $f$ as a net in $\R^3$, the spheres $Z$ define a unit normal
field $n$ at the vertices of $f$, which satisfies the trapezoid property
of \href[ref.sch07]{[22,~Sect.~3]}.

We summarize these observations in the following\footnote{We insist
on a consistent orientation of the spheres of an enveloped sphere
congruence.}

\proclaim\htag[def.Envelope]{Cor.~\& Def.~3.10}.
If $\pcq(\lambda)=\lambda^NZ+\dots+Q$ is a polynomial conserved quantity
of $f$ with $|Z|^2=1$, then $f$ envelops the discrete sphere congruence
$Z$:
we say that $f:M\to S^3$ \define{envelops} a discrete sphere congruence
$S:M\to S^{3,1}$ if
\item{\rm(i)} $S_i\in f_i^\perp$ for all $i\in M$ (incidence)
   and
\item{\rm(ii)} $S_j=S_i \bmod f_i\oplus f_j$ for each edge $(ij)$ of $M$
   (touching);
\item[0pt]{} the common sphere $S_{ij}$ of the two contact elements
$S_i+f_i$ and $S_j+f_j$ given by $S_i$ and $S_j$ at the endpoints
of an edge will be called a \define{curvature sphere}.

Note that the condition for a sphere congruence to be enveloped by a
net is a condition on the congruence of contact elements defined by
the sphere congruence --- hence, if $f$ envelops $S$ and $F$ is any
light cone lift of $f$, then any sphere congruence $S+hF$, $h:M\to\R$,
is also enveloped by $f$.

Before focusing on this geometrically interesting configuration we shall,
for the rest of this section, investigate the degenerate case where we
can, without loss of generality, take $Z=F$ as a canonical light cone
lift\footnote{We neglect the case where $Z$ may have zeroes, i.e.,
where the degree of $\pcq(\lambda)$ may not be constant.} of $f$.

\proclaim\htag[thm.degPCQ]{Lemma 3.11}.
If $\pcq(\lambda)=\lambda^NF+\dots+Q$ is a polynomial conserved quantity
of $f$, then $F$ is a Moutard lift of $f$.

\proof
From (ii) of \href[thm.EnvelopingZ]{Cor.~3.9} we learn that
$$
     \langle F_i,F_j\rangle + a_{ij}\langle Y_j,F_j\rangle
   = \langle F_i,F_j\rangle + a_{ij}\langle Y_i,F_i\rangle
   = 0
$$
since $F_i$ and $F_j$ are linearly independent;
thus\footnote{This we already knew from \href[thm.pcqNorm]{Lemma 3.5},
because $\langle Y,F\rangle$ is the $\lambda^{2N-1}$-coefficient of
$|\pcq(\lambda)|^2$.} $\langle Y,F\rangle\equiv const=:c$ and $F$ satisfies
\reqn{MoutardCond} with the cross ratio factorizing function $\tilde a:=-c\,a$.
\endproof

As an example we seek an isothermic net with a degenerate degree $1$
polynomial conserved quantity: this is the lowest degree possible
since the top term of a degree $0$ polynomial conserved quantity
is constant and can therefore not be a (Moutard) lift of an
isothermic net.

Now consider, as before,
$$
   |(x_0,\dots,x_4)|^2 = -x_0^2 + \sum_{i=1}^4x_i^2
$$
as the quadratic form of the Minkowski scalar product of $\R^{4,1}$
and the isothermic net
$$
   \{-1,0,1\}^2 \ni (m,n) \mapsto f_{(m,n)} :
   = (\eta m,{1+\alpha\over2}+(-1)^n{1-\alpha\over2},\beta n) \in \R^3,
\eqno explMain$$
where $\alpha\in(0,1)$ and $\beta,\eta>0$, and let
$$
   Z_{(m,n)} := (-1)^nF_{(m,n)}
   \quad{\rm and}\quad
   Q := -{4\over1-\alpha}({1+\alpha\over2},0,1,0,-{1+\alpha\over2}),
\eqno explDegPCQ$$
where $F=({1+|f|^2\over2},f,{1-|f|^2\over2})$ is a Euclidean lift%
\footnote{Cf.~\reqn{EuclideanLift} and \reqn{QEuclideanLift}:
this is the Euclidean lift with respect to $Q_0=(1,0,0,0,-1)$,
defining a flat quadric of constant curvature via \reqn{SpaceForm}.}.
Note that $|Q|^2>0$ so that $Q$ describes a sphere\footnote{%
Suppose that $|Q|^2\leq0$ for a degenerate linear conserved quantity
$\pcq(\lambda)=\lambda F+Q$ of an isothermic net $f$;
then
$$\textstyle
   a_{ij} = \langle F_i,F_j\rangle = -{1\over2}|dF_{ij}|^2 < 0
$$
for any edge $(ij)$ since $dF_{ij}\perp Q$ by
\href[thm.pcqNorm]{Lemma~3.5} and $f$ is regular.
Hence the cross ratio of any face becomes positive;
thus, if we seek an isothermic net with embedded faces,
then $Q$ necessarily describes a sphere.}
and $f$ appears to be a perturbation of a net on that sphere.

Since
$$
   \langle F_i,F_j\rangle = -{1\over2}|f_j-f_i|^2
$$
and $f$ has rectangular faces, so that the cross ratio \reqn{ProjCrossRatio}
on an elementary quadrilateral becomes
$$
   [f_i;f_j;f_k;f_l] = -{|f_i-f_j|^2\over|f_i-f_l|^2},
$$
$a_{ij}:={1\over2}\langle Z_i,Z_j\rangle$ provides a cross ratio factorizing
function, i.e., $Z$ is a Moutard lift of $f$.
Moreover
$$
   \langle\lambda Z+Q,Z\rangle
   = \langle Q,Z\rangle
   \equiv -2,
$$
so that
$$
   d(\lambda Z+Q)_{ij}
   - {\lambda a_{ij}\over\langle Z_i,Z_j\rangle}\{
      \langle\lambda Z+Q,Z\rangle_jZ_i-\langle\lambda Z+Q,Z\rangle_iZ_j \}
   = 0
$$
and $\pcq(\lambda):=\lambda Z+Q$ is a linear conserved quantity for $f$
by \href[thm.PCQcondition]{Lemma 3.7}.

Note that $Z=\pcq(\infty)$ is a lift of $f$ and that
$$
   \pcq({4\over(1-\alpha)^2})
   = {4\over(1-\alpha)^2}\{Z-2{\langle Z,Q\rangle\over|Q|^2}\,Q\}
$$
is a lift of a M\"obius equivalent net or, more precisely,
of the ``antipodal'' net in the quadric of constant curvature
given by $Q$, see \href[ref.imdg]{[18,~Section~1.4]}.
We shall come back to this observation later.

Finally observe that our restriction of the domain to $\{-1,0,1\}^2$ was
not necessary: $f$, as given by \reqn{explMain}, can be extended to all of
$\Z^2$ while \reqn{explDegPCQ} keeps defining a linear conserved quantity for
$f$.
However, this restriction will be convenient later when we shall
recycle this example to demonstrate another aspect of our theory.

\h2 Special isothermic nets.
As in the smooth case, see \href[ref.sa07]{[21, Sect.~2.2]},
we can now use the existence of a polynomial conserved quantity
to define a special class of discrete isothermic nets, ordered by the
(minimal) degree of an associated polynomial conserved quantity.
However, in contrast to the smooth case, where a polynomial
conserved quantity is essentially unique because its top degree
coefficient encodes the conformal Gauss map of the underlying isothermic
surface, the space of polynomial conserved quantities of minimal degree
may be higher dimensional\footnote{However, we do expect uniqueness of a
(normalized) polynomial conserved quantity of minimal degree for generic
isothermic nets of a sufficient size.} and may contain elements with null
top degree coefficient, as in the above example, even though we require the
existence of a normalized polynomial conserved quantity of minimal degree:

\proclaim\htag[def.special]{Def.~3.12}.
A polynomal conserved quantity $\pcq(\lambda)=\lambda^NZ+\dots+Q$ of
an isothermic net $f$ will be called \define{normalized} if $|Z|^2\equiv1$;
we say that $f$ is a \define{special isothermic net of type $N$} if it
has a normalized polynomial conserved quantity of degree $N$,
but not of any lower degree.

As a direct consequence of this definition and the $1$-parameter group
property \reqn{1ParameterT} of the Calapso transformations we obtain
stability of the class of special isothermic nets of a fixed type $N$
under the Calapso transformation:

\proclaim\htag[thm.SpecialCalapso]{Thm.~3.13}.
If $f$ is special isothermic of type $N$
then so are its Calapso transforms $f^\mu=T^\mu f$.

\proof
$T^{\mu,\lambda}=T^{\mu+\lambda}(T^\mu)^{-1}$ are the Calapso transformations
of $f^\mu$ by \reqn{1ParameterT}, see also \href[ref.imdg]{[18,~\S5.7.30]}.

Hence, if $\pcq(\lambda)=\lambda^NZ+\dots+Q$ is a polynomial conserved
quantity of $f$, then
$$
   \pcq^\mu(\lambda) := T^\mu\pcq(\mu+\lambda)
$$
defines a polynomial conserved quantity of $f^\mu$ of the same degree;
moreover,
$$
   |\pcq^\mu(\lambda)|^2
   = |\pcq(\lambda+\mu)|^2
   = \lambda^{2N}|Z|^2 + \lambda^{2N-1}\dots + \dots
$$
showing that $\pcq^\mu(\lambda)$ is normalized as soon as $\pcq(\lambda)$ is.
Thus $f^\mu$ is special isothermic of type $\leq N$.

On the other hand $f=T^{\mu,-\mu}f^\mu$ is a Calapso transform of $f^\mu$,
by \reqn{1ParameterT} again, so that the same argument shows that
$f^\mu$ is special isothermic of type $\geq N$.
\endproof

So far it is rather unclear how restrictive the definition of special
isothermic of a given type is: as we shall see later (and as one might
expect after an equation count), the condition on a discrete isothermic
net of being special isothermic of type $N$ is not a local condition
--- in particular, we will see that every isothermic $3\times3$-net
has plenty of linear conserved quantities and is, therefore, special
isothermic of type $1$.
On the other hand, the condition of being special isothermic of type $0$
does already impose a condition on a $3\times3$-grid, and the corresponding
constant conserved quantity is generically unique:

\proclaim\htag[thm.type0]{Thm.~3.14}.
$f$ is special isothermic of type $0$ if and only if $f$ takes values
in a $2$-sphere.

\proof
First suppose that $f$ is special isothermic of type $0$ with
(necessarily\footnote{Remember that a constant conserved quantity
cannot be degenerate by \href[thm.EnvelopingZ]{Cor.~3.9 (i)}, since $f$
is not constant.} normalized) constant conserved quantity $\pcq(\lambda)=Q$.
Then $Q$ defines a fixed $2$-sphere, since $|Q|^2=1$, that contains
the points of the net by \href[thm.EnvelopingZ]{Cor.~3.9 (i)}.
This also shows that a (normalized) constant conserved quantity is unique
as soon as the isothermic net does not take values in a circle.

Conversely, if $f$ takes values in a fixed $2$-sphere $S\subset S^3$ then the
Calapso transformations $T^\lambda$ of $f$ take (up to a constant M\"obius
transformation) values in the M\"obius group $\Mob(S)$ of this $2$-sphere
(see \href[ref.je00]{[17,~Thm.~3.14]} or \href[ref.imdg]{[18,~\S5.7.22]}).
That is, the Calapso transformations $T^\lambda$ of $f$ can be chosen to fix
a unit vector $Q$ defining the $2$-sphere $S$: this yields a normalized
constant conserved quantity so that $f$ is special of type $0$.
\endproof

As a more involved example we take up surfaces of revolution and
investigate the symmetry of a corresponding polynomial conserved
quantity:
suppose $F_{(m,n)}=(-1)^m(M_m+\Phi_nC)$ is the Moutard lift of a
discrete surface of revolution, see \reqn{RevolutionNet},
with a rotationally symmetric polynomial conserved quantity, i.e.,
we assume that $\Phi_n^{-1}\pcq_{(m,n)}(\lambda)=:\hat\pcq_m(\lambda)$
does not depend on $n$.
Then, clearly, $$
   p_{(m,n)}(\lambda) :
   = \langle\pcq_{(m,n)}(\lambda),F_{(m,n)}\rangle
   = (-1)^m\langle\hat\pcq_m(\lambda),M_m+C\rangle
$$
is independent of $n$; we shall see that the converse also holds:

\proclaim\htag[thm.symmetricPCQ]{Lemma 3.15}.
A polynomial conserved quantity $\pcq(\lambda)$ of a discrete net
$f$ of revolution with canonical lift $F_{(m,n)}=(-1)^m\Phi_n(M_m+C)$ is
rotationally symmetric, $\pcq_{(m,n)}(\lambda)=\Phi_n\hat\pcq_m(\lambda)$,
if and only if
$$
   p(\lambda) := \langle\pcq(\lambda),F\rangle
$$
does not depend on $n$.

\proof
We already know that $p(\lambda)$ depends only on $m$ if $\pcq(\lambda)$
is rotationally symmetric;
to prove the converse we assume that $p(\lambda)$ does not depend
on $n$ and write
$$
   \pcq_{(m,n)}
   = \pcq^\perp_{(m,n)} + \Phi_n(\alpha_{(m,n)}C+\beta_{(m,n)}C^\perp)
   \in \R^{2,1} \oplus \R^2,
\eqno symmPCQdecomposition$$
where $C^\perp$ complements $C$ to form an orthonormal basis of $\R^2$
and $\alpha$ and $\beta$ are suitable real valued functions.
First note that
$$
   (-1)^mp_m = \langle\pcq^\perp_{(m,n)},M_m\rangle + \alpha_{(m,n)}.
$$
Now we fix $m$ and, in order to simplify notation, consider all functions
as functions of $n$ only.
From \reqn{dPCQ} we get
$$
   0 = d\pcq_{n,n+1}(\lambda) + \lambda\,p(\lambda)\,dF_{n,n+1};
$$
hence, in particular,
$$
   d\pcq^\perp_{n,n+1} = 0
   \quad{\rm and}\quad
   d\alpha_{n,n+1} = -\langle d\pcq^\perp_{n,n+1},M\rangle = 0.
$$
For a function $g$ on the vertices let $g_{n,n+1}:={g_n+g_{n+1}\over2}$
denote the associated function on the edges, and let $d\varphi_{n,n+1}$
be the rotation angle along the edge, i.e., the rotation angle of
$\Phi_{n+1}\Phi_n^{-1}$;
observe that
$$\matrix{
   d(\Phi C^\perp)_{n,n+1}
      &=& -2\tan{d\varphi_{n,n+1}\over2}\,(\Phi C)_{n,n+1}
      &\perp& dF_{n,n+1}, \cr
   2(\Phi C^\perp)_{n,n+1}
      &=&   \cot{d\varphi_{n,n+1}\over2}\,d(\Phi C)_{n,n+1}
      &\parallel& dF_{n,n+1}. \cr
}$$
Hence the $\R^2$-part of
$0=d\pcq_{n,n+1}(\lambda)+\lambda\,p(\lambda)\,dF_{n,n+1}$
yields
$$
   0 = \beta_{n,n+1}(\lambda)
   \quad{\rm and}\quad
   0 = \alpha(\lambda)+(-1)^m\lambda p(\lambda)
   + {1\over2}d\beta_{n,n+1}(\lambda)\cot{d\varphi_{n,n+1}\over2}.
$$
Thus, considering two consecutive edges,
$0=\beta_n\sin{d\varphi_{n-1,n}+d\varphi_{n,n+1}\over2}$
so that $\beta$ vanishes and
$$
   \pcq_{(m,n)}(\lambda)
   = \pcq^\perp_m(\lambda) + (-1)^{m+1}\lambda p_m(\lambda)\,\Phi_n C
\eqno symmetricPCQ$$
is clearly rotationally symmetric.
\endproof

As a simple consequence we see that a linear conserved quantity
of a discrete surface of revolution $f$ is rotationally symmetric
as soon as $f$ is rotationally symmetric in the space form defined
by the constant term:

\proclaim\htag[thm.symmetricLCQ]{Cor.~3.16}.
A linear conserved quantity $\pcq(\lambda)=\lambda Z+Q$ of a net $f$
of revolution is rotationally symmetric if and only if $Q$ is.

\proof
By \reqn{symmetricPCQ}, the constant term of a rotationally symmetric
polynomial conserved quantity has no $\R^2$-component in the decomposition
\reqn{symmPCQdecomposition};
for the converse observe that, in the case of a linear conserved quantity,
$p_{(m,n)}=\langle Q,F_{(m,n)}\rangle$.
\endproof

\h1 The B\"acklund transformation.
Bianchi's B\"acklund transformation of smooth constant mean curvature
surfaces in Euclidean space turned out to be a special case of the
Darboux transformation of isothermic surfaces:
considering the Darboux transformation as an initial value problem
depending on a real (spectral) parameter, the Darboux transforms of
a constant mean curvature surface turn out to have constant mean
curvature as soon as a certain relation between the initial value
and the parameter is satisfied, see \href[ref.imdg]{[18,~\S5.4.15]};
in fact, they turn out to be the B\"acklund transforms of the surface,
see \href[ref.jepe97]{[15,~Thm.~7]} and \href[ref.inko05]{[20,~Thm.~4.4]}.
A similar fact holds true for discrete constant mean curvature surfaces
in the sense of \href[ref.jehopi99]{[16]}.

In \href[ref.sa07]{[21, Thm.~3.2]}, this condition for a Darboux transform
of a constant mean curvature surface to have (the same) constant mean
curvature again is shown to be a special case of a similar condition
for (smooth) special isothermic surfaces of given type $N$.
Here we shall analyze the situation for discrete special isothermic nets,
giving rise to what we will call the ``B\"acklund transformation'' for
special isothermic nets.

To this end we shall first recall the Darboux transformation for discrete
isothermic nets,
cf.~\href[ref.jehopi99]{[16]} or \href[ref.imdg]{[18,~\S\S5.7.12~\&~5.7.19]}:

\proclaim\htag[def.darboux]{Def.~4.1}.
Let $f:M\to S^3$ be a discrete isothermic net with its family of Calapso
transformations $T^\lambda$.
Then a \define{Darboux transform} of $f$ is a discrete net
$$
   \hat f:M\to S^3 \quad\hbox{\sl so that}\quad
   \exists\mu\in\R:T^\mu\hat f \equiv const.
$$

Equivalently, we can characterize a Darboux transform by the condition
$$
   \hat f_i
   = \Gamma^\mu_{ij}\hat f_j
   = \Gamma^{1-a_{ij}\mu}_{f_i,f_j}\hat f_j
$$
for all edges $(ij)$ of $M$: remember that the isothermic family of connections
\reqn{IsothConnection} descends to a family of connections on $M\times S^3$;
thus a Darboux transform $\hat f$ can be thought of as a parallel section
of a connection $\Gamma^\mu$ on $M\times S^3$ of the isothermic family of
connections.
Using a cross ratio identity (see \href[ref.imdg]{[18,~\S4.9.11]}),
this condition can be reformulated as the cross ratio condition%
\footnote{Using quaternions or the Clifford algebra of $\R^3$ and thinking
of $S^3=\R^3\cup\{\infty\}$, this is equivalent to the (discrete) Riccati
type equation
$$
   d\hat f_{ij} = \mu\,(\hat f-f)_j\,df^\ast_{ij}\,(\hat f-f)_i,
$$
where $df^\ast_{ij}=a_{ij}(df_{ij})^{-1}$ is the derivative of the Christoffel
transform of $f$, cf.~\href[ref.imdg]{[18,~\S5.7.7]};
the condition $\hat f_i=\Gamma^\mu_{ij}\hat f_j$ appears as the usual
linearization from this Riccati equation: as ``Darboux's linear system''.}
(cf.~\reqn{CrossRatioTrafo})
$$
   [f_i;f_j;\hat f_j;\hat f_i]
   = 1-[\hat f_j;f_i;\hat f_i;f_j]
   = a_{ij}\mu.
\eqno DarbouxCrossRatio$$

Clearly, a Darboux transform $\hat f$ of an isothermic net has a light cone
lift $\hat F:M\to L^4$ so that
$$
   T^\mu\hat F\equiv const.
   \quad\Leftrightarrow\quad
   \hat F_i = \Gamma^\mu_{ij}\hat F_j
\eqno DTcondition$$
on all edges $(ij)$ of $M$, that is, $\hat F$ is a $\Gamma^\mu$-parallel
lift of $\hat f$.
Note that this is exactly the conserved quantity condition \reqn{PCQcondition}
for a light cone map $\hat F$ and a fixed value $\mu$ of the spectral
parameter;
we shall come back to this point later.

The crucial observation is that the Darboux transformation produces
discrete isothermic nets from isothermic nets,
see \href[ref.jehopi99]{[16]} or \href[ref.imdg]{[18,~\S5.7.12]}, where
a proof relying on the hexahedron lemma \href[ref.imdg]{[18,~\S4.9.13]}
was given; here we shall again give an alternative proof, relying on
\href[thm.dic]{Lemma~2.5}, that also provides us with a useful formula
for the Calapso transformations of a Darboux transform,
cf.~\href[ref.je00]{[17]} or \href[ref.imdg]{[18,~\S5.7.35]}:

\proclaim\htag[thm.darbouxT]{Lemma 4.2}.
A Darboux transform $\hat f$, $T^\mu\hat f\equiv const$, of a discrete
isothermic net $f$ with cross ratio factorizing function $a$ is isothermic
with the same cross ratio factorizing function $\hat a=a$ and with Calapso
transformations
$$
   \hat T^\lambda = T^\lambda\Gamma^{1-\lambda/\mu}_{f,\hat f}.
\eqno TDT$$

\proof
We shall show that the isothermic family of connections \reqn{IsothConnection}
of $\hat f$ is given by
$$
   \hat\Gamma^\lambda_{ij} :
   = \Gamma^{1-\lambda a_{ij}}_{\hat f_i,\hat f_j}
   = (\Gamma_{f,\hat f}^{1-\lambda/\mu})^{-1}_i
      \Gamma^\lambda_{ij}
      (\Gamma_{f,\hat f}^{1-\lambda/\mu})_j,
$$
which is then, with $\Gamma^\lambda$, clearly flat so that $\hat f$
is isothermic with cross ratio factorizing function $\hat a=a$ by
\href[thm.dic]{Lemma 2.5};
further, the formula for the Calapso transformations $\hat T^\lambda$
also follows directly.

Thus we have to show that
$$
   \Gamma_{f_i,\hat f_i}^{1-\lambda/\mu}
      \Gamma^{1-\lambda a_{ij}}_{\hat f_i,\hat f_j}
   =
   \Gamma^{1-\lambda a_{ij}}_{f_i,f_j}
      \Gamma_{f_j,\hat f_j}^{1-\lambda/\mu}.
$$
But, with \reqn{DarbouxCrossRatio}, this follows directly from
\href[thm.dicCor]{Lemma 2.6}.
\endproof

At this point we are now in perfect shape to discuss:

\h2 Darboux transforms of special isothermic nets.
As outlined at the beginning of this section, we aim to obtain the
B\"acklund transformations for special isothermic nets of type $N$
as particular classes of Darboux transformations.
The following theorem provides the essential criterion:

\proclaim\htag[thm.darbouxPCQ]{Thm.~4.3}.
Let $\hat f$ be a Darboux transform, $T^\mu\hat f\equiv const$,
of a special isothermic net $f$ of type $N$ with normalized polynomial
conserved quantity $\pcq(\lambda)$ of degree $N$.
Then:
\item{\rm(i)} $\hat f$ is special isothermic of type $\leq N+1$ and
\item{\rm(ii)} $\hat f$ is special isothermic of type $\leq N$
   as soon as $\pcq(\mu)\in\hat f^\perp$.

Note that, if $\pcq_i(\mu)\in\hat f^\perp_i$ for some $i\in M$,
then this holds true for all $i\in M$.
Namely,
if $\hat F$ is a light cone lift of $\hat f$ as in \reqn{DTcondition},
then $\pcq(\mu)$ and $\hat F$ are both parallel sections of the metric
connection $\Gamma^\mu_{ij}$ on $M\times\R^{4,1}$ so that their scalar
product $\langle\pcq(\mu),\hat F\rangle\equiv const$.

\proof
We define
$$
   \hat\pcq(\lambda) :
   = (\lambda-\mu)\,\Gamma^{1-\lambda/\mu}_{\hat f,f}\pcq(\lambda).
\eqno pcqDarboux$$
Clearly $\hat T^\lambda\hat\pcq(\lambda)\equiv const$ by \reqn{TDT}
and, using \reqn{CrossRatioTrafo}, we obtain
$$
   \hat\pcq(\lambda)
   = (\lambda-\mu)\pcq(\lambda) - {1\over\langle F,\hat F\rangle}\{
      {\lambda(\lambda-\mu)\over\mu}\langle\pcq(\lambda),F\rangle\hat F
      + \lambda \langle\pcq(\lambda),\hat F\rangle F
      \},
$$
which is polynomial of degree $\leq N+1$ as $\langle\pcq(\lambda),F\rangle$
has degree $\leq N-1$ by \href[thm.EnvelopingZ]{Cor.~3.9~(i)}.
Finally, writing
$$
   \pcq(\lambda) = \lambda^NZ+\dots+Q
   \quad{\rm and}\quad
   \hat\pcq(\lambda) = \lambda^{N+1}\hat Z+\dots+\hat Q,
$$
we find that
$$
   \lambda^{2N+2}|\hat Z|^2+\dots
   = |\hat\pcq(\lambda)|^2
   = (\lambda-\mu)^2\,|\pcq(\lambda)|^2
   = \lambda^{2N+2}|Z|^2+\dots
$$
for all $\lambda$, so that $|\hat Z|^2=|Z|^2$.
Hence $\hat\pcq(\lambda)$ is a normalized polynomial conserved quantity
of degree $N+1$ and, therefore, $\hat f$ is special isothermic of type
$\leq N+1$, proving (i).

To prove (ii) note that
$$
   \hat\pcq(\mu)
   = -\mu{\langle\pcq(\mu),\hat F\rangle\over\langle F,\hat F\rangle}\,F
   = 0
$$
when $\pcq(\mu)\in\hat f^\perp$.
Hence $\hat f$ has a polynomial conserved quantity of degree $\leq N$ by
\href[thm.pcqReduction]{Lemma~3.3}; in particular,
$$
   \hat\pcq(\lambda) :
   = \Gamma^{1-\lambda/\mu}_{\hat f,f}\pcq(\lambda)
\eqno pcqBacklund$$
provides a normalized polynomial conserved quantity of degree $N$ for
$\hat f$ showing that $\hat f$ is special isothermic of type $\leq N$.
\endproof

\proclaim\htag[def.backlund]{Def.~4.4}.
A \define{B\"acklund transform} $\hat f$ of a special isothermic net $f$
of type $N$ with polynomial conserved quantity $\pcq(\lambda)$ is
a Darboux transform, that is, $T^\mu\hat f\equiv const$,
so that $\pcq(\mu)\in\hat f^\perp$.

The B\"acklund transformation between special isothermic nets of the same
type\footnote{We know that the B\"acklund does not increase the type but,
in special circumstances, if may decrease the type (by $1$).  We shall
come back to this point later.}
is symmetric, that is, if $\hat f$ is a type $N$ B\"acklund transform of a
special isothermic net $f$ of type $N$ with polynomial conserved quantity
$\pcq(\lambda)$ then $f$ is a B\"acklund transform of $\hat f$, where
$\hat\pcq(\lambda)$ is given by \reqn{pcqBacklund}:
firstly, $f$ is a Darboux transform of $\hat f$ since the cross ratio
condition \reqn{DarbouxCrossRatio} is symmetric in $f$ and $\hat f$,
$$
   [\hat f_i;\hat f_j;f_j;f_i]
   = [f_i;f_j;\hat f_j;\hat f_i]
   = a_{ij}\mu,
$$
by a cross ratio identity, see \href[ref.imdg]{[18,~\S4.9.11]}; and
secondly, using \reqn{CrossRatioTrafo}, we find that $f$ satisfies the
B\"acklund condition\footnote{The assumption that $\hat f$ is special
isothermic of the same type $N$ as $f$ ensures that $\hat\pcq(\lambda)$
has, with $\pcq(\lambda)$, minimal degree and, in particular, that
$\hat\pcq(\mu)\neq0$ so that the condition is meaningful.}
$\hat\pcq(\mu)\in f^\perp$:
$$
   \hat\pcq(\mu)
   = \lim_{\lambda\to\mu}\hat\pcq(\lambda)
   = \pcq(\mu)
      - {\langle\pcq(\mu),F\rangle\over\langle F,\hat F\rangle}\,\hat F
      \bmod f.
$$

Recall that constructing a Darboux transform $\hat f$ of a given
isothermic net amounts to determining a parallel isotropic section
$\hat F$ of $M\times\R^{4,1}$ equipped with a connection $\Gamma^\mu$
in the isothermic family, see \reqn{DTcondition}:
after choosing an initial value there is a unique solution.
As the condition on a Darboux transform of a special isothermic net 
to become a B\"acklund transform is preserved by the propagation it
is sufficient to choose the initial value of the Darboux transform
appropriately in order to obtain a B\"acklund transform.

Now note that \reqn{pcqBacklund} has exactly the same structure%
\footnote{This is a common phenomenon in discrete differential
geometry, often referred to as ``multidimensional consistency'',
cf.~\href[ref.bosu06a]{[6]} or \href[ref.bosu06]{[7]}:
the transformations of a discrete net are governed by exactly
the same conditions as the net itself.}
as the conserved quantity condition \reqn{PCQcondition}, with ${1\over\mu}$
taking the role of the $a_{ij}$.
Thus, completely analogous to \href[thm.PCQcondition]{Lemma~3.7},
$$
   \hat\pcq(\lambda) + {\lambda\over\mu}
      {\langle\pcq(\lambda),F\rangle\over\langle F,\hat F\rangle}\,\hat F
   = \pcq(\lambda) + {\lambda\over\mu}
      {\langle\hat\pcq(\lambda),\hat F\rangle\over\langle F,\hat F\rangle}\,F,
\eqno DarbouxPCQ$$
where $F$ and $\hat F$ are any%
\footnote{A canonical choice would be to take $F$ and $\hat F$ to satisfy
\reqn{MoutardCond} and $\langle F,\hat F\rangle\equiv{1\over\mu}$:
this is the ``3D-consistency'' of the condition to be discrete isothermic,
see \href[ref.bosu06]{[7]}.}
light cone lifts of $f$ and $\hat f$, respectively.
In particular, writing
$$
   \pcq(\lambda) = \lambda^NZ + \lambda^{N-1}Y + \dots + Q
   \quad{\rm and}\quad
   \hat\pcq(\lambda) = \lambda^N\hat Z + \lambda^{N-1}\hat Y + \dots + \hat Q,
$$
we find that,
analogous to \href[thm.ConstantQ]{Cor.~3.8}, the constant coefficients
   $\hat Q=Q$ and that,
analogous to the geometric interpretation of \href[thm.EnvelopingZ]{Cor.~3.9},
$$
   S :
   = Z + {\langle\hat Y,\hat F\rangle\over\mu\langle F,\hat F\rangle}\,F
   = \hat Z + {\langle Y,F\rangle\over\mu\langle F,\hat F\rangle}\,\hat F
\eqno DarbouxSphere$$
yields a sphere congruence\footnote{Indeed, a similar fact can be proved
for Darboux transforms in general: given a sphere congruence $Z$ enveloped
by an isothermic net and a Darboux transform $\hat f$ of $f$, the sphere
congruence
$$\textstyle
   S := Z - {\langle Z,\hat F\rangle\over\langle F,\hat F\rangle}\,F
$$
will be enveloped by both $f$ and $\hat f$.
This is a fact about the existence of Ribaucour transforms for discrete
curvature line nets in Lie geometry, cf.~\href[ref.bosu06a]{[6]}.
Note that this is a different ``Darboux sphere congruence'' than the
one discussed in \href[ref.jehopi99]{[16]}, which ``lives'' on the
faces of the domain $M$.} that is enveloped by both $f$ and $\hat f$
in the sense of \href[def.Envelope]{Def.~3.10}.
Thus our B\"acklund transformation for special isothermic nets is a special
type of Ribaucour transformation for discrete curvature line nets in Lie
sphere geometry, see \href[ref.bosu06a]{[6,~Def.~22]}.

\proclaim\htag[thm.pcqDarboux]{Thm.~4.5}.
Let $\hat f$ be a B\"acklund transform of a special isothermic net $f$
of type $N$, with respective polynomial conserved quantities
$$
   \pcq(\lambda)
   = \lambda^NZ + \lambda^{N-1}Y + \dots + Q
   \quad{\sl and}\quad
   \hat\pcq(\lambda)
   = \Gamma^{1-\lambda/\mu}_{\hat f,f}\pcq(\lambda)
   = \lambda^N\hat Z + \lambda^{N-1}\hat Y + \dots + \hat Q.
$$
Then
\item{\rm(i)}
   $\hat Q=Q$, and
\item{\rm(ii)}
   $ Z + {\langle\hat Y,\hat F\rangle\over\mu\langle F,\hat F\rangle}\,F
    = \hat Z + {\langle Y,F\rangle\over\mu\langle F,\hat F\rangle}\,\hat F$;
    and,
\item{\rm(iii)}
   in particular, $Z+f$ and $\hat Z+\hat f$ give rise to a Ribaucour pair
   in the Lie geometric sense.

\h2 Bianchi permutability.
A key feature of Darboux-B\"acklund type transformations for smooth
or discrete classes of surfaces is Bianchi permutability:
given two transforms of a surface there is a fourth (often unique) surface,
which is a simultaneous transform of the two initial transforms --- thus
providing the combinatorics of a quadrilateral for the transformation.
We will refer to such quadrilaterals as \define{Bianchi quadrilaterals}.

In particular, such a theorem holds true for the Darboux transformation
of (discrete) isothermic surfaces,
see \href[ref.jehopi99]{[16]} or \href[ref.imdg]{[18,~\S5.7.28]},
where the fourth surface is uniquely determined\footnote{This is in
contrast to the Ribaucour transformation of (discrete) principal nets,
where a $1$-parameter family of fourth surfaces exists,
see \href[ref.buje06]{[10]} or \href[ref.boje01]{[5]}.}
by a cross ratio condition\footnote{Note how repeated application
of this theorem builds up a discrete isothermic net.}:

\proclaim\htag[thm.Dperm]{Thm.~4.6}.
If $\hat f_1$ and $\hat f_2$ are two Darboux transforms,
$T^{\mu_i}\hat f_i\equiv const$, of a discrete isothermic net $f$,
then
$$
   f_{12} := \Gamma^{\mu_2/\mu_1}_{\hat f_1,\hat f_2}f
$$
is a simultaneous Darboux transform of $\hat f_1$ and $\hat f_2$:
$$
   \hat T_1^{\mu_2}f_{12}\equiv const
   \quad{\sl and}\quad
   \hat T_2^{\mu_1}f_{12}\equiv const.
$$

We wish to prove a similar theorem for the B\"acklund transformation
of special isothermic nets of a given type $N$.

Thus let $f$ be a special isothermic net with polynomial conserved
quantity $\pcq(\lambda)$ and let $\hat f_i$ be two B\"acklund transforms
of $f$, i.e.,
$$
   T^{\mu_i}\hat f_i \equiv const
   \quad{\rm and}\quad
   \pcq(\mu_i) \in \hat f_i^\perp.
$$
Clearly, by the permutability theorem \href[thm.Dperm]{Thm.~4.6} for
the Darboux transformation and \reqn{pcqBacklund}, the simultaneous
B\"acklund transform $f_{12}$ of $\hat f_1$ and $\hat f_2$ is
necessarily
$$
   f_{12} = \Gamma^{\mu_2/\mu_1}_{\hat f_1,\hat f_2}f
   \quad{\rm with}\quad
   \pcq_{12}(\lambda)
   = \Gamma^{1-\lambda/\mu_2}_{f_{12},\hat f_1}\hat\pcq_1(\lambda)
   = \Gamma^{1-\lambda/\mu_1}_{f_{12},\hat f_2}\hat\pcq_2(\lambda)
\eqno PCQperm$$
as a polynomial conserved quantity
--- just as in the proof of \href[thm.darbouxT]{Lemma~4.2} it follows
from \href[thm.dicCor]{Lemma~2.6} that
$$
   \Gamma^{1-\lambda/\mu_2}_{f_{12},\hat f_1}
   \Gamma^{1-\lambda/\mu_1}_{\hat f_1,f}
   =
   \Gamma^{(1-\lambda/\mu_1)/(1-\lambda/\mu_2)}_{\hat f_1,\hat f_2}
   =
   \Gamma^{1-\lambda/\mu_1}_{f_{12},\hat f_2}
   \Gamma^{1-\lambda/\mu_2}_{\hat f_2,f},
$$
so that $\pcq_{12}(\lambda)$ is well defined by \reqn{PCQperm}.
Now it is straightforward to see that
$$
   \pcq_{12}(\mu_1)
   = \lim_{\lambda\to\mu_1}
      \Gamma^{(1-\lambda/\mu_1)/(1-\lambda/\mu_2)}_{\hat f_1,\hat f_2}
      \pcq(\lambda)
   = \pcq(\mu_1) - {\langle\pcq(\mu_1),\hat F_2\rangle
      \over\langle\hat F_1,\hat F_2\rangle}\,\hat F_1 + \dots\hat F_2
   \in \hat f_2^\perp,
$$
showing that $\hat f_2$ is a B\"acklund transform of $f_{12}$;
similarly, $\hat f_1$ is also a B\"acklund transform of $f_{12}$.
The symmetry of the B\"acklund transformation then completes the
proof of the following

\proclaim\htag[thm.BBperm]{Thm.~4.7 (Bianchi permutability)}.
Given two B\"acklund transforms $\hat f_1$ and $\hat f_2$ with parameters
$\mu_1$ and $\mu_2$, respectively, of a special isothermic net $f$ of type $N$,
there is a net $f_{12}$ so that the four nets form a Bianchi quadrilateral:
that is, $f_{12}$ is a B\"acklund transform of $\hat f_1$ with parameter
$\mu_2$ and of $\hat f_2$ with parameter $\mu_1$.

Note that we have not used any new arguments to prove this theorem
---
indeed, using the similarity of the polynomial conserved quantity equations
and the conditions governing the B\"acklund transformation, we could have
formulated a proof based on the fact that ``3D-consistency'' of a
``2D-system'' implies higher dimensional consistency,
see \href[ref.bosu06]{[7,~Thm.~7]}:
in the case at hand we were interested in ``4D-consistency''.
Thus, any higher dimensional permutability theorems can now be proved
by purely combinatorial arguments\footnote{This is in contrast with
the Ribaucour transformation, where $3$-dimensional permutability,
i.e., a ``Bianchi cube'' theorem, is the critical case as there
is no uniqueness in the Bianchi quadrilateral.}.
For example, we can now argue that a ``Bianchi cube'' can be (uniquely)
constructed from a special isothermic net and three B\"acklund transforms:
the existence of the eighth B\"acklund transform is ensured by the very same
fact that ensured the existence of a Darboux transform and the compatibility
of the B\"acklund transformation with the construction, as discussed above.

\h2 Complementary nets.
As we already noticed earlier, the equation \reqn{DTcondition} on a light
cone map $\hat F$ to provide a Darboux transform of an isothermic net $f$
is exactly the conserved quantity equation \reqn{PCQcondition} for a fixed
parameter $\mu$.
Consequently, any zero $\mu$ of $|\pcq(\lambda)|^2$ provides a (light cone
lift of a) Darboux transform $\hat F=\pcq(\mu)$ of an isothermic net $f$
with polynomial conserved quantity $\pcq(\lambda)$.
Moreover, since $\langle\hat F,\pcq(\mu)\rangle=|\pcq(\mu)|^2=0$,
this Darboux transform will, in fact, be a B\"acklund transform of $f$.
In this section we shall discuss the role of these special B\"acklund
transforms of a special isothermic net.

\proclaim\htag[def.complementary]{Def.~4.8}.
Let $f$ be a special isothermic net with polynomial conserved quantity
$\pcq(\lambda)$;
those B\"acklund transforms $\hat f$ of $f$ given by $\hat F=\pcq(\mu)$,
where $|\pcq(\mu)|^2=0$, are the \define{complementary nets} of $f$.

Clearly, a special isothermic net of type $N$ has at most $2N$ complementary
nets and, as $|\pcq(\lambda)|^2$ is an even degree polynomial, there may
be no (real) complementary nets:
for example, let $f$ be a type $1$ special isothermic net with linear
conserved quantity
$$
   \pcq(\lambda) = \lambda Z+Q
   \quad{\rm and}\quad
   H := -\langle Z,Q\rangle, \quad \kappa := -|Q|^2;
\eqno MeanCurvature$$
then the complementary nets of $f$ are given by\footnote{Note that
$\pcq(H\pm\sqrt{H^2+\kappa})\neq0$ since $f$ is special isothermic of
type $1$ so that $\pcq(\lambda)$ cannot have zeroes.}
$$
   \pcq(H \pm \sqrt{H^2+\kappa})
   = (H \pm \sqrt{H^2+\kappa})\,Z + Q
\eqno cmcComplementary$$
so that the number of (real) complementary nets depends on the sign
of $H^2+\kappa$.

On the other hand, if enough complementary nets are known, then the
corresponding polynomial conserved quantity can be reconstructed:
first we observe that, given $N+1$ distinct parameter values
$\lambda=\mu_0,\dots,\mu_N$,
$$
   \pcq(\lambda)
   = \sum_{n=0}^N\pcq(\mu_n)\prod_{m\neq n}{\lambda-\mu_m\over\mu_n-\mu_m},
\eqno pcqInterpol$$
where $\pcq(\mu_n)$ are $\Gamma^{\mu_n}$-parallel and,
with $\alpha_n:=\prod_{m\neq n}{1\over\mu_n-\mu_m}$,
the leading coefficient
$$
   Z = \sum_{n=0}^N\alpha_n\pcq(\mu_n) \in f^\perp
$$
of $\pcq(\lambda)$ is a (constant) linear combination of the $\pcq(\mu_n)$.
These are the assumptions that will allow us to reconstruct a polynomial
conserved quantity from $N+1$ suitable Darboux transforms:

\proclaim\htag[thm.pcqReconstruction]{Lemma 4.9}.
Let $\hat F^n$, $n=0,\dots,N$, be $\Gamma^{\mu_n}$-parallel sections
of $M\times\R^{4,1}$ for pairwise distinct $\mu_n$ and suppose that
$$
   Z_i = \sum_{n=0}^N\alpha_n\hat F^n_i \in f_i^\perp
$$
for some constants $\alpha_n\in\R$ and all $i\in M$.
Then
$$
   \pcq(\lambda) :
   = \sum_{n=0}^N\alpha_n\hat F^n\prod_{m\neq n}(\lambda-\mu_m)
\eqno pcqReconstruction$$
is a degree $N$ polynomial conserved quantity for $f$ with top degree
coefficient $Z$;
if $|Z_i|^2=1$ at some $i\in M$, then \reqn{pcqReconstruction}
defines a normalized polynomial conserved quantity.

Note that, if $|Z|^2>0$, then $f$ envelops the sphere congruence $Z$
in the sense of \href[def.Envelope]{Def.~3.10}, as desired:
incidence is given by assumption and touching is guaranteed because
the sections $\hat F^n$ are $\Gamma^{\mu_n}$ parallel so that
$$
   dZ_{ij}
   = \sum_{n=0}^N\alpha_n\,d\hat F^n_{ij}
   = 0 \bmod f_i\oplus f_j.
$$
This also implies directly that $|Z|^2\equiv const$ on $M$.

\proof
It remains to show that $\pcq(\lambda)$ is indeed a polynomial conserved
quantity, that is, we wish to show that
$$
   \Gamma^\lambda_{ij}\pcq_j(\lambda) - \pcq_i(\lambda)
   = d\pcq_{ij}(\lambda)
   + {a_{ij}\lambda\over\langle F_i,F_j\rangle}\{
      {1\over1-a_{ij}\lambda}\langle\pcq_j(\lambda),F_i\rangle F_j
      -\langle\pcq_j(\lambda),F_j\rangle F_i
      \}
   = 0
$$
for all $\lambda$ and each edge $(ij)$ of $M$.
Clearly, this equality holds true for $\lambda=\mu_n$, $n=0,\dots,N$,
since
$$
   \Gamma^{\mu_n}_{ij}\pcq_j(\mu_n)
   = \alpha_n\Gamma^{\mu_n}\hat F^n_j\prod_{m\neq n}(\mu_n-\mu_m)
   = \alpha_n\hat F_i\prod_{m\neq n}(\mu_n-\mu_m)
   = \pcq_i(\mu_n)
$$
so that it holds for all $\lambda$ as soon as we know that
$\Gamma^\lambda_{ij}\pcq_j(\lambda)-\pcq_i(\lambda)$ is a
degree $N$ polynomial.
By definition $\pcq_i(\lambda)$ and $\pcq_j(\lambda)$ are degree $N$
polynomials, and $\langle\pcq_j(\lambda),F_j\rangle$ is a degree
$N-1$ polynomial by incidence, $Z_j\perp F_j$.
Moreover,
$$\matrix{
   \langle\pcq_j({1\over a_{ij}}),F_i\rangle
   &=& {1\over a_{ij}^N}\prod_{m=0}^N(1-a_{ij}\mu_m)
     \sum_{n=0}^N{\alpha_n\over1-a_{ij}\mu_n}\langle\hat F^n_j,F_i\rangle
     \hfill\cr
   &=& {1\over a_{ij}^N}\prod_{m=0}^N(1-a_{ij}\mu_m)
     \sum_{n=0}^N\alpha_n\langle\hat F^n_i,F_i\rangle \hfill\cr
   &=& {1\over a_{ij}^N}\prod_{m=0}^N(1-a_{ij}\mu_m)
     \langle Z_i,F_i\rangle \hfill\cr
   &=& 0 \hfill\cr
}$$
so that ${1\over1-a_{ij}\lambda}\langle\pcq_j(\lambda),F_i\rangle$
is also a degree $N-1$ polynomial.
Hence the claim follows.
\endproof

Note that, in \href[thm.pcqReconstruction]{Lemma~4.9}, we did not require
the $\hat F^n$ to be isotropic sections of $M\times\R^{4,1}$, that is,
we did not require them to be lifts of Darboux transforms of $f$.

Now suppose that $\mu$ is a simple zero of $|\pcq(\lambda)|^2$ and let
$\hat F=\pcq(\mu)$ denote (a lift of) the corresponding complementary
net of $f$.
Then, from \reqn{pcqBacklund},
$$
   \hat\pcq(\lambda)
   = \pcq(\lambda)
   - {\lambda\over\mu}
      {\langle\pcq(\lambda),F\rangle\over\langle\pcq(\mu),F\rangle}\,\pcq(\mu)
   - {\lambda\over\lambda-\mu}
      {\langle\pcq(\lambda),\pcq(\mu)\rangle\over\langle F,\pcq(\mu)\rangle}\,F
$$
such that
$$
   \hat\pcq(\mu)
   = -{\mu\over\langle\pcq(\mu),F\rangle}
      \lim_{\lambda\to\mu}
      {\langle\pcq(\lambda),\pcq(\mu)\rangle\over\lambda-\mu}\,F
   = -{\mu\over2\langle\pcq(\mu),F\rangle}
      \lim_{\lambda\to\mu}
      {|\pcq(\lambda)|^2\over\lambda-\mu}\,F
$$
since
$$
   {|\pcq(\lambda)|^2\over\lambda-\mu}
   - 2{\langle\pcq(\lambda),\pcq(\mu)\rangle\over\lambda-\mu}
   = (\lambda-\mu)\,\left|{\pcq(\lambda)-\pcq(\mu)\over\lambda-\mu}\right|^2
   \longrightarrow 0\cdot|\pcq'(\mu)|^2
   = 0
$$
as $\lambda\to\mu$.
Consequently,
if $\mu$ is a simple zero of $|\pcq(\lambda)|^2$,
then $f$ is a complementary net of $\hat f$,
that is, the notion of complementary nets is symmetric.

If, on the other hand, $\mu$ is a higher order zero of $|\pcq(\lambda)|$,
then $\hat\pcq(\mu)=0$ so that $\hat f$ is special isothermic of lower type
than $f$ by \href[thm.pcqReduction]{Lemma~3.3}.
Indeed, all type lowering B\"acklund transformations arise in this way:
suppose that $\hat f$ is a type $N-1$ B\"acklund transform of a special
isothermic net $f$ of type $N$ --- or, otherwise said, $f$ is a Darboux
transform of $\hat f$ which is not a B\"acklund transform.
Then, their polynomial conserved quantities are related by \reqn{pcqDarboux}:
$$
   \pcq(\lambda)
   = (\lambda-\mu)\Gamma^{1-\lambda/\mu}_{f,\hat f}\hat\pcq(\lambda)
   = (\lambda-\mu)(\hat\pcq(\lambda)-{\lambda\over\mu}
      {\langle\hat\pcq(\lambda),\hat F\rangle\over\langle F,\hat F\rangle}\,F)
   - \lambda{\langle\hat\pcq(\lambda),F\rangle\over\langle F,\hat F\rangle}\,
      \hat F
$$
so that
$$
   \pcq(\mu)
      = -\mu{\langle\hat\pcq(\mu),F\rangle\over\langle F,\hat F\rangle}\,\hat F
$$
spans $\hat f$, which is therefore a complementary net of $f$.

Note that $\mu$ is a higher order zero of
$|\pcq(\lambda)|^2=(\lambda-\mu)^2|\hat\pcq(\lambda)|^2$.

We summarize these results:

\proclaim\htag[thm.cmcOne]{Lemma 4.10}.
Let $f$ be special isothermic of type $N$ with polynomial conserved quantity
$\pcq(\lambda)$.
\item{\rm(i)} If $\mu$ is a higher order zero of $|\pcq(\lambda)|^2$,
   then $\hat F=\pcq(\mu)$ defines a type $N-1$ B\"acklund transform of $f$.
\item{\rm(ii)} If $\hat f$ is a type $N-1$ B\"acklund transform of $f$,
   then $\hat f$ is a complementary net of $f$, $\hat f\ni\pcq(\mu)$ for
   some $\mu$, where $\mu$ is a higher order zero of $|\pcq(\lambda)|^2$.

As a consequence of this lemma, a special isothermic net of type $N$ can
have at most $N$ B\"acklund transforms of type $N-1$ --- generically, a
B\"acklund transformation is between special isothermic nets of the same
type, and is therefore symmetric as discussed above.
For example, consider a special isothermic net $f$ of type $1$:
by \reqn{cmcComplementary}, $f$ is a Darboux transform of a type $0$ net,
that is, of a spherical net (see \href[thm.type0]{Lemma~3.14}), if and only
if $H^2+\kappa=0$.

In the remainder of this section we shall discuss geometric properties
of complementary nets of special isothermic nets of type $1$ and $2$.

First consider a type $1$ special isothermic net with two complementary
nets, i.e., $H^2+\kappa>0$ in \reqn{cmcComplementary};
also, we assume $\kappa\neq0$, excluding the degenerate case, where one
of the complementary nets becomes constant.
Then
$$
   \hat F^\pm := Z + {1\over\mu_\pm}\,Q,
   \quad{\rm where}\quad
   \mu_\pm := H \pm \sqrt{H^2+\kappa},
\eqno cmcClifts$$
provides $\Gamma^{\mu_\pm}$-parallel light cone lifts of the two
complementary nets so that
$$
   \hat F^\pm
   = \hat F^\mp \pm {2\sqrt{H^2+\kappa}\over\kappa}\,Q
   = \hat F^\mp - 2{\langle\hat F^\mp,Q\rangle\over|Q|^2}\,Q.
$$
Thus the two nets are M\"obius equivalent and, more precisely,
they are ``antipodal'' in the quadric of constant curvature
given by $Q$, see \href[ref.imdg]{[18,~Section~1.4]}.
Moreover, the ``orthogonal circles''
$$
   \hat c^\pm := \span\{F,Z,\hat F^\pm\}
$$
of the Ribaucour pair\footnote{In the smooth case, these orthogonal circles
form a \define{cyclic system}, that is, they have a $1$-parameter family
of orthogonal surfaces so that any two orthogonal surfaces are Ribaucour
transforms of each other.}
$(f,\hat f^\pm)$, see \href[thm.pcqDarboux]{Thm.~4.5},
coincide in corresponding points:
$$
   \hat c^\pm_i
   = \span\{F_i,Z_i,\hat F^\pm_i\}
   = \span\{F_i,Z_i,Q\}.
$$

For a type $2$ special isothermic net $f$ we obtain two similar properties:
first let
$$
   \hat F^n = \pcq(\mu_n) = \mu_n^2 Z + \mu_n Y + Q
   \quad(n=1,2)
$$
denote ($\Gamma^{\mu_n}$-parallel lifts of) two complementary nets of $f$;
then the planes (in the quadric $Q$)
$$
   \hat e^n := \span\{F,Z,\hat F^n,Q\}
$$
of the orthogonal circles $\hat c^n$ coincide at corresponding points:
$$
   \hat e^n_i
   = \span\{F_i,Z_i,\hat F^n_i,Q\}
   = \span\{F_i,Z_i,Y_i,Q\}.
$$

On the other hand, if
$$
   \hat F^n = \pcq(\mu_n) = \mu_n^2 Z + \mu_n Y + Q
   \quad(n=1,2,3)
$$
provide three complementary nets of $f$, then
$$
   Q \in \hat c_i := \span\{\hat F^n_i\,|\,n=1,2,3\}
$$
for all $i\in M$, that is, the circles through corresponding points of
the three complementary nets are, in fact, straight lines in the quadric
of constant curvature given by $Q$.

Reversing these observations yields four constructions for linear or quadratic
conserved quantities from geometric configurations of Darboux transforms:

\proclaim\htag[thm.DoubleDarboux]{Thm.~4.11}.
Let $\hat f^n$, $n=1,2$, be two Darboux transforms with different parameters
$\mu_n$ of a discrete isothermic net $f$ so that the circles
$$
	\hat c^n_{ij} :
	= \span\{F_i,F_j,\hat F^n_i\}
	= \span\{F_i,F_j,\hat F^n_j\}
$$
on the edges of $M$ do not coincide for $n=1,2$.
Suppose that the $\hat f^n$ are antipodal in a suitable non-Euclidean
space form.
Then $f$ has a normalized linear conserved quantity.

\proclaim\htag[thm.TripleDarboux]{Thm.~4.12}.
Let $\hat f^n$, $n=1,2,3$, be three Darboux transforms with different parameters
$\mu_n$ of a discrete isothermic net $f$ so that the circles\footnote{Note that
$\span\{F_i,F_j,\hat F^n_i\}=\span\{F_i,F_j,\hat F^n_j\}$ since corresponding
edges of a Darboux (or, more generally, Ribaucour pair) have concircular
endpoints.}
$$
	\hat c^n_{ij} :
	= \span\{F_i,F_j,\hat F^n_i\}
	= \span\{F_i,F_j,\hat F^n_j\}
$$
on the edges of $M$ are not cospherical for $n=1,2,3$.
Suppose that the circles
$$
   \hat c_i := \span\{\hat F^n_i\,|\,n=1,2,3\}
$$
are straight lines in a suitable space form geometry.
Then $f$ has a normalized quadratic conserved quantity.

\proof
We prove both theorems.
Thus let $N=2$ or $N=3$ and let $\hat F^n$, $n=1,\dots,N$ denote
$\Gamma^{\mu_n}$-parallel light cone lifts of the Darboux transforms
$\hat f^n$.
Note that the non-degeneracy assumption for the edge circles $c^n_{ij}$
ensures that the $\hat F^n_i$ or $\hat F^n_j$ are linearly independent
$\bmod f_i\oplus f_j$.

Now let $Q\in\R^{4,1}\setminus\{0\}$ denote the vector defining the space form.
Then
$$
   Q = \sum_{n=1}^N\alpha_n\hat F^n
$$
with suitable functions $\alpha_n:M\to\R$.
First, we see that the $\alpha_n$ are constant:
$$
   0 = dQ_{ij}
   = \sum_{n=1}^Nd(\alpha_n)_{ij}\hat F^n_i+(\alpha_n)_jd\hat F^n_{ij}
	= \sum_{n=1}^Nd(\alpha_n)_{ij}\hat F^n_i \bmod f_i\oplus f_j,
$$
so that $d(\alpha_n)_{ij}=0$ since $(\hat F^1,\dots,\hat F^N)$ are linearly
independent $\bmod f_i\oplus f_j$.
Now,
$$
	0 = dQ_{ij}
	= \sum_{n=1}^N\alpha_nd\hat F^n_{ij}
	= \sum_{n=1}^N\alpha_n\mu_n\{
		\langle\hat F^n,F\rangle_j\,F_i - \langle\hat F^n,F\rangle_i\,F_j
		\},
$$
showing that
$$
	Z := \sum_{n=1}^N\alpha_n\mu_n\hat F^n \perp f
$$
defines a sphere congruence enveloped by $f$:
note that $|Z|^2>0$ since the $\hat F^n$ and $F$ are linearly independent.
Hence $f$ has a linear or quadratic conserved quantity if $N=2$ or $N=3$,
respectively, by \href[thm.pcqReconstruction]{Lemma~4.9}.
\endproof

Note that there are no corresponding theorems for the existence of higher
degree polynomial conserved quantities for codimension $1$ isothermic nets:
the non-degeneracy assumption on the $\hat F^n_i$'s, $F_i$ and $F_j$ being
linearly independent in $\R^{4,1}$ restricts $N$ to numbers not greater
than $3$.

\proclaim\htag[thm.DarbouxPlanes]{Thm.~4.13}.
Let $\hat f^n$, $n=1,2$, be two Darboux transforms with different parameters
$\mu_n$ of a discrete isothermic net $f$ and let
$$
	M \ni i \mapsto \hat c^n_i :
	= \span\{F_i,Z_i,\hat F^n_i\}
$$
denote the orthogonal circle congruences of the Ribaucour pairs $(f,\hat f^n)$,
$n=1,2$, with respect to an enveloped sphere congruence\footnote{That is: $Z$
defines the normal direction at each vertex of the net so that the notion of
an ``orthogonal circle'' is well defined --- thus, only the contact element
$Z_i+\R\,F_i$ at each point $i$ is needed.} $Z$ and let
$$
	(ij) \mapsto c_{ij} :
	= \span\{F_i,Z_i,F_j\}
	= \span\{F_i,Z_j,F_j\}
$$
denote the orthogonal edge circles\footnote{Note that
$\span\{F_i,Z_i,F_j\}=\span\{F_i,Z_j,F_j\}$ by the enveloping
condition \href[def.Envelope]{Def.~3.10}.}.
Further, let $Q\in\R^{4,1}$ define a quadric of constant curvature
so that neither the circles $c_{ij}$ nor the circles $\hat c_i$ are
straight lines and let
$$
	e_{ij} :
	= \span\{F_i,Z_i,F_j,Q\}
	= \span\{F_i,Z_j,F_j,Q\}
	\quad{\sl and}\quad
	\hat e^n_i :
	= \span\{F_i,Z_i,\hat F_i,Q\}
$$
denote the corresponding circle planes.
Assume that $e_{ij}\neq\hat e^n_i$ for every edge and
suppose that the circle planes $\hat e^1=\hat e^2$.
Then $f$ has a quadratic conserved quantity.

\proclaim\htag[thm.DarbouxCircles]{Thm.~4.14}.
Let $\hat f^n$, $n=1,2$, be two Darboux transforms with different parameters
$\mu_n$ of a discrete isothermic net $f$ and let
$$
   i \mapsto \hat c^n_i
   \quad{\sl and}\quad
   (ij) \mapsto c_{ij}
$$
denote the orthogonal circle congruences as in the
\href[thm.DarbouxPlanes]{previous theorem}.
Assume that $c_{ij}\neq\hat c^n_i$ for every edge and
suppose that $\hat c^1=\hat c^2$.
Then $f$ has a normalized linear conserved quantity.

The first of these two theorems is a discrete version of a famous
classical theorem,
cf.~\href[ref.da99]{[13]},
\href[ref.ei23]{[14,~\S84]} and
\href[ref.sa07]{[21, Thms.~2.1 \& 2.33]}.

\proof
Again, we prove both theorems as the proofs are very similar.

Let $\hat F^n$ denote $\Gamma^{\mu_n}$-parallel lifts of the two Darboux
transforms $\hat f^n$ and note that $Q$ is $\Gamma^0$-parallel.

First consider the second theorem:
as $\hat c^1=\hat c^2$ we have
$$
   \hat F^2 = \alpha_1\hat F^1 + \alpha\,F + \beta\,Z
$$
with suitable functions $\alpha_1,\alpha,\beta$,
where $\alpha_1\neq0$ since $\hat f^2\neq f$.
Hence
$$
   d(\alpha_1)_{ij}\hat F^1_i
   = d\hat F^2_{ij} - d(\alpha\,F+\beta\,Z)_{ij} - (\alpha_1)_jd\hat F^1_{ij}
   = 0 \bmod c_{ij},
$$
so that $\alpha_1$ is constant since $\hat c^1_i\neq c_{ij}$.
Moreover, as $\hat F^1$, $\hat F^2$ and $F$ are linearly independent,
$$
   \hat F^2 - \alpha_1\hat F^1 = \alpha\,F + \beta\,Z
$$
defines an enveloped sphere congruence and the claim follows
from \href[thm.pcqReconstruction]{Lemma~4.9}.

In the situation of the \href[thm.DarbouxPlanes]{first theorem} we have
$$
   \hat F^2 = \alpha_0Q + \alpha_1\hat F^1 + \alpha\,F + \beta\,Z
$$
and deduce
$$
   d(\alpha_1)_{ij}\hat F^1_i
   = d\hat F^2_{ij} - d(\alpha\,F+\beta\,Z)_{ij}
   - d(\alpha_0)_{ij}Q - (\alpha_1)_jd\hat F^1_{ij}
   = 0 \bmod e_{ij},
$$
so that, again, $\alpha_1$ is constant since $\hat e^1_i\neq e_{ij}$.
Then
$$
   d(\alpha_0)_{ij}Q
   = d\hat F^2_{ij} - d(\alpha\,F+\beta\,Z)_{ij} - (\alpha_1)_jd\hat F^1_{ij}
   = 0 \bmod c_{ij},
$$
showing that $\alpha_0$ is constant as well, because $c_{ij}$ was assumed
to be a proper circle.
Now
$$
   \hat F^2 - \alpha_1\hat F^1 - \alpha_0 Q
   = \alpha\,F + \beta\,Z,
$$
so that, again, the claim follows
from \href[thm.pcqReconstruction]{Lemma~4.9}.
\endproof

Note that, in the first theorem, we obtain a normalized quadratic conserved
quantity as soon as we assume that the circles $\span\{\hat F^1,\hat F^2,F\}$
through corresponding points of $f$, $\hat f^1$ and $\hat f^2$ do not
become straight lines in the quadric given by $Q$.

\h1 Discrete cmc nets in space forms.
We are now prepared to define discrete cmc nets in space forms.
A smooth isothermic surface $f$ in the conformal $3$-sphere has constant
mean curvature $H$ in a quadric of constant curvature $\kappa=-|Q|^2$,
given by $Q\in\R^{4,1}$, if and only if it has a linear conserved quantity
$$
   \pcq(\lambda) = \lambda Z+Q,
$$
where $Z$ is the mean curvature sphere congruence\footnote{The mean curvature
sphere congruence of $f$, consisting of spheres touching $f$ that have the same
mean curvature as the surface at the touching points, can be defined using
{\it any\/} ambient space form geometry: it can be characterized as
the {\it conformal Gauss map\/} of $f$, i.e.,
   the unique enveloped sphere congruence that induces the same conformal
   structure as $f$,
or as
the {\it central sphere congruence\/}, i.e.,
   the congruence of spheres that exchange the curvature spheres
   (via inversion) or, equivalently, that have second order contact
   with the surface in orthogonal directions.}
of $f$, see \href[ref.buca07]{[11]} or \href[ref.sa07]{[21, Thm.~2.27]}:
for smooth isothermic surfaces in the conformal $3$-sphere (of codimension $1$)
it turns out that the top coefficient of a polynomial conserved quantity is
necessarily its conformal Gauss map\footnote{Hence it follows directly that
an isothermic surface with linear conserved quantity has constant mean
curvature $H=-\langle Z,Q\rangle$ in the space form given by $Q$.}.
Note that, in the case of a spherical surface, the conformal Gauss map
$Z$ of $f$ is constant and $f$ has a constant conserved quantity $Q=Z$
(as in the discrete case: see \href[thm.type0]{Thm.~3.14}),
hence a linear conserved quantity, e.g.,
$$
   \pcq(\lambda) = \lambda Z + Q = (\lambda+1)\,Z.
$$

In the discrete setting we use this characterization as a definition:

\proclaim\htag[def.cmc]{Def.~5.1}.
A discrete isothermic net $f$ will be called a \define{discrete cmc net}
if it is special isothermic of type $N\leq1$.
In particular, if $f$ is special isothermic of type $1$ with normalized
linear conserved quantity
$$
   \pcq(\lambda) = \lambda Z + Q,
$$
we say that $Z$ is the \define{mean curvature sphere congruence} of $f$ 
in the quadric
$$
   {\cal Q} = \{Y\in L^4\,|\,\langle Y,Q\rangle=-1\}
$$
of constant curvature $\kappa=-|Q|^2$
and that $f$ has (constant) \define{mean curvature}%
\footnote{A change $a\to\tilde a={a\over c}$ of the cross ratio factorizing
function results in a change of the equation \reqn{PCQcondition} for a
linear conserved quantity, hence of the linear conserved quantity,
see \href[thm.pcqScale]{Lemma~3.6};
in particular, $\tilde Z=Z$ and $\tilde Q=cQ$.
Hence we obtain the effect of an \htag[key.AmbHom]{ambient homothety}:
$\tilde\kappa=c^2\kappa$ and $\tilde H=cH$.}
$$
   H := -\langle Z,Q\rangle.
$$

Thus discrete cmc nets in space forms are special isothermic nets and a
transformation theory is readily available to us, cf.~\href[ref.je00]{[17]}:
the B\"acklund transformation, see \href[def.backlund]{Def.~4.4}, yields a
transformation for discrete cmc nets in a given space form preserving the
mean curvature,
$$
   \tilde Q = Q
   \quad{\rm and}\quad
   \tilde H
   = -\langle\tilde Z,\tilde Q\rangle
   = -\langle Z,Q\rangle
   = H
$$
by \href[thm.pcqDarboux]{Thm.~4.5}, and
satisfying Bianchi permutability by \href[thm.BBperm]{Thm.~4.7};
the Calapso transformation provides a Lawson correspondence $f\mapsto f^\mu$
for discrete cmc nets, where both the mean and ambient curvature change,
$$
   \kappa^\mu = -|Q^\mu|^2 = \kappa + 2\mu H - \mu^2
   \quad{\rm and}\quad
   H^\mu = -\langle Z^\mu,Q^\mu\rangle = H - \mu
$$
since $Z^\mu=T^\mu Z$ and $Q^\mu=T^\mu(\mu Z+Q)$ by 
\href[thm.SpecialCalapso]{Thm.~3.13},
but
$$
   (H^\mu)^2+\kappa^\mu = \kappa + H^2
$$
remains invariant.

\proclaim\htag[def.lawson]{Def.~5.2}.
The Calapso transformation for discrete cmc nets will also be called
\define{Lawson correspondence}.

The main goal of this section will be to show that our definition
generalizes and truly extends previous definitions\footnote{Note that
the mean curvature in \href[ref.sch07]{[22]} is defined on the faces of
a principal net and therefore different from the mean curvature defined here,
living on the vertices.}
from
\href[ref.bopi96]{[3]},
\href[ref.jehopi99]{[16]},
\href[ref.bopi99]{[4]} and
\href[ref.je00]{[17]}.

\h2 Uniqueness and existence questions.
However, before addressing the relation of our definition with previous
approaches we shall discuss the construction and uniqueness of linear
conserved quantities for a given discrete isothermic net.

Clearly, given the value of a linear conserved quantity at one point of a
discrete cmc net $f$ with cross ratio factorizing function $a$, the linear
conserved quantity is uniquely determined as it is a parallel section of the
isothermic family of (flat) connections:
given an initial value of a linear (or polynomial) conserved quantity
$\pcq(\lambda)$ we can use \reqn{PCQcondition} to determine the values
of $\pcq(\lambda)$ at any point.
Thus, starting with a linear quantity $\lambda Z+Q$ at some point of
a discrete isothermic net, there is a unique parallel section of the
isothermic family of connections\footnote{Here we use the flatness of
the connections in the family.};
but the parallel section may fail to be linear or even to be polynomial
at other points of the isothermic net --- hence we may not obtain a linear
conserved quantity if we use the ``wrong'' initial value or if the isothermic
net is not cmc.

Using instead the equivalent condition \reqn{dPCQ}
   from \href[thm.PCQcondition]{Lemma~3.7}
(see also \href[thm.ConstantQ]{Cor.~3.8}
   and \href[thm.EnvelopingZ]{Cor.~3.9~(i-ii)}),
we learn that $\pcq(\lambda)=\lambda Z+Q$ is a linear conserved quantity for
an isothermic net $f$ with cross ratio factorizing function $a$ if and
only if
$$
   Q\equiv const,\quad
   Z\perp F\quad{\rm and}\quad
   Z_j = Z_i + {a_{ij}\over\langle F_i,F_j\rangle}
      \{ \langle Q,F_j\rangle\,F_i-\langle Q,F_i\rangle\,F_j \}
\eqno LCQconditions$$
for any edge $(ij)$ of the domain graph $M$.
Again, we obtain a propagation formula that fixes the linear conserved
quantity uniquely once an initial value of $\pcq(\lambda)$ is given at
some point of a discrete isothermic net;
however, now the existence of $Z$ as well as the incidence relation
$Z\perp F$ are only satisfied if the isothermic net was, in fact, cmc
and if the initial value of $\pcq(\lambda)$ was chosen correctly.

First we address the integrability of the difference equation for
the mean curvature sphere congruence $Z$,
$$
   Z_j = Z_i + {a_{ij}\over\langle F_i,F_j\rangle}
      \{ \langle Q,F_j\rangle\,F_i-\langle Q,F_i\rangle\,F_j \}.
\eqno PropagatingZ$$
This equation does clearly not depend on the choice of lift $F$
of the isothermic net $f$; to simplify the computation we may choose
a Moutard lift, satisfying \reqn{MoutardCond}, so that
$$\matrix{
   Z_k
   &=& Z_i
    +  \langle Q,F_k-F_i\rangle\,F_j
    -  \langle Q,F_j\rangle\,(F_k-F_i) \hfill\cr
   &=& Z_i + {a_{ij}-a_{il}\over\langle F_j,F_l\rangle}
      \{ \langle Q,F_j\rangle\,F_l-\langle Q,F_l\rangle\,F_j \} \hfill\cr
}\eqno Zdiagonal$$
on an elementary quadrilateral $(ijkl)$, where we used the Moutard
equation \reqn{MoutardEqn}.
The last expression is symmetric in $j$ and $l$;
hence the propagation equation \reqn{PropagatingZ} is integrable.

Next we address the incidence relation $Z\perp F$.
First consider an edge $(ij)$: from \reqn{PropagatingZ}
$$
   \langle Z_j,F_j\rangle
   = \langle Z_i,F_j\rangle + a_{ij}\langle Q,F_j\rangle.
$$
Consequently, incidence determines $Z_i$ at the center $i=i_{(0,0)}$ of
a non-spherical vertex star\footnote{Here we use the fact that any
vertex has four neighbours; in a more general quad-graph this argument
only works at vertices of degree $4$, i.e., not at ``umbilics'' of a
discrete net.}:
if we let $i_{(m,n)}$, $m,n\in\{-1,0,1\}$, denote the vertices of
a $3\times3$-grid then the equations
$$
   \langle Z_{i_{(0,0)}}, F_{i_{(0,0)}}\rangle
   = 0
   \quad{\rm and}\quad
   \langle Z_{i_{(0,0)}}, F_{i_{(m,n)}}\rangle
   = -a_{i_{(0,0)}i_{(m,n)}}\langle Q,F_{i_{(m,n)}}\rangle,
\eqno Zstar$$
where $m^2+n^2=1$, have a unique solution $Z_{i_{(0,0)}}$ since
the vertices $F_{i_{(0,0)}}$ and $F_{i_{(m,n)}}$ of a non-spherical
vertex star form a basis of $\R^{4,1}$.
Thus an appropriate choice of the initial value $Z_{i_{(0,0)}}$ for
the mean curvature sphere congruence at the center of a vertex star
ensures that incidence is satisfied on all five vertices of the vertex
star when using \reqn{PropagatingZ} to define $Z$ on the corresponding
$3\times3$-grid.
At the diagonal vertices $F_{i_{(m,n)}}$, $m,n=\pm1$, we obtain incidence
without further conditions:
let $(ijkl)$ denote an elementary quadrilateral; then, using \reqn{Zdiagonal},
\reqn{MoutardEqn} and \reqn{MoutardCond}, we get
$$\matrix{
   \langle Z_k,F_k\rangle
   &=& \langle 
   Z_i + {a_{ij}-a_{il}\over\langle F_j,F_l\rangle}
      \{ \langle Q,F_j\rangle\,F_l-\langle Q,F_l\rangle\,F_j \},
   F_i + {a_{ij}-a_{il}\over\langle F_j,F_l\rangle}
      \{ F_j-F_l \}
   \rangle \hfill\cr
   &=& \langle Z_i,F_i\rangle + {a_{ij}-a_{il}\over\langle F_j,F_l\rangle}\{
      \langle Z_i+a_{ij}Q,F_j\rangle
      - \langle Z_i+a_{il}Q,F_l\rangle
      \} \hfill\cr
   &=& \langle Z_i,F_i\rangle + {a_{ij}-a_{il}\over\langle F_j,F_l\rangle}\{
      \langle Z_j,F_j\rangle
      - \langle Z_l,F_l\rangle
      \} \hfill\cr
   &=& 0, \hfill\cr
}$$
that is, incidence of $Z_k$ and $f_k$.
Thus we have proved the following, cf.~\href[thm.pcqUniq]{Cor.~3.4}:

\proclaim\htag[thm.uniqLCQ]{Lemma 5.3}.
Let $f:\{(m,n)\,|\,m,n\in\{-1,0,1\}\}\to S^3$ be a non-spherical discrete
isothermic $3\times3$-net and let $Q\in\R^{4,1}\setminus\{0\}$.
Then $f$ has a unique linear conserved quantity $\pcq(\lambda)=\lambda Z+Q$.

Note that, if the vertex star of $f_{(0,0)}$ is cospherical, then the
corresponding $3\times3$-net is necessarily also cospherical since
$f$ is a discrete principal net, i.e., the vertices of its faces
are concircular.
Hence, assuming that $f$ is non-spherical in \href[thm.uniqLCQ]{Lemma~5.3},
we have that the vertex star used to define $Z$ at its center is
non-spherical.
Moreover, since any two adjacent vertex stars (i.e., vertex stars at the
endpoints of an edge) of a discrete isothermic net have two face circles
in common, they lie on the same sphere if they are cospherical
--- thus, an isothermic net is either spherical, hence type $0$, or it
has a non-spherical vertex star.

As a (degenerate) example consider
the Moutard lift $F$ of an isothermic $3\times3$-net, i.e.,
let $F$ satisfy $\langle F_i,F_j\rangle=a_{ij}$ as in \reqn{MoutardCond},
and choose $Q$ so that
$$
   Q\perp F_{(1,0)}-F_{(-1,0)},F_{(0,1)}-F_{(0,-1)},F_{(0,1)}-F_{(1,0)}
   \quad\Leftrightarrow\quad
   \langle Q,F_{(m,n)}\rangle=c_0
\eqno badQ$$
for $m^2+n^2=1$ and some $c_0\in\R$;
further let $c_1:=\langle Q,F_{(0,0)}\rangle$ and observe that,
from \reqn{MoutardEqn},
$$
   \langle F_{(m,n)}-F_{(0,0)},Q\rangle
   = {a_{(0,0)(m,0)}-a_{(0,0)(0,n)}\over\langle F_{(m,0)},F_{(0,n)}\rangle}
      \langle F_{(m,0)}-F_{(0,n)},Q\rangle
   = 0
$$
for $m,n=\pm1$, so that
$$
   \langle F_{(m,n)},Q\rangle = \cases{
      c_0 & if $m^2+n^2=1$, \cr
      c_1 & if $m^2-n^2=0$. \cr
      }
$$
Now, \reqn{Zstar} yields $Z_{(0,0)}+c_0F_{(0,0)}\perp F_{(0,0)},F_{(m,n)}$
for $m^2+n^2=1$, that is, $Z_{(0,0)}+c_0F_{(0,0)}$ is orthogonal to a basis
of $\R^{4,1}$, hence vanishes.
Consequently, \reqn{PropagatingZ} yields
$$
   Z_{(m,n)} = \cases{
      -c_0F_{(m,n)} & if $m^2-n^2=0$, \cr
      -c_1F_{(m,n)} & if $m^2+n^2=1$. \cr
   }
$$
Thus $Z\parallel F$ becomes isotropic and yields another Moutard lift of $f$.
Note that, if we choose $Q$ so that $c_0=0$, i.e., so that $Q$ represents the
(unique) sphere containing the outer points of the vertex star, then $Z$
vanishes for $m^2-n^2=0$, that is, we are in the situation that we
excluded from consideration in \href[thm.degPCQ]{Lemma~3.11}.
This last situation is a worst-case scenario for a non-spherical isothermic
net and that, by \reqn{PropagatingZ}, $Z$ does not vanish as soon as
$\langle F,Q\rangle$ does not, that is, as soon as $f$ does not
hit the infinity boundary of the space form given by $Q$.

Now recall from \href[thm.degPCQ]{Lemma~3.11} that, if $Z$ is null but does
not vanish, then it is necessarily a Moutard lift of $f$.
In particular,
$$
   Z_{(m,n)} = \cases{
      -c_0F_{(m,n)} & if $m^2-n^2=0$, \cr
      -c_1F_{(m,n)} & if $m^2+n^2=1$ \cr
   }
$$
for some $c_0,c_1\in\R$ and a Moutard lift $F$ satisfying \reqn{MoutardCond}
as any two Moutard lifts are related in this way.
Hence, for $m^2+n^2=1$,
$$
   \langle F_{(m,n)},Q\rangle
   = -{1\over c_1}\langle Z_{(m,n)},Q\rangle
   \equiv const,
$$
so that we are back in the situation of \reqn{badQ}.
Thus avoiding \reqn{badQ} we can ensure that $Z$ becomes spacelike and
hence a suitable rescaling of the constructed linear conserved quantity
will leave us with a {\it normalized\/} linear conserved quantity;
hence, according to \href[def.cmc]{Def.~5.1}:

\proclaim\htag[thm.cmcQ]{Lemma 5.4}.
Let $f:\{(m,n)\,|\,m,n\in\{-1,0,1\}\}\to S^3$ be a non-spherical discrete
isothermic $3\times3$-net and choose $Q\in\R^{4,1}$ to satisfy,
for $m,n\in\{-1,0,1\}$,
$$
   Q \not\perp F_{(m,n)}
   \quad{\sl and}\quad
   Q \not\in
      \{F_{(1,0)}-F_{(-1,0)},F_{(0,1)}-F_{(0,-1)},F_{(0,1)}-F_{(1,0)}\}^\perp,
\eqno goodQ$$
where $F$ is a Moutard lift of $f$.
Then $f$ is a discrete cmc net in a space form defined by a suitable
rescaling of $Q$.

Thus, as \reqn{goodQ} imposes only open conditions on $Q$, any given
isothermic $3\times3$-net is cmc in a $4$-parameter family of possible
space forms.
Enlarging the net and using \reqn{PropagatingZ} to propagate $Z$ will
add more incidence conditions, hence conditions on $Q$;
it is therefore natural to expect that a large enough generic isothermic net,
will have a unique linear conserved quantity and being cmc becomes a condition.
In particular, the conditions \reqn{LCQconditions} for a linear conserved
quantity on an isothermic $5\times5$-net
$$
   f:\{(m,n)\,|\,m,n\in\{-2,-1,0,1,2\}\} \to S^3
$$
can be reduced to the incidence conditions on an extended vertex star
$$
   \{f_{(0,0)}, f_{(\pm1,0)}, f_{(0,\pm1)}, f_{(\pm2,0)}, f_{(0,\pm2)}\}
$$
by the arguments that proved \href[thm.uniqLCQ]{Lemma 5.3}:
this yields nine linear equations for the linear conserved quantity
at the center vertex $f_{(0,0)}$ of the net --- hence we get existence
of a linear conserved quantity on any isothermic $5\times5$-net, and
we expect uniqueness up to scaling generically, i.e., using the
scaling freedom to normalize the obtained linear conserved
quantity we expect a generic isothermic $5\times5$-net to be
a discrete cmc net in a unique way.

However, the following example shows that even ``arbitrarily large''
isothermic nets can be cmc in different space forms --- and even with
their respective mean curvature sphere congruences defining the same
contact elements at each vertex.

For this purpose we reconsider our example \reqn{explMain}.
However, instead of thinking of $f$ as part of a ``zigzag-plane'' as in
\reqn{explMain}, we now think of it as part of a discrete circular cylinder
by letting
$$
   \alpha := \cos\varphi \quad{\rm and}\quad \beta := \sin\varphi
$$
for some $\varphi\in(0,{\pi\over2})$, so that\footnote{Remember that
$m,n\in\{-1,0,1\}$.}
$$
   f_{(m,n)}
   = (\eta m,{1+\cos\varphi\over2}+(-1)^n{1-\cos\varphi\over2},n\sin\varphi)
   = (\eta m,\cos n\varphi,\sin n\varphi).
\eqno explAltMain$$
Further we introduce\footnote{This is the (parallel) Christoffel transform
of $f$, which is the net of centers of the mean curvature sphere congruence
in Euclidean space.}
$$
   f^\ast_{(m,n)} := (\eta m,-\cos n\varphi,-\sin n\varphi).
$$
Now let
$$
   Z := {1\over2}F^\ast-Q
   \quad{\rm and}\quad
   Q := (1,0,0,0,-1),
\eqno explNormPCQ$$
where $F^\ast=({1+|f^\ast|^2\over2},f^\ast,{1-|f^\ast|^2\over2})$
denotes again the Euclidean lift of $f^\ast$, and note that
$$
   \langle F,F^\ast\rangle = -{1\over2}|f-f^\ast|^2 \equiv -2,
$$
so that $\langle\lambda Z+Q,F\rangle\equiv-1$.
Now observe that $|f^\ast|^2=|f|^2$ do not depend on $n$ and hence
$$
   dF^\ast_{(m,n)(m+1,n)}-dF_{(m,n)(m+1,n)}
   = dF^\ast_{(m,n)(m,n+1)}+dF_{(m,n)(m,n+1)}
   = 0;
$$
consequently, $Z$ and $Q$ define a linear conserved quantity
by \href[thm.PCQcondition]{Lemma~3.7}:
for all edges $(ij)$
$$
   d(\lambda Z+Q)_{ij} - {\lambda a_{ij}\over\langle F_i,F_j\rangle}
      \{\langle\lambda Z+Q,F\rangle_jF_i-\langle\lambda Z+Q,F\rangle_iF_j\}
   = 0,
$$
where $a_{ij}=\pm{1\over2}\langle F_i,F_j\rangle$ as in the example
\reqn{explMain}.
Since $|Z|^2\equiv1$ the linear conserved quantity is normalized,
characterizing $f$ as a discrete net of constant mean curvature
$$
   H = -\langle Z,Q\rangle = {1\over2}
$$
in Euclidean space (as $|Q|^2=0$) according to \href[def.cmc]{Def.~5.1},
as one would expect.

Superposition of the linear conserved quantities given by \reqn{explDegPCQ}
and by \reqn{explNormPCQ} then yields a $1$-parameter family of normalized
linear conserved quantities for $f$, given by
$$
   Z_\alpha := {1\over2}(F^\ast+\alpha(-1)^nF)-Q_0
   \quad{\rm and}\quad
   Q_\alpha := Q_0 - {2\alpha\over1-\cos\varphi}
      ({1+\cos\varphi\over2},0,1,0,-{1+\cos\varphi\over2}),
$$
where $Q_0=(1,0,0,0,-1)$ and $\alpha\in\R$.
Thus the isothermic net has constant mean curvatures
$$
   H_\alpha = -\langle Z_\alpha,Q_\alpha\rangle
   = {1\over2}(1+\alpha^2)-{1+\cos\varphi\over1-\cos\varphi}\alpha
$$
in the hyperbolic spaces
$$
   {\cal Q}_\alpha = \{Y\in L^4\,|\,\langle Y,Q_\alpha\rangle=-1\}
$$
of curvatures $\kappa_\alpha=-{4\alpha^2\over(1-\cos\varphi)^2}$.
Note that all $Z_\alpha$ define the same contact element at a vertex of $f$,
that is, they all define the same ``normal direction'' at a vertex.

Finally observe that, as for the discrete net \reqn{explMain} with degenerate
linear conserved quantity \reqn{explDegPCQ}, the restriction of the domain
to $\{-1,0,1\}^2$ is again not necessary:
the circular cylinder in \reqn{explAltMain} can be defined on all of
$\Z^2$ with linear conserved quantity given by \reqn{explNormPCQ}.
As the two nets, the ``zigzag-plane'' and the circular cylinder, coincide
on $\Z\times\{-1,0,1\}$ we obtain an example of an isothermic $N\times3$-net,
$N\in\N$ arbitrary, that is cmc in different space forms.

\h2 Relation with previous approaches.
At this point we are prepared to link the present approach to discrete
cmc nets to previous approaches; in particular, we shall discuss how
our definition relates to:
\item{$\bullet$} the notion of discrete minimal surface in Euclidean space
   introduced in \href[ref.bopi96]{[3,~Def.~7]};
\item{$\bullet$} the notion of discrete cmc net in Euclidean space
   introduced in \href[ref.jehopi99]{[16,~Sect.~5]};
\item{$\bullet$} the notion of discrete horospherical net (cmc $1$ net)
   in hyperbolic space in \href[ref.je00]{[17,~Def.~4.3]};
\item{$\bullet$} the definition of a net of constant mean curvature $H$
   in a space form of curvature $\kappa$, where $H^2+\kappa\geq0$,
   suggested in \href[ref.je00]{[17]}.

In this context we shall also discuss the relation of
our notion of ``mean curvature sphere'' with
   that of \href[ref.bopi99]{[4,~Sect.~4.5]}
and the notion of ``central sphere congruence''
   in \href[ref.bosu06]{[7,~Sect.~3]}.

First we wish to make contact with the definition of a discrete cmc net
in Euclidean space given in \href[ref.jehopi99]{[16]}:
recall that a discrete isothermic net $f:M^2\to\R^3$ is called cmc
in \href[ref.jehopi99]{[16]} if it has a parallel (isothermic) net,
i.e., a simultaneous Christoffel and Darboux transform $f^\ast$.
The (constant) distance of this parallel net yields (up to sign) the
mean curvature of both nets,
$$
   H := {1\over|f^\ast-f|}.
$$

Note that any constant multiple of $\langle df,df^\ast\rangle$
is a cross ratio factorizing function\footnote{Given a cross ratio
factorizing function $a$ of an isothermic net one defines the Christoffel
transform $f^\ast$ by
$$\textstyle
   df^\ast_{ij} = a_{ij}(df_{ij})^{-1} = -{a_{ij}\over|df_{ij}|^2}\,df_{ij};
$$
as $a$ is only defined up to constant multiples, so is $df^\ast$.
In the case of the parallel net of a discrete cmc net, however, there
is a canonical scaling for $f^\ast$.}, so that, without loss of generality,
$$
   a_{ij} = -{H\over2}\,\langle df_{ij},df^\ast_{ij}\rangle,
   \quad\hbox{\rm that is,}\quad
   df^\ast_{ij} = -{2a_{ij}\over H|df_{ij}|^2}\,df_{ij}.
$$
Further, since $|f^\ast-f|^2\equiv const$, \reqn{Leibniz} yields
$$
   (f^\ast-f)_{ij} \perp df_{ij}, df^\ast_{ij},
$$
that is, $H\,(f^\ast-f)$ defines a (unit) normal field for $f$ as well as for
$f^\ast$ in the sense of \href[ref.sch07]{[22]}, cf.~\href[ref.jehopi99]{[16]}.
Consequently,
$$
   d(|f^\ast|^2)_{ij}
   = 2\langle f^\ast_{ij},df^\ast_{ij}\rangle
   = -{4a_{ij}\over H|df_{ij}|^2}\langle f_{ij},df_{ij}\rangle
   = -{2a_{ij}\over H|df_{ij}|^2}\,d(|f|^2)_{ij}.
$$

Now suppose $f:M^2\to\R^3$ is a discrete cmc net in this sense,
with parallel cmc net $f^\ast$ and (constant) mean curvature $H$.
Let $F$ and $F^\ast$ denote the respective Euclidean lifts
and let
$$
   Z := HF^\ast - {1\over2H}Q
   \quad{\rm and}\quad
   Q = (1,0,0,0,-1)
$$
as in the above example of a discrete circular cylinder,
see \reqn{explNormPCQ}.
Then
$$
   \langle Z,F\rangle
   = -{H\over2}|f^\ast-f|^2 + {1\over2H}
   = 0
$$
and
$$
   dZ_{ij} - {a_{ij}\over\langle F_i,F_j\rangle}\,dF_{ij}
   = HdF^\ast_{ij} + {2a_{ij}\over|df_{ij}|^2}\,dF_{ij}
   = 0,
$$
so that $\lambda Z+Q$ defines a (normalized) linear conserved quantity
of $f$ by \href[thm.PCQcondition]{Lemma~3.7}.
Moreover,
$$
   \kappa = -|Q|^2 = 0
   \quad{\rm and}\quad
   -\langle Z,Q\rangle = H.
$$

Conversely, suppose that $f$ is a discrete isothermic net with cross ratio
factorizing function $a$ and normalized linear conserved quantity
$\lambda Z+Q$ so that
$$
   \kappa = -|Q|^2 = 0
   \quad{\rm and}\quad
   H = -\langle Z,Q\rangle \neq 0.
$$
Let
$$
   F^\ast := {1\over H}(Z+{1\over2H}Q) = {1\over2H^2}(2H\,Z+Q)
$$
denote the complementary net of $f$, see \href[def.complementary]{Def.~4.8}.
Clearly, $F^\ast$ defines a Darboux transform of $f$,
$$
   T^{2H}F^\ast \equiv const,
$$
and $\langle F^\ast,Q\rangle \equiv -1$.
Choosing a Euclidean lift $F$ of $f$, i.e., $\langle F,Q\rangle\equiv-1$,
\reqn{dPCQ} yields
$$
   dF^\ast_{ij} = {a_{ij}\over H\langle F_i,F_j\rangle}\,dF_{ij}.
\eqno EuclideanDual$$

Now let, without loss of generality, $Q=(1,0,0,0,-1)$ so that
$$
   F^\ast = ({1+|f^\ast|^2\over2},f^\ast,{1-|f^\ast|^2\over2})
   \quad{\rm and}\quad
   F = ({1+|f|^2\over2},f,{1-|f|^2\over2}).
$$
Then \reqn{EuclideanDual} gives
$$
   df^\ast_{ij} = -{2a_{ij}\over H|df_{ij}|^2}\,df_{ij}.
$$
Hence $f^\ast$ is also the Christoffel transform of $f:M^2\to\R^3$,
so that $f$ is a cmc net with parallel cmc net $f^\ast$
in the sense of \href[ref.jehopi99]{[16]}.

Thus we have proved the following

\proclaim\htag[thm.cmcEuclidean]{Thm.~5.5}.
A discrete isothermic net in $\R^3$ is cmc in Euclidean space in the
sense of \href[ref.jehopi99]{[16]} if and only if it is cmc in the
sense of \href[def.cmc]{Def.~5.1}.

In \href[ref.bopi99]{[4,~Sect.~4.5]}, a notion of mean curvature sphere
for a discrete isothermic net in Euclidean space is introduced.
We shall see that this mean curvature sphere is the same as our
mean curvature sphere $Z$ in the case of a discrete cmc net.

First recall \href[ref.bopi99]{[4,~Def.~12]}:
given a vertex star $\{f_{i_{(m,n)}}\in\R^3\,|\,m^2+n^2\leq1\}$ of a discrete
isothermic net with constant (negative) cross ratio function
$q_{ijkl}={a_{ij}\over a_{il}}$
there is a unique point $c_{i_{(0,0)}}$ so that
$$
   |f_{i_{(1,0)}}-c_{i_{(0,0)}}| = |f_{i_{(-1,0)}}-c_{i_{(0,0)}}|
   \quad{\rm and}\quad
   |f_{i_{(0,1)}}-c_{i_{(0,0)}}| = |f_{i_{(0,-1)}}-c_{i_{(0,0)}}|
\eqno bopiA$$
and
$$
     {|f_{i_{(1,0)}}-c_{i_{(0,0)}}|^2-|f_{i_{(0,0)}}-c_{i_{(0,0)}}|^2
      \over a_{i_{(0,0)}i_{(1,0)}}}
   = {|f_{i_{(0,1)}}-c_{i_{(0,0)}}|^2-|f_{i_{(0,0)}}-c_{i_{(0,0)}}|^2
      \over a_{i_{(0,0)}i_{(0,1)}}}.
\eqno bopiB$$
This point $c_{i_{(0,0)}}$ is the center of the mean curvature sphere at
$f_{i_{(0,0)}}$ and its radius
$$
   r_{i_{(0,0)}} := |f_{i_{(0,0)}}-c_{i_{(0,0)}}|.
\eqno bopiC$$

Now suppose that $f$ is a discrete cmc net in Euclidean space, with $Q$ as
above and mean curvature sphere congruence $Z$, and go back to \reqn{Zstar}:
let $F$ denote the Euclidean lift for $f$ as before and write
$$
   Z = {1\over r}({1+|c|^2-r^2},c,{1-|c|^2+r^2})
$$
in terms of its center $c$ and radius $r$;
then
$$
   0 = \langle Z_{i_{(0,0)}},F_{i_{(0,0)}}\rangle
   = -{1\over2r_{i_{(0,0)}}}
      (|f_{i_{(0,0)}}-c_{i_{(0,0)}}|^2-r_{i_{(0,0)}}^2)
$$
is equivalent to \reqn{bopiC}, while the remaining four equations of
\reqn{Zstar} read
$$
   1 = {\langle Z_{i_{(0,0)}},F_{i_{(m,n)}}\rangle
         \over a_{i_{(0,0)}i_{(m,n)}}}
   = -{1\over2r_{i_{(0,0)}}}
      {|f_{i_{(m,n)}}-c_{i_{(0,0)}}|^2-r_{i_{(0,0)}}^2
         \over a_{i_{(0,0)}i_{(m,n)}}},
$$
clearly implying \reqn{bopiB}; the equations \reqn{bopiA} follow
since the cross ratio function is constant, so that
$$
   a_{i_{(0,0)}i_{(1,0)}}=a_{i_{(0,0)}i_{(-1,0)}}
   \quad{\rm and}\quad
   a_{i_{(0,0)}i_{(0,1)}}=a_{i_{(0,0)}i_{(0,-1)}}.
$$
Hence we have proved:

\proclaim\htag[thm.mcsphere]{Thm.~5.6}.
The mean curvature sphere $Z$ of a discrete cmc net $f:M\to\R^3$ in Euclidean
space in the sense of \href[def.cmc]{Def.~5.1} is the mean curvature sphere of
$f$ in the sense of \href[ref.bopi99]{[4,~Sect.~4.5]}.

In fact, we have seen slightly more:
the equations \reqn{Zstar} define\footnote{Besides the sphere, the
equations \reqn{Zstar} also determine the scaling of its representative
in Minkowski space.} the mean curvature sphere of \href[ref.bopi99]{[4]}
for any isothermic net in Euclidean space with constant cross ratio
function\footnote{Recently, a new approach has come into focus,
where the ``mean curvature'' of a discrete principal net in $\R^3$
is defined, as a function on faces, via the area change of a face
when varying through parallel nets,
see \href[ref.sch07]{[22]} and \href[ref.boetal07]{[8]}.
This approach leads to the same class of discrete minimal or constant mean
curvature nets as the one discussed here \href[ref.bo07talk]{[9]}.}.

This suggests to use \reqn{Zstar} to define the \define{mean curvature sphere}
of an isothermic net in any space form.
Note that these mean curvature spheres of an isothermic net are, in contrast
to the smooth case, {\it not\/} M\"obius invariant: different choices of the
ambient space form given by $Q$ will, in general, lead to different mean
curvature spheres $Z$.
In particular, this also shows that the mean curvature sphere defined in this
way is generally\footnote{Using a Moutard lift satisfying \reqn{MoutardCond},
as for \reqn{badQ}, \reqn{MoutardEqn} and \reqn{Zstar} yield
$$
   \langle Z_{i_{(0,0)}},F_{i_{(m,n)}}\rangle = 0
   \quad\Leftrightarrow\quad
     a_{i_{(0,0)}i_{(m,0)}}\langle Q,F_{i_{(m,0)}}\rangle
   = a_{i_{(0,0)}i_{(0,n)}}\langle Q,F_{i_{(0,n)}}\rangle
$$
for $m,n=\pm1$.
Hence the ``mean curvature sphere'' defined by \reqn{Zstar} is the ``central
sphere'' of \href[ref.bosu06]{[7]} (cf.~\href[thm.din]{Lemma~2.2}) if and
only if, for $m^2+n^2=1$,
$$
   \langle Q, a_{i_{(0,0)}i_{(m,n)}}F_{i_{(m,n)}}\rangle = const.
$$}
different from the central sphere of \href[ref.bosu06]{[7]}, which is
M\"obius invariantly related to the isothermic net as it is characterized
by incidence only.

Following ideas from \href[ref.je00]{[17]}, the Lawson correspondence
of \href[def.lawson]{Def.~5.2} can be used to define discrete cmc nets
in space forms with
$$
   H^2 + \kappa > 0
$$
as Calapso transforms of cmc nets in Euclidean space.
Using permutability theorems the above characterization of cmc nets in
Euclidean space by the existence of a simultaneous Darboux and Christoffel
transform, i.e., the existence of a parallel cmc net can then be carried
over to other space forms to obtain an alternative characterization,
similarly to the way in which horospherical nets in hyperbolic space
can be characterized in two ways,
see \href[ref.je00]{[17,~Lemma~4.2]} or \href[ref.imdg]{[18,~\S5.7.37]}.

Namely, let $f^\ast$ denote the parallel cmc net of a discrete cmc net $f$
in Euclidean space and consider their Calapso transforms $f^\lambda$ and
$(f^\ast)^\lambda$
--- which are only determined up to M\"obius transformation.
Then $f^\lambda$ and $(f^\ast)^\lambda$ can be positioned in $S^3$ so that:
\item{(i)} they form a Darboux pair with parameter $-\lambda$ since
   $f$ and $f^\ast$ form a Christoffel pair,
   see \href[ref.je00]{[17,~Cor.~3.24]} or \href[ref.imdg]{[18,~\S5.7.34]};
\item{(ii)} they form a Darboux pair with parameter $\mu-\lambda$ since
   $f$ and $f^\ast$ form a Darboux pair with some parameter $\mu$,
   see \href[ref.je00]{[17,~Cor.~3.27]} or \href[ref.imdg]{[18,~\S5.7.35]}.
\par
That is, there are two ways to position $(f^\ast)^\lambda$ in $S^3$ so that it
is a Darboux transform of $f^\lambda$ with different parameters or, otherwise
said, $f^\lambda$ has a pair of M\"obius equivalent Darboux transforms.
More precisely, given a Calapso transform of a discrete cmc net in Euclidean
space, \reqn{cmcClifts} provides two antipodal Darboux transforms:
the antipodal map identifies the ambient space form.

\href[thm.DoubleDarboux]{Thm.~4.11} provides the precise formulation of the
converse:
normalizing the linear conserved quantity from the proof,
$$
   (1+\lambda\mu_1)\alpha_1\hat F^1 + (1+\lambda\mu_2)\alpha_2\hat F^2,
$$
we obtain $|\pcq(\lambda)|^2={(1+\lambda\mu_1)(1+\lambda\mu_2)\over\mu_1\mu_2}$
as its squared norm, hence
$$
   H^2 + \kappa
   = ({\mu_1+\mu_2\over2\mu_1\mu_2})^2 - {1\over\mu_1\mu_2}
   = {(\mu_1-\mu_2)^2\over4\mu_1^2\mu_2^2}
   > 0
$$
since $\mu_1\neq\mu_2$, so that the isothermic net $f$ is indeed a
Calapso transform of a cmc net in Euclidean space.

Finally we discuss minimal nets in Euclidean space and horospherical
nets in hyperbolic space, that is, the discrete constant mean curvature
nets with
$$
   H^2 + \kappa = 0.
$$
First recall that horospherical nets in hyperbolic space can be
defined as Darboux transforms of their hyperbolic Gauss maps
\href[ref.je00]{[17,~Def.~4.3]},
that is, as Darboux transforms of a spherical net:
hence \href[thm.cmcOne]{Lemma~4.10} identifies them as the discrete
cmc nets with
$$
   H^2 + \kappa = 0
   \quad{\rm and}\quad
   \kappa \neq 0
$$
in the sense of \href[def.cmc]{Def.~5.1}.
On the other hand, horospherical nets can equivalently be characterized
as Calapso transforms of a discrete minimal net in Euclidean space in the
sense of \href[ref.bopi96]{[3]} and, conversely, any discrete minimal net
in Euclidean space gives rise to a Lawson family of horospherical nets,
see \href[ref.je00]{[17,~Lemma 4.2]}.
Hence we can reverse the argument to conclude that the discrete minimal
nets in the sense of \href[ref.bopi96]{[3]} are those discrete cmc nets
with
$$
   H^2 + \kappa = 0
   \quad{\rm and}\quad
   \kappa = 0
$$
in the sense of \href[def.cmc]{Def.~5.1}.
Note that \href[ref.bopi96]{[3,~Thm.~8]} provides two equivalent
characterizations of discrete minimal nets in Euclidean space:
\item{(i)} as Christoffel transforms of spherical isothermic nets
   (their ``Gauss maps'') and
\item{(ii)} by the fact that their mean curvature sphere congruence
   consists of planes\footnote{Assuming, as in \href[ref.bopi96]{[3]},
   that $a_{ij}=\pm\delta$ a very similar computation as the one leading to
   \href[thm.mcsphere]{Thm.~5.6} shows that our statement here is equivalent
   to \href[ref.bopi96]{[3,~Def.~7]}:
   with the Euclidean lift $F$ and $Z=(d,n,-d)$, where $n$ is the unit
   normal of the plane described by $Z$ and $d=\langle n,f\rangle$ its
   distance from the origin, the equations \reqn{Zstar} become
   $$
      \langle n_{i_{(0,0)}},f_{i_{(m,n)}}-f_{i_{(0,0)}}\rangle
      = \langle Z_{i_{(0,0)}},F_{i_{(m,n)}}\rangle
      = a_{i_{(0,0)}i_{(m,n)}}
      = \pm\delta.
   $$}.

We summarize these discussions in the following

\proclaim\htag[thm.minimal]{Thm.~5.7}.
A discrete isothermic net is
\item{\rm(i)} minimal in $\R^3$ in the sense of \href[ref.bopi96]{[3]}
   iff it is cmc with $H=\kappa=0$ in the sense of \href[def.cmc]{Def.~5.1};
\item{\rm(ii)} horospherical in the sense of \href[ref.je00]{[17]}
   iff it is cmc with $-H^2=\kappa<0$ in the sense of \href[def.cmc]{Def.~5.1}.

\h2 Discrete cmc surfaces of revolution.
We conclude by discussing our example of a discrete surface of revolution
in more detail: here we will be interested in the construction of discrete
cmc surfaces of revolution with prescribed mean curvature in a given
ambient space form.

Thus we consider a discrete surface of revolution
$$
   (m,n) \mapsto F_{(m,n)} = (-1)^m(M_m+\Phi_nC)
   \in \R^{2,1}\oplus\R^2,
\eqno LCQrotnet$$
see \reqn{RevolutionNet}, with cross ratio factorizing function
$$
   a_{ij} = \alpha\,\langle F_i,F_j\rangle
\eqno LCQrotfact$$
since $F$ is a Moutard lift\footnote{We shall use the scaling freedom $\alpha$
in the cross ratio factorizing function later to normalize the linear conserved
quantity that we will construct, see \href[thm.pcqScale]{Lemma~3.6}.}.
By \href[thm.symmetricLCQ]{Cor.~3.16} a linear conserved quantity
$\pcq(\lambda)=\lambda Z+Q$ has the same rotational symmetry as the net if
and only if the vector $Q$ defining the ambient space form has this symmetry;
in particular, the linear conserved quantity of an equivariant cmc net is
rotationally symmetric with
$$
   \pcq_{(m,n)}(\lambda)
   = \pcq^\perp_m(\lambda) + (-1)^{m+1}\lambda p_m(\lambda)\Phi_nC
   = \lambda\{Z^\perp_m-\alpha\langle Q,M_m\rangle\,\Phi_nC\} + Q,
$$
where $Z^\perp,Q\in\R^{2,1}$,
see \reqn{symmetricPCQ} and \href[thm.pcqScale]{Lemma~3.6}.
Hence
$$\matrix{
   Z_{(m,n)}
   &=& Z^\perp_m + \alpha\langle Q,M_m\rangle\,M_m
   &-& \alpha\langle Q,M_m\rangle\,(M_m+\Phi_nC) \hfill\cr
   &=& S_m
   &-& \alpha\langle Q,F_{(m,n+1)}\rangle\,F_{(m,n)}, \hfill\cr
}\eqno LCQrotZ$$
where $S_m=S_{(m,n)(m,n+1)}$ is the curvature sphere \reqn{CurvatureSphere}
on an edge in the rotational direction: note that this family of curvature
spheres does only depend on $m$, so that our discrete surface of revolution
is the envelope of a $1$-parameter family of spheres in the sense of
\href[def.Envelope]{Def.~3.10}.

Our first aim is to formulate the condition \reqn{dPCQ} for a linear
conserved quantity of the given form:
clearly
$$
   0 = \langle Z,F\rangle = (-1)^m\langle S,M\rangle,
$$
that is, the incidence relation is again expressed by orthogonality;
further,
$$\matrix{
   0
   &=& dZ_{(m,n)(m+1,n)} + \alpha\{
      \langle Q,F_{(m,n)}\rangle\,F_{(m+1,n)}
    - \langle Q,F_{(m+1,n)}\rangle\,F_{(m,n)}\} \hfill\cr
   &=& dS_{m,m+1} - 2\alpha\langle Q,M_{m,m+1}\rangle\,dM_{m,m+1} \hfill\cr
}$$
and
$$\matrix{
   0
   &=& dZ_{(m,n)(m,n+1)} + \alpha\{
      \langle Q,F_{(m,n)}\rangle\,F_{(m,n+1)}
    - \langle Q,F_{(m,n+1)}\rangle\,F_{(m,n)}\} \hfill\cr
}$$
are identically satisfied.
Hence we obtain:

\proclaim\htag[thm.rotLCQcondition]{Lemma 5.8}.
Let
$Q\in\R^{2,1}$ and $Z_{(m,n)}:=S_m-\alpha\langle Q,F_{(m,n)}\rangle F_{(m,n)}$,
where $F$ is a Moutard lift \reqn{LCQrotnet} of a discrete surface of
revolution and $m\mapsto S_m\in\R^{2,1}$ is a discrete $1$-parameter
family of spheres. 
Then $\pcq(\lambda):=\lambda Z+Q$ defines a linear conserved quantity
for $f$ with respect to \reqn{LCQrotfact} as a cross ratio factorizing
function if and only if
$$
   0 \equiv \langle S,M\rangle \quad{\sl and}\quad
   dS = 2\alpha\langle Q,M\rangle\,dM.
\eqno LCQrotScondition$$

Note that the equations \reqn{LCQrotScondition} determine $S_m$ and $S_{m+1}$
(and, hence, $Z_m$ and $Z_{m+1}$) up to a common orthogonal offset from the
plane\footnote{Recall that $M$ in \reqn{RevolutionNet} takes values in
one component of the hyperbolic quadric $|Y|^2=-1$ in $\R^{2,1}$.}
spanned by $M_m$ and $M_{m+1}$.
Also prescribing the mean curvature $H$,
$$
   \langle S,Q\rangle
   = \langle Z,Q\rangle + \alpha\langle M,Q\rangle^2
   = -H + \alpha\langle M,Q\rangle^2,
\eqno LCQrotHcondition$$
we obtain six equations\footnote{Note that
$d(\langle M,Q\rangle^2)=2\langle M,Q\rangle\langle dM,Q\rangle$
by \reqn{Leibniz} so that \reqn{LCQrotHcondition} only adds one
real equation to \reqn{LCQrotScondition}.}
$$
\vcenter{\halign{$#$ \hfill& $#$ \hfill&\quad $#$ \hfill&$#$ \hfill\cr
   \langle S_m,Q\rangle &= -H+\alpha\langle M_m,Q\rangle^2 &
   \langle S_{m+1},Q\rangle &= -H+\alpha\langle M_{m+1},Q\rangle^2 \cr
   \langle S_m,M_m\rangle &= 0 &
   \langle S_{m+1},M_m\rangle &
      = -\alpha\langle Q,M_{m,m+1}\rangle\,|dM_{m,m+1}|^2 \cr
   \langle S_m,M_{m+1}\rangle &
      = -\alpha\langle Q,M_{m,m+1}\rangle\,|dM_{m,m+1}|^2 &
   \langle S_{m+1},M_{m+1}\rangle &= 0 \cr
}}
$$
that determine $S_m$ and $S_{m+1}$ as soon as $(Q,M_m,M_{m+1})$ is a basis
of $\R^{2,1}$.

In order to understand this restriction geometrically we analyze what
happens when $M_m$, $M_{m+1}$ and $Q$ are linearly dependent, that is,
when $Q$ is in the $(1,1)$-plane spanned by $M_m$ and $M_{m+1}$.
First recall that the light cone in $\R^{2,1}$ is the axis of our
discrete surface of revolution \reqn{RevolutionNet};
similarly, for fixed $n$,
$$
   \span\{F_{(m,n)},F_{(m+1,n)},\Phi_nC\}
   = \span\{M_m,M_{m+1},\Phi_nC\}
   =: c_{m,m+1}
$$
defines the circle through $F_{(m,n)}$ and $F_{(m+1,n)}$ that intersects
the axis orthogonally --- which becomes a geodesic in the quadric ${\cal Q}$
of constant curvature \reqn{SpaceForm} given by $Q$ as soon as
$$
   Q \in c_{m,m+1}.
$$
Thus prescribing an edge of a meridian curve, $M_m$ and $M_{m+1}$,
and a space form, $Q$, so that the orthogonal circle $c_{m,m+1}$ of
the axis passing through $F_{(m,n)}$ and $F_{(m+1,n)}$ does {\it not\/}
become a straight line, any choice of $H$ will lead to a unique solution
$S_m$ and $S_{m+1}$ of the equations \reqn{LCQrotScondition} and
\reqn{LCQrotHcondition}, hence to a unique linear conserved quantity
for the discrete net obtained by rotating the edge.

The top coefficient $Z$ of a linear conserved quantity constructed from
a solution $S$ of \reqn{LCQrotScondition} and \reqn{LCQrotHcondition} is,
with $S$, spacelike since $S\perp M$ and has constant length by
\href[thm.pcqNorm]{Lemma~3.5}.
The idea is then to use our scaling freedom $\alpha$ in the cross ratio
factorizing function \reqn{LCQrotfact} to obtain a {\it normalized\/}
linear conserved quantity, as sought in \href[def.cmc]{Def.~5.1}:
note that, in contrast to \href[thm.pcqScale]{Lemma~3.6}, where only the
conserved quantity condition played a role, the prescribed mean curvature
equation \reqn{LCQrotHcondition} causes a more complicated dependence of
$S$ on $\alpha$.
In particular\footnote{Note that the symmetry of the formula confirms
that $|S_m|^2=|S_{m+1}|^2$.},
$$
   |S|^2 = {C^2\over A}\,H^2 - {A\over\Delta}(\alpha+{B\over A}H)^2,
$$
where we let
$$\matrix{
   \Delta &:=& |Q\wedge M_m\wedge M_{m+1}|^2, \hfill\cr
   C &:=& |dM_{m,m+1}|^2\langle M_{m,m+1},Q\rangle, \hfill\cr
   B &:=& -|dM_{m,m+1}|^2\{
      |M_{m,m+1}|^2\langle M_m,Q\rangle\langle M_{m+1},Q\rangle
      + 2\langle M_{m,m+1},Q\rangle^2
      \}, \hfill\cr
   A &:=& -{B^2-\Delta C^2\over|M_m\wedge M_{m+1}|^2}. \hfill\cr
   }
$$
Note that $|Q\wedge M_m\wedge M_{m+1}|^2,|M_m\wedge M_{m+1}|^2<0$.
Also, the assumption $M_m,M_{m+1}\not\perp Q$ that the endpoints of
our meridian edge do not lie in the infinity boundary of the space
form defined by $Q$ implies that $B$ and $C$ do not simultaneously
vanish; as a consequence $A>0$.

Thus the equation $|S|^2=1$ can be solved for $\alpha$ if and only if
$C^2H^2\leq A$.

Moreover, as $\alpha$ is a factor of our cross ratio factorizing function,
we seek a non-zero solution of the equation:
since $A=C^2H^2$ clearly implies $H\neq0$ we need to exclude
$$
   A-C^2H^2 = B = 0
   \quad\Rightarrow\quad
   H^2 = {\Delta\over|M_m\wedge M_{m+1}|^2}.
$$
Thus the equation $|S|^2=1$ has a non-zero solution $\alpha$ if and only if
$$
   C^2H^2 \leq A
   \quad{\rm and}\quad
   C^2H^2 = A \Rightarrow H^2 \neq {\Delta\over|M_m\wedge M_{m+1}|^2}.
\eqno{LCQrotHconstraintA}$$
Note that, in the case of strict inequality in \reqn{LCQrotHconstraintA},
$C^2H^2<A$, our construction will, in general, provide two different linear
conserved quantities\footnote{As a consequence, our discrete net will have a
polynomial conserved quantity of degree $0$ by \href[thm.pcqUniq]{Cor.~3.4}
--- which is not too surprising since the net is clearly spherical,
cf.~\href[thm.type0]{Thm.~3.14}.}.
However, in the case $H=0$ of a minimal surface, the equation $|S|^2=1$
has always exactly one non-zero solution $\alpha$.

\proclaim\htag[thm.rotLCQedge]{Lemma 5.9}.
Prescribing an ambient space form and a mean curvature,
the equivariant discrete surface of revolution obtained by rotating
a single edge in the space form has a normalized linear conserved quantity
as soon as:
\item{\rm(i)} the orthogonal circle of the axis of revolution
   passing through the endpoints of the edge is not a straight
   line in the ambient constant curvature geometry, and
\item{\rm(ii)} the mean curvature is not chosen too large;
   more precisely, the mean curvature $H$ satisfies the
   constraint \reqn{LCQrotHconstraintA}.

Having equipped an initial edge of a meridian curve for a discrete
equivariant cmc net with prescribed mean curvature $H$ and ambient
space form $Q$ with a suitable mean curvature sphere at both
endpoints,
we shall now investigate how to propagate the meridian curve $M$
and its enveloped $1$-parameter family of spheres $S$ to ``build''
a larger equivariant cmc net.
That is, we aim to construct $M_{m+1}$ and $S_{m+1}$ from the data
at the other endpoint of the edge, $M_m$, $S_m$, $Q$.

As a discrete analogue of a constant speed parametrization in the hyperbolic
plane of the meridian curve\footnote{Recall that a hyperbolic constant speed
parametrization of the meridian curve of a (smooth) surface of revolution
leads to a conformal curvature line parametrization of the surface up to
constant rescaling of the parameters.}
we prescribe a constant cross ratio factorizing function $a$ along the
meridian curve so that
$$
   \langle M_{m+1},M_m\rangle
   = -(1+{a_{(m,n),(m+1,n)}\over\alpha})
   \equiv -(1+{c\over\alpha})
$$
for all $m$ and some $c\in\R$, see \reqn{LCQrotfact}.
Note that, since we seek $M_m$ and $M_{m+1}$ to take values in the same
hyperbolic plane,
$$
   {c\over\alpha} = {1\over2}|dM_{m,m+1}|^2 > 0.
$$
A second equation obtained from \reqn{LCQrotScondition},
$$
   0 = \langle M_{m+1},S_m\rangle
     + \alpha\langle Q,M_{m,m+1}\rangle|dM_{m,m+1}|^2
     = \langle M_{m+1},S_m+cQ\rangle
     + c\langle M_m,Q\rangle,
$$
then confines $M_{m+1}$ to a line in $\R^{2,1}$ as soon as $M_m$
and $S_m+cQ$ are linearly independent or, equivalently, as soon as
$$
   X_m :
   = S_m + c\,(Q+\langle Q,M_m\rangle\,M_m) \neq 0.
$$
Observe that, requiring the meridian curve to not cross the axis of
rotation into the other hyperbolic half plane nor to hit the infinity
boundary $Q^\perp$ of the ambient space form\footnote{We do, however,
not exclude the possibility of the meridian curve crossing the infinity
boundary.},
$$
   0 \neq c\,\langle Q,M_m\rangle\{\langle M_{m+1},M_m\rangle-1\}
   = \langle M_{m+1},X_m\rangle,
$$
so that $X_m\neq0$ describes a sphere intersecting the axis of rotation
orthogonally and containing the point of the meridian curve given by
$M_m$ since $X_m\perp M_m$.
In particular, $X\neq0$ at both endpoints of a ``proper'' meridian curve
edge.

Thus
$$
   M_{m+1} = \{
   {\alpha+c\over\alpha}\,M_m
   - {((\alpha+c)^2-\alpha^2)\langle Q,M_m\rangle\over\alpha}
     {X_m\over|X_m|^2} \} + t\,{Y_m\over|Y_m|^2}
$$
for a suitable $t\in\R$, where $Y_m$ is orthogonal to $X_m$
in the (Euclidean) plane $M_m^\perp$ in $\R^{2,1}$ and has the
same length,
$$
   Y_m \perp M_m,X_m
   \quad{\rm and}\quad
   |Y_m|^2 = |X_m|^2.
$$
As we wish $M$ to take values in the hyperbolic plane
$$\textstyle
   0 \buildrel!\over= 1 + |M_{m+1}|^2 
   = {1\over|X_m|^2}\{
      t^2 - {((\alpha+c)^2-\alpha^2)(1-2cH-c^2\kappa)\over\alpha^2}
      \},
$$
where we have used \reqn{LCQrotHcondition},
so that we obtain two candidates for $M_{m+1}$ as soon as
$$
   1 - 2cH - c^2\kappa > 0.
\eqno LCQrotHconstraintB$$
When propagating the meridian curve with constant cross ratio factorizing
function, then one of the two solutions must give the predecessor $M_{m-1}$
of $M_m$;
hence the propagation of a meridian curve is unique if it is possible.
Observe that the condition \reqn{LCQrotHconstraintB} does not depend on $m$;
hence it is automatically satisfied as soon as a ``seed'' meridian curve
contains more than one edge\footnote{In case we have {\it one\/} edge
of a seed meridian curve we infer that the quadratic equation has
at least one solution, i.e., $$1-2cH-c^2\kappa\geq0.$$
In the case $1-2cH-c^2\kappa=0$ of one solution, the described
propagation procedure will then alternate the two endpoints of one edge.}

Once $M_{m+1}$ is constructed \reqn{LCQrotScondition} yields
$$
   S_{m+1} := S_m + 2\alpha\langle Q,M_{m,m+1}\rangle\,dM_{m,m+1};
$$
it is then straightforward to verify incidence\footnote{Recall that the
equations \reqn{LCQrotScondition} also ensure that $|S|^2$ is constant.}
$$
   \langle S_{m+1},M_{m+1}\rangle
   = \langle M_{m+1},S_m+cQ\rangle + c\langle Q,M_m\rangle
   = 0
$$
and the mean curvature \reqn{LCQrotHcondition} being constant,
$$
   \langle Q,S_{m+1}\rangle - \alpha\langle Q,M_{m+1}\rangle
   = \langle Q,S_m\rangle - \alpha\langle Q,M_m\rangle.
$$
Hence, by \href[thm.rotLCQcondition]{Lemma~5.8}, we have succeeded in
propagating the meridian curve of a discrete cmc net of revolution:

\proclaim\htag[thm.rotLCQvertex]{Lemma~5.10}.
Let a point and unit normal of the meridian curve of a discrete
cmc net of revolution in a space form be given.
Prescribing a (constant) cross ratio factorizing function $c$ and
a Moutard factor $\alpha$, the meridian curve can be propagated
uniquely in either of two directions to obtain a discrete cmc
surface of revolution as long as:
\item{\rm(i)} ${c\over\alpha}>0$ and the point of the meridian curve
   does not lie in the infinity boundary of the space form, and
\item{\rm(ii)} $1-2cH-c^2\kappa>0$, where the mean curvature $H$ is
   given by \reqn{LCQrotHcondition} and $\kappa$ is the ambient curvature.

\vskip 3ex\hbox to \hsize{\hfil\vbox{\hsize=250pt
\hbox to \hsize{\psfile(250x246+145-35@79){s3tor}}\vglue 2ex
\hbox to \hsize{\hfil{\small
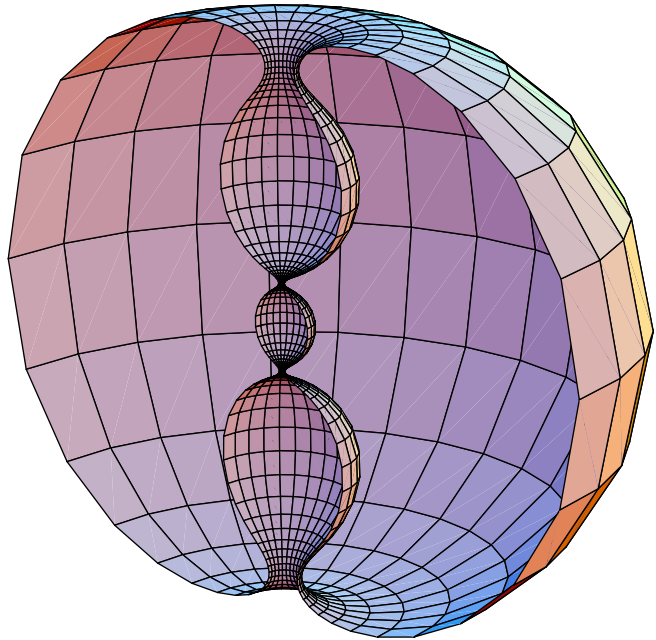
   \htag[fig.s3torus]{{\bf Fig.~2.}}\enspace Discrete cmc torus in $S^3$
   }\hfil}}\hfil}\vskip 2ex

These two lemmas now provide a method of construction for discrete cmc
nets of revolution in a prescribed ambient space form and with prescribed
mean curvature:

\proclaim\htag[thm.rotLCQconstruction]{Construction}.
Choose an ambient space form $Q\in\R^{2,1}$ and a mean curvature $H\in\R$.
\item[0pt]{} Choose an initial edge $M_0,M_1\in H^2\subset\R^{2,1}$
   so that:
   \item{\rm(i)} $\langle Q,M_0\rangle,\langle Q,M_1\rangle\neq0$, that is, 
      $M_0$ and $M_1$ do not lie in the infinity boundary of the space
      form\footnote{Note that this condition is sufficient but not
      necessary, as \href[fig.babo]{Fig.~1} suggests.} $Q$;
   \item{\rm(ii)} $(Q,M_0,M_1)$ is a basis of $\R^{2,1}$, that is,
      the straight line in $H^2$ through $M_0$ and $M_1$ is not
      straight in the chosen space form $Q$; and
   \item{\rm(iii)} the constraint \reqn{LCQrotHconstraintA} is satisfied
      with the chosen mean and ambient curvatures.
\item[0pt]{} Next, construct the spheres $S_0,S_1\in S^2\subset\R^{2,1}$
   to satisfy the conserved quantity and mean curvature conditions
   \reqn{LCQrotScondition} and \reqn{LCQrotHcondition},
   see \href[thm.rotLCQedge]{Lemma~5.9};
   generically, there is a choice of two such sets of spheres.
\item[0pt]{} Then, propagate this ``seed'' meridian curve using a
   constant cross ratio factorizing function,
   see \href[thm.rotLCQvertex]{Lemma~5.10};
   this yields a unique and proper propagation of the meridian curve
   if \reqn{LCQrotHconstraintB} is satisfied.

A variant of this construction was used to obtain the discrete cmc
torus in $S^3$ shown in \href[fig.s3torus]{Fig.~2} as well as to
construct the discrete minimal net in two copies of $H^3$ shown
in \href[fig.babo]{Fig.~1}: the plane shown in the figure indicates
the common infinity boundary of the two copies of the ambient $H^3$.

\vfill\eject
\h11 References.
{\frenchspacing

\item{\htag[ref.babo93]{1.}}
   M. Babich, A. Bobenko:
   {\it Willmore tori with umbilic lines and minimal surfaces in
   hyperbolic space\/};
   Duke Math. J. 72, 151--185 (1993)


\item{\htag[ref.bi04]{2.}}
   L. Bianchi:
   {\it Ricerche sulle superficie isoterme e sulla deformazione
   delle quadriche\/};
   Ann. Mat. 11, 93--157 (1904)

\item{\htag[ref.bopi96]{3.}}
   A. Bobenko, U. Pinkall:
   {\it Discrete isothermic surfaces\/};
   J. Reine Angew. Math. 475, 187--208 (1996)

\item{\htag[ref.bopi99]{4.}}
   A. Bobenko, U. Pinkall:
   {\it Discretization of surfaces and integrable systems\/};
   Oxf. Lect. Ser. Math. Appl. 16, 3--58 (1999)

\item{\htag[ref.boje01]{5.}}
   A. Bobenko, U. Hertrich-Jeromin:
   {\it Orthogonal nets and Clifford algebras\/};
   T\^ohoku Math. Publ. 20, 7--22 (2001)

\item{\htag[ref.bosu06a]{6.}}
   A. Bobenko, Y. Suris:
   {\it On Organizing Principles of Discrete Differential Geometry.
   Geometry of spheres\/};
   EPrint math.DG/0608291 (2006)

\item{\htag[ref.bosu06]{7.}}
   A. Bobenko, Y. Suris:
   {\it Isothermic surfaces in sphere geometries as Moutard nets\/};
   EPrint math.DG/0610434 (2006)

\item{\htag[ref.boetal07]{8.}}
   A. Bobenko, Y. Liu, H. Pottmann, J. Wallner, W. Wang:
   {\it Geometry of multi-layer structures for architecture\/};
   ACM Trans. Graphics 26 (2007), SIGGRAPH 2007

\item{\htag[ref.bo07talk]{9.}}
   A. Bobenko:
   {\it Discrete minimal surfaces from quadrilaterals\/};
   Talk at the International Congress on Minimal and Constant Mean Curvature
   Surfaces, Buzios, August 19--24, 2007 (2007)

\item{\htag[ref.buje06]{10.}}
   F. Burstall, U. Hertrich-Jeromin:
   {\it The Ribaucour transform in Lie sphere geometry\/};
   Differ. Geom. Appl. 24, 503--520 (2006)

\item{\htag[ref.buca07]{11.}}
   F. Burstall, D. Calderbank:
   {\it Conformal submanifold geometry\/};
   Manuscript 2007

\item{\htag[ref.ch67]{12.}}
   E. Christoffel:
   {\it Ueber einige allgemeine Eigenschaften der Minimumsfl\"achen\/};
   Crelle's J. 67, 218--228 (1867)

\item{\htag[ref.da99]{13.}}
   G. Darboux:
   {\it Sur les surfaces isothermiques\/};
   Ann. Sci. \'Ec. Norm. Sup. 16, 491--508 (1899)

\item{\htag[ref.ei23]{14.}}
   L.P. Eisenhart:
   {\it Transformations of surfaces\/};
   Princeton University Press, Princeton (1923)

\item{\htag[ref.jepe97]{15.}}
   U. Hertrich-Jeromin, F. Pedit:
   {\it Remarks on the Darboux transform of isothermic surfaces\/};
   Doc. Math. J. DMV 2, 313--333 (1997)

\item{\htag[ref.jehopi99]{16.}}
   U. Hertrich-Jeromin, T. Hoffmann, U. Pinkall:
   {\it A discrete version of the Darboux transform for isothermic surfaces\/};
   Oxf. Lect. Ser. Math. Appl. 16, 59--81 (1999)

\item{\htag[ref.je00]{17.}}
   U. Hertrich-Jeromin:
   {\it Transformations of discrete isothermic nets and discrete cmc-1
   surfaces in hyperbolic space\/};
   Manuscr. Math. 102, 465--486 (2000)

\item{\htag[ref.imdg]{18.}}
   U. Hertrich-Jeromin:
   {\it Introduction to M\"obius differential geometry\/};
   London Math. Soc. Lect. Note Ser. 300,
   Cambridge University Press (2003)

\item{\htag[ref.ho00]{19.}}
   T. Hoffmann:
   {\it Discrete curves and surfaces\/};
   PhD thesis, Technische Universit\"at Berlin (2000)

\item{\htag[ref.inko05]{20.}}
   J. Inoguchi, S. Kobayashi:
   {\it Characterizations of Bianchi-B\"acklund transformations
   of constant mean curvature surfaces\/};
   Int. J. Math. 16, 101--110 (2005)

\item{\htag[ref.sa07]{21.}}
   S. Santos:
   {\it Special isothermic surfaces\/};
   PhD thesis, University of Bath (2008)

\item{\htag[ref.sch07]{22.}}
   W. Schief:
   {\it On a maximum principle for minimal surfaces and their integrable
   discrete counterparts\/};
   Manuscript 2007

}

\vfill
\bgroup
\small
\def\\{\hfill\break}
\def\addwd{\hsize=.32\hsize}
\def\fran{\vbox{\addwd
   Francis E.\ Burstall\\
   Department of Mathematical Sciences\\
   University of Bath\\
   Bath, BA2 7AY (United Kingdom)\\
   Email: f.e.burstall@bath.ac.uk
   }}
\def\udo{\vbox{\addwd
   Udo Hertrich-Jeromin\\
   Department of Mathematical Sciences\\
   University of Bath\\
   Bath, BA2 7AY (United Kingdom)\\
   Email: u.hertrich-jeromin@bath.ac.uk
   }}
\def\wayne{\vbox{\addwd
   Wayne Rossman\\
   Department of Mathematics\\
   Kobe University\\
   Rokko, Kobe 657-8501 (Japan)\\
   Email: wayne@math.kobe-u.ac.jp
   }}
\def\susana{\vbox{\addwd
   Susana Santos\\
   Department of Mathematical Sciences\\
   University of Bath\\
   Bath, BA2 7AY (United Kingdom)\\
   Email: mapsdcs@bath.ac.uk
   }}
\hbox to \hsize{\hfil \fran \hfil \udo \hfil}\vskip 4ex
\hbox to \hsize{\hfil \wayne \hfil \susana \hfil}
\egroup
\bye